\documentclass[12pt]{amsart}
\usepackage{amsmath}
\usepackage{amsxtra}
\usepackage{amscd}
\usepackage{amsthm}
\usepackage{amsfonts}
\usepackage{amssymb}
\usepackage{eucal}
\usepackage{hyphenat}
\usepackage{verbatim}
\usepackage[all]{xy}
\usepackage{mathrsfs}
\usepackage{color}
\definecolor{deepjunglegreen}{rgb}{0.0, 0.29, 0.29}
\definecolor{darkspringgreen}{rgb}{0.09, 0.45, 0.27}
\definecolor{Red}{rgb}{0.7, 0,0}
\usepackage[normalem]{ulem}
\usepackage{stmaryrd}

\makeatletter
\usepackage{etex,etoolbox,needspace}
\pretocmd\section{\Needspace*{4\baselineskip}}{}{}
\makeatletter

\usepackage[pdftex,bookmarks=false,colorlinks=true,citecolor=darkspringgreen,debug=true,
  naturalnames=true,pdfnewwindow=true]{hyperref}

\newtheorem{thm}{Theorem}[subsection]

\newtheorem{cor}[thm]{Corollary}

\newtheorem{lem}[thm]{Lemma}
\newtheorem{prop}[thm]{Proposition}
\newtheorem{conj}[thm]{Conjecture}

\theoremstyle{definition}
\newtheorem{defn}[thm]{Definition}

\theoremstyle{remark}
\newtheorem{rem}[thm]{Remark}

\newcommand{\nc}{\newcommand}
\nc{\renc}{\renewcommand} \nc{\ssec}{\subsection}
\nc{\sssec}{\subsubsection}

\nc{\on}{\operatorname} \nc{\wh}{\widehat}
\nc\ol{\overline} \nc\ul{\underline} \nc\wt{\widetilde}

\newcommand{\red}[1]{{\color{Red}#1}}

\nc{\BA}{{\mathbb{A}}} \nc{\BC}{{\mathbb{C}}} \nc{\BF}{{\mathbb{F}}}
\nc{\BD}{{\mathbb{D}}} \nc{\BG}{{\mathbb{G}}} \nc{\BQ}{{\mathbb{Q}}}
\nc{\BM}{{\mathbb{M}}} \nc{\BN}{{\mathbb{N}}} \nc{\BO}{{\mathbb{O}}}
\nc{\BP}{{\mathbb{P}}} \nc{\BR}{{\mathbb{R}}}
\nc{\BZ}{{\mathbb{Z}}} \nc{\BS}{{\mathbb{S}}} \nc{\BW}{{\mathbb{W}}}\nc{\BY}{{\mathbb{Y}}}

\nc{\CA}{{\mathcal{A}}} \nc{\CB}{{\mathcal{B}}} \nc{\CalC}{{\mathcal{C}}} \nc{\CalD}{{\mathcal{D}}}
\nc{\CE}{{\mathcal{E}}} \nc{\CF}{{\mathcal{F}}}
\nc{\cF}{{\mathcal{F}}}
\nc{\CG}{{\mathcal{G}}} \nc{\CH}{{\mathcal{H}}}
\nc{\CI}{{\mathcal{I}}} \nc{\CK}{{\mathcal{K}}} \nc{\CL}{{\mathcal{L}}}
\nc{\CM}{{\mathcal{M}}} \nc{\CN}{{\mathcal{N}}}
\nc{\CO}{{\mathcal{O}}} \nc{\CP}{{\mathcal{P}}}
\nc{\cP}{{\mathcal{P}}}
\nc{\CQ}{{\mathcal{Q}}} \nc{\CR}{{\mathcal{R}}}
\nc{\CS}{{\mathcal{S}}} \nc{\CT}{{\mathcal{T}}}
\nc{\CU}{{\mathcal{U}}} \nc{\CV}{{\mathcal{V}}}  
\nc{\cV}{{\mathcal{V}}}  
\nc{\CY}{{\mathcal Y}}
\nc{\CW}{{\mathcal{W}}} \nc{\CZ}{{\mathcal{Z}}}

\nc{\scrM}{{\mathscr M}}
\nc{\cM}{{\check{\mathcal M}}{}} \nc{\csM}{{\check{\mathcal A}}{}}
\nc{\oM}{{\overset{\circ}{\mathcal M}}{}}
\nc{\obM}{{\overset{\circ}{\mathbf M}}{}}
\nc{\oCA}{{\overset{\circ}{\mathcal A}}{}}
\nc{\obA}{{\overset{\circ}{\mathbf A}}{}}
\nc{\ooM}{{\overset{\circ}{M}}{}}
\nc{\osM}{{\overset{\circ}{\mathsf M}}{}}
\nc{\vM}{{\overset{\bullet}{\mathcal M}}{}}
\nc{\nM}{{\underset{\bullet}{\mathcal M}}{}}
\nc{\oD}{{\overset{\circ}{\mathcal D}}{}}
\nc{\obD}{{\overset{\circ}{\mathbf D}}{}}
\nc{\oA}{{\overset{\circ}{\mathbb A}}{}}
\nc{\op}{{\overset{\bullet}{\mathbf p}}{}}
\nc{\cp}{{\overset{\circ}{\mathbf p}}{}}
\nc{\oU}{{\overset{\bullet}{\mathcal U}}{}}
\nc{\ofZ}{{\overset{\circ}{\mathfrak Z}}{}}
\nc{\oCY}{{\overset{\circ}{\CY}}{}}

\nc{\fa}{{\mathfrak{a}}} \nc{\fb}{{\mathfrak{b}}}
\nc{\fd}{{\mathfrak{d}}} \nc{\fe}{{\mathfrak{e}}} \nc{\ff}{{\mathfrak{f}}}
\nc{\fg}{{\mathfrak{g}}} \nc{\fgl}{{\mathfrak{gl}}}
\nc{\fh}{{\mathfrak{h}}} \nc{\fri}{{\mathfrak{i}}}
\nc{\fj}{{\mathfrak{j}}} \nc{\fk}{{\mathfrak{k}}} \nc{\fl}{{\mathfrak{l}}}
\nc{\fm}{{\mathfrak{m}}} \nc{\fn}{{\mathfrak{n}}}
\nc{\ft}{{\mathfrak{t}}} \nc{\fu}{{\mathfrak{u}}} \nc{\fv}{{\mathfrak{v}}}
\nc{\fw}{{\mathfrak{w}}} \nc{\fz}{{\mathfrak{z}}}
\nc{\fp}{{\mathfrak{p}}} \nc{\fq}{{\mathfrak{q}}} \nc{\frr}{{\mathfrak{r}}}
\nc{\fs}{{\mathfrak{s}}} \nc{\fsl}{{\mathfrak{sl}}}
\nc{\fso}{{\mathfrak{so}}} \nc{\fsp}{{\mathfrak{sp}}} \nc{\osp}{{\mathfrak{osp}}}
\nc{\hsl}{{\widehat{\mathfrak{sl}}}}
\nc{\hgl}{{\widehat{\mathfrak{gl}}}}
\nc{\hg}{{\widehat{\mathfrak{g}}}}
\nc{\chg}{{\widehat{\mathfrak{g}}}{}^\vee}
\nc{\hn}{{\widehat{\mathfrak{n}}}}
\nc{\chn}{{\widehat{\mathfrak{n}}}{}^\vee}

\nc{\fA}{{\mathfrak{A}}} \nc{\fB}{{\mathfrak{B}}} \nc{\fC}{{\mathfrak{C}}}
\nc{\fD}{{\mathfrak{D}}} \nc{\fE}{{\mathfrak{E}}}
\nc{\fF}{{\mathfrak{F}}} \nc{\fG}{{\mathfrak{G}}} \nc{\fH}{{\mathfrak{H}}}
\nc{\fI}{{\mathfrak{I}}} \nc{\fJ}{{\mathfrak{J}}}
\nc{\fK}{{\mathfrak{K}}} \nc{\fL}{{\mathfrak{L}}}
\nc{\fM}{{\mathfrak{M}}} \nc{\fN}{{\mathfrak{N}}}
\nc{\frP}{{\mathfrak{P}}} \nc{\fQ}{{\mathfrak{Q}}}
\nc{\fS}{{\mathfrak{S}}} \nc{\fT}{{\mathfrak{T}}} \nc{\fU}{{\mathfrak{U}}}
\nc{\fV}{{\mathfrak{V}}} \nc{\fW}{{\mathfrak{W}}}
\nc{\fX}{{\mathfrak{X}}} \nc{\fY}{{\mathfrak{Y}}}
\nc{\fZ}{{\mathfrak{Z}}}

\nc{\ba}{{\mathbf{a}}}
\nc{\bb}{{\mathbf{b}}} \nc{\bc}{{\mathbf{c}}} \nc{\be}{{\mathbf{e}}}
\nc{\bg}{{\mathbf{g}}} \nc{\bj}{{\mathbf{j}}} \nc{\bm}{{\mathbf{m}}}
\nc{\bn}{{\mathbf{n}}} \nc{\bp}{{\mathbf{p}}}
\nc{\bq}{{\mathbf{q}}} \nc{\br}{{\mathbf{r}}} \nc{\bt}{{\mathbf{t}}}
\nc{\bfu}{{\mathbf{u}}} \nc{\bv}{{\mathbf{v}}}
\nc{\bx}{{\mathbf{x}}} \nc{\by}{{\mathbf{y}}} \nc{\bz}{{\mathbf{z}}}
\nc{\bw}{{\mathbf{w}}} \nc{\bA}{{\mathbf{A}}}
\nc{\bB}{{\mathbf{B}}} \nc{\bC}{{\mathbf{C}}}
\nc{\bD}{{\mathbf{D}}} \nc{\bF}{{\mathbf{F}}} \nc{\bG}{{\mathbf{G}}}
\nc{\bH}{{\mathbf{H}}} \nc{\bI}{{\mathbf{I}}} \nc{\bJ}{{\mathbf{J}}}
\nc{\bK}{{\mathbf{K}}} \nc{\bM}{{\mathbf{M}}} \nc{\bN}{{\mathbf{N}}}
\nc{\bO}{{\mathbf{O}}} \nc{\bS}{{\mathbf{S}}} \nc{\bT}{{\mathbf{T}}}
\nc{\bU}{{\mathbf{U}}} \nc{\bV}{{\mathbf{V}}} \nc{\bW}{{\mathbf{W}}}
\nc{\bX}{{\mathbf{X}}}
\nc{\bY}{{\mathbf{Y}}} \nc{\bP}{{\mathbf{P}}}
\nc{\bZ}{{\mathbf{Z}}} \nc{\bh}{{\mathbf{h}}}

\nc{\sA}{{\mathsf{A}}} \nc{\sB}{{\mathsf{B}}}
\nc{\sC}{{\mathsf{C}}} \nc{\sD}{{\mathsf{D}}}
\nc{\sE}{{\mathsf{E}}} \nc{\sF}{{\mathsf{F}}} \nc{\sG}{{\mathsf{G}}} \nc{\sH}{{\mathsf{H}}}
\nc{\sI}{{\mathsf{I}}} \nc{\sK}{{\mathsf{K}}} \nc{\sL}{{\mathsf{L}}}
\nc{\sfm}{{\mathsf{m}}} \nc{\sM}{{\mathsf{M}}} \nc{\sN}{{\mathsf{N}}}
\nc{\sO}{{\mathsf{O}}} \nc{\sQ}{{\mathsf{Q}}} \nc{\sP}{{\mathsf{P}}}
\nc{\sT}{{\mathsf{T}}} \nc{\sZ}{{\mathsf{Z}}}
\nc{\sV}{{\mathsf{V}}} \nc{\sW}{{\mathsf{W}}}
\nc{\sfp}{{\mathsf{p}}} \nc{\sq}{{\mathsf{q}}} \nc{\sr}{{\mathsf{r}}}
\nc{\sfs}{{\mathsf{s}}} \nc{\st}{{\mathsf{t}}} \nc{\sfb}{{\mathsf{b}}}
\nc{\sfc}{{\mathsf{c}}} \nc{\sd}{{\mathsf{d}}}
\nc{\sz}{{\mathsf{z}}}
\nc{\sY}{{\mathsf{Y}}}

\nc{\tA}{{\widetilde{\mathbf{A}}}}
\nc{\tB}{{\widetilde{\mathcal{B}}}}
\nc{\tg}{{\widetilde{\mathfrak{g}}}} \nc{\tG}{{\widetilde{G}}}
\nc{\TM}{{\widetilde{\mathbb{M}}}{}}
\nc{\tO}{{\widetilde{\mathsf{O}}}{}}
\nc{\tU}{{\widetilde{\mathfrak{U}}}{}} \nc{\TZ}{{\tilde{Z}}}
\nc{\tx}{{\tilde{x}}} \nc{\tbv}{{\tilde{\bv}}}
\nc{\tfP}{{\widetilde{\mathfrak{P}}}{}} \nc{\tz}{{\tilde{\zeta}}}
\nc{\tmu}{{\tilde{\mu}}}

\nc{\urho}{\underline{\rho}} \nc{\uB}{\underline{B}}
\nc{\uC}{{\underline{\mathbb{C}}}} \nc{\ui}{\underline{i}}
\nc{\uj}{\underline{j}} \nc{\ofP}{{\overline{\mathfrak{P}}}}
\nc{\oB}{{\overline{\mathcal{B}}}}
\nc{\og}{{\overline{\mathfrak{g}}}} \nc{\oI}{{\overline{I}}}

\nc{\eps}{\varepsilon} \nc{\hrho}{{\hat{\rho}}} \nc{\balpha}{{\boldsymbol{\alpha}}}
\nc{\blambda}{{\boldsymbol{\lambda}}} \nc{\bmu}{{\boldsymbol{\mu}}} \nc{\bnu}{{\boldsymbol{\nu}}}
\nc{\btheta}{{\boldsymbol{\theta}}} \nc{\bzeta}{{\boldsymbol{\zeta}}} \nc{\bta}{{\boldsymbol{\eta}}}
\nc{\bbeta}{{\boldsymbol{\beta}}} \nc{\bkappa}{{\boldsymbol{\kappa}}} \nc{\bomega}{{\boldsymbol{\omega}}}

\nc{\one}{{\mathbf{1}}} \nc{\two}{{\mathbf{t}}}

\nc\DMO\DeclareMathOperator
\DMO\Sym{Sym}
\nc{\Tot}{{\mathop{\operatorname{\rm Tot}}}}
\nc{\Spec}{\mathop{\operatorname{\rm Spec}}}
\nc{\Ker}{\mathop{\operatorname{\rm Ker}}}
\nc{\Isom}{{\mathop{\operatorname{\rm Isom}}}}
\nc{\Hilb}{{\mathop{\operatorname{\rm Hilb}}}}
\nc{\deeq}{{\mathop{\operatorname{\rm deeq}}}}
\nc{\End}{{\mathop{\operatorname{\rm End}}}}
\nc{\Ran}{{\mathop{\operatorname{\rm Ran}}}}
\nc{\Ext}{{\mathop{\operatorname{\rm Ext}}}}
\nc{\Hom}{{\mathop{\operatorname{\rm Hom}}}}
\nc{\FS}{{\mathop{\operatorname{\rm FS}}}}
\nc{\CHom}{{\mathop{\operatorname{{\mathcal{H}}\it om}}}}
\nc{\GL}{{\mathop{\operatorname{\rm GL}}}}
\nc{\AGL}{{\mathop{\operatorname{\rm AGL}}}}
\nc{\PGL}{{\mathop{\operatorname{\rm PGL}}}}
\nc{\SL}{{\mathop{\operatorname{\rm SL}}}}
\nc{\SO}{{\mathop{\operatorname{\rm SO}}}}
\nc{\Sp}{{\mathop{\operatorname{\rm Sp}}}}
\nc{\PSp}{{\mathop{\operatorname{\rm PSp}}}}
\nc{\GSp}{{\mathop{\operatorname{\rm GSp}}}}
\nc{\GSpin}{{\mathop{\operatorname{\rm GSpin}}}}
\nc{\OSp}{{\mathop{\operatorname{\rm SOSp}}}}
\nc{\gr}{{\mathop{\operatorname{\rm gr}}}}
\nc{\Id}{{\mathop{\operatorname{\rm Id}}}}
\nc{\perf}{{\mathop{\operatorname{\rm perf}}}}
\nc{\defi}{{\mathop{\operatorname{\rm def}}}}
\nc{\length}{{\mathop{\operatorname{\rm length}}}}
\nc{\supp}{{\mathop{\operatorname{\rm supp}}}}
\nc{\oG}{{\mathop{\overline{G}}}}
\nc{\HC}{{\mathcal H}{\mathcal C}}

\nc{\pr}{{\operatorname{pr}}}
\nc{\Cliff}{{\mathsf{Cliff}}}
\nc{\loc}{{\operatorname{loc}}} \nc{\lc}{{\operatorname{lc}}}
\nc{\Fl}{{\mathbf{Fl}}} \nc{\Ffl}{{\mathcal{F}\ell}}
\nc{\Fib}{{\mathsf{Fib}}}
\nc{\Coh}{{\mathsf{Coh}}} \nc{\FCoh}{{\mathsf{FCoh}}}
\nc{\Perf}{{\mathsf{Perf}}}
\nc{\mPerf}{{\on{-}\mathsf{Perf}}}
\nc{\wtimes}{\mathbin{\widetilde\times}}

\nc{\reg}{{\text{\rm reg}}} \nc{\ren}{{\text{\rm ren}}}
\nc{\self}{{\text{\rm self}}}
\nc{\gvee}{{\mathfrak g}^{\!\scriptscriptstyle\vee}}
\nc{\tvee}{{\mathfrak t}^{\!\scriptscriptstyle\vee}}
\nc{\nvee}{{\mathfrak n}^{\!\scriptscriptstyle\vee}}
\nc{\bvee}{{\mathfrak b}^{\!\scriptscriptstyle\vee}}
       \nc{\rhovee}{\rho^{\!\scriptscriptstyle\vee}}

\nc{\cplus}{{\mathbf{C}_+}} \nc{\cminus}{{\mathbf{C}_-}}
\nc{\cthree}{{\mathbf{C}_*}} \nc{\Qbar}{{\bar{Q}}}

\newcommand{\oC}{\vphantom{j^{X^2}}\smash{\overset{\circ}{\vphantom{\rule{0pt}{0.55em}}\smash{C}}}}
\newcommand{\oCI}{\vphantom{j^{X^2}}\smash{\overset{\circ}{\vphantom{\rule{0pt}{0.55em}}\smash{\mathcal I}}}}
\newcommand{\oZ}{\vphantom{j^{X^2}}\smash{\overset{\circ}{\vphantom{\rule{0pt}{0.55em}}\smash{Z}}}}

\newcommand\iso{\mathbin{\vphantom{j^{X^2}}\smash{\overset{\sim}{\vphantom{\rule{0pt}{0.20em}}\smash{\longrightarrow}}}}}
\nc{\Gtimes}{\vphantom{j^{X^2}}\smash{\overset{G}{\vphantom{\rule{0pt}{0.30em}}\smash{\times}}}}
\nc{\sGtimes}{\vphantom{j^{X^2}}\smash{\overset{\mathsf G}{\vphantom{\rule{0pt}{0.30em}}\smash{\times}}}}

\nc{\svee}{{\!\scriptscriptstyle\vee}}

\nc{\bOmega}{{\overline{\Omega}}}

\nc{\seq}[1]{\stackrel{#1}{\sim}}

\nc{\aff}{{\operatorname{aff}}}
\nc{\fin}{{\operatorname{fin}}}
\nc{\mir}{{\operatorname{mir}}}
\nc{\triv}{{\operatorname{triv}}}
\nc{\ext}{{\operatorname{ext}}}
\nc{\righ}{{\operatorname{right}}}
\nc{\lef}{{\operatorname{left}}}
\nc{\forg}{{\operatorname{forg}}}
\nc{\fid}{{\operatorname{fd}}}
\nc{\odd}{{\operatorname{odd}}}
\nc{\even}{{\operatorname{even}}}
\nc{\modu}{{\operatorname{-mod}}}
\nc{\Gr}{{\operatorname{Gr}}}
\nc{\FT}{{\operatorname{FT}}}
\nc{\Mat}{{\operatorname{Mat}}}
\nc{\MSt}{{\operatorname{MSt}}}
\nc{\sph}{{\operatorname{sph}}}
\nc{\GR}{{\mathbf{Gr}}}
\nc{\Perv}{{\operatorname{Perv}}}
\nc{\Rep}{{\operatorname{Rep}}}
\nc{\Ind}{{\operatorname{Ind}}}
\nc{\IC}{{\operatorname{IC}}}
\nc{\Bun}{{\operatorname{Bun}}}
\nc{\Proj}{{\operatorname{Proj}}}
\nc{\Stab}{{\operatorname{Stab}}}
\nc{\pt}{{\operatorname{pt}}}
\nc{\bfmu}{{\boldsymbol{\mu}}}
\nc{\bfomega}{{\boldsymbol{\omega}}}
\nc{\calM}{\mathcal M}
\nc{\calA}{\mathcal A}
\nc{\calO}{\mathcal O}
\nc{\CC}{\mathcal C}
\nc{\cC}{\mathcal C}
\nc{\calN}{\mathcal N}
\nc{\grg}{\mathfrak g}
\nc{\dslash}{/\!\!/}
\nc{\tslash}{/\!\!/\!\!/}
\nc\grt{\mathfrak t}
\nc\bfM{\mathbf M}
\nc\bfN{\mathbf N}
\nc\Sig{\Sigma}
\nc\ZZ{\mathbb{Z}}
\nc\calC{\mathcal C}
\nc\calF{\mathcal F}
\nc\calX{\mathcal X}
\nc\calY{\mathcal Y}
\nc\QCoh{\operatorname{QCoh}}
\nc\Sch{\operatorname{Sch}}
\nc\IndCoh{\operatorname{IndCoh}}
\nc\Maps{\operatorname{Maps}}
\nc\Dmod{D-\operatorname{mod}}
\newcommand\Hecke{\operatorname{Hecke}}
\nc{\calD}{\mathcal D}
\nc\bfO{\mathbf O}

\nc\GG{\mathbb G}
\nc\calK{\mathcal K}
\nc{\calG}{\mathcal G}
\nc\RHom{\operatorname{RHom}}
\nc\Res{\operatorname{Res}}
\nc\Av{\operatorname{Av}}
\nc{\RH}{{\operatorname{RH}}}
\nc{\RT}{{\operatorname{RT}}}
\nc{\DR}{{\operatorname{DR}}}

\nc\grs{\mathfrak s}
\nc{\tilX}{\widetilde X}
\nc\calB{\mathcal B}
\nc\calS{\mathcal S}
\nc\calT{\mathcal T}
\nc\calZ{\mathcal Z}
\nc\LS{\operatorname{LocSys}}
\nc\bfL{\on{\mathbf L}}

\newcommand*\circled[1]
{*{#1}}
\makeatletter
\newcommand{\raisemath}[1]{\mathpalette{\raisem@th{#1}}}
\newcommand{\raisem@th}[3]{\raisebox{#1}{$#2#3$}}
\nc{\binlim}[2][]{\def\@tempa{#1}\@ifnextchar^{\@binlim{#2}}{\@binlim{#2}^{}}}
\def\@binlim#1^#2{\mathbin{\@ifempty{#2}{\mathop{#1}}{\mathop{#1}\@xp\displaylimits\@tempa^{#2}}}}

\makeatother

\nc\cX{{\mathcal X}}

\nc\Gm{{\mathbb G_m}}
\nc{\isoto}
    {\stackrel{\displaystyle\raisemath{-.2ex}{\sim}}\to}
\nc\cD\CalD
\nc\setm\setminus
\nc\cE\CE
\renc\Hecke{\mathit{\CH\kern-.2ex ecke}}
\nc\bsl\backslash
\nc\Fq{\mathbb F_q}
\nc\x\times
\nc\bGO{{\bG_\bO}}
\nc\bla\blambda
\nc\blt\bullet
\nc\opp{{\on{op}}}
\nc\srel\stackrel
\nc\tbx{\binlim{\widetilde\boxtimes{}}}
\nc\into\hookrightarrow
\nc\otto\leftrightarrow
\renc\setminus\smallsetminus
\nc\phitau{\varphi\tau}
\newenvironment{i-ii-iii}{%
\begin{enumerate}
}%
{\end{enumerate}}
\nc\ceil[1]{\lceil#1\rceil}  \nc\floor[1]{\lfloor#1\rfloor}
\nc\Lie{\on{Lie}}
\nc\sS{{\mathsf S}}
\nc\vvv{\ensuremath{\red\surd}}
\nc\la\lambda
\def\arxiv#1{\href{http://arxiv.org/abs/#1}{\tt arXiv:#1}} \let\arXiv\arxiv
\nc\kap{\kappa}
\nc\gra{\mathfrak a}
\nc\diag{\mathrm{diag}}
\nc\gl{\mathfrak{gl}}
\nc\sTr{\operatorname{sTr}}
\nc\hatG{\widehat{G}}
\nc\calL{\mathcal L}
\nc\Whit{\operatorname{Whit}}
\nc\KL{\operatorname{KL}}

\makeatletter
\renewcommand{\subsection}{\@startsection{subsection}{2}{0pt}{-3ex
plus -1ex minus -0.2ex}{-2mm plus -0pt minus
-2pt}{\normalfont\bfseries}} \makeatother

\numberwithin{equation}{subsection}
\allowdisplaybreaks
\nc\mto{\mapsto }
\nc\en{\enspace }

\begin{document}

\author[M.~Finkelberg]{Michael Finkelberg}
\address{Einstein Institute of Mathematics, The Hebrew University of Jerusalem,
  Edmond J. Safra Campus, Giv’at Ram, Jerusalem, 91904, Israel;
\newline  National Research University Higher School of Economics;
\newline Skolkovo Institute of Science and Technology}
\email{fnklberg@gmail.com}

\author[R.Travkin]{Roman Travkin}
\address{Skolkovo Institute of Science and Technology, Moscow, Russia}
\email{roman.travkin2012@gmail.com}

\author[R.Yang]{Ruotao Yang}
\address{Skolkovo Institute of Science and Technology, Moscow, Russia}
\email{yruotao@gmail.com}

\title
{Gaiotto equivalence for $\Rep_q(\on{F}(4))$}
\dedicatory{To Alexander Braverman on his 50th birthday}

\thanks{{\bf Mathematics Subject Classification (2010).}
Primary 22E67, 58B32; secondary 17B37.}

\thanks{{\bf Key words.} Affine Grassmannian, quantum supergroup.}

\begin{abstract}
This paper is a part of the series proving the Gaiotto conjecture for basic
classical quantum supergroups. The previous part~\cite{bft1,ty} proved the Gaiotto conjecture
for the general linear quantum supergroups $U_q(\fgl(N|M))$. Here we deal with the exceptional
quantum supergroup $U_q(\ff(4))$.
\end{abstract}

\maketitle

\tableofcontents

\section{Introduction}

\subsection{Gaiotto conjecture} This paper is a part of the series proving the Gaiotto conjecture
for basic classical quantum supergroups. The previous part~\cite{bft1,ty} proved the Gaiotto
conjecture for the general linear quantum supergroups $U_q(\fgl(N|M))$. Here we deal with the
exceptional quantum supergroup $U_q(\ff(4))$. Recall the conjecture from~\cite[\S3.3]{bft2}.

We consider a nilpotent element $e\in\fsp(6)$ of Jordan type $(3,3)$, so that $e$ lies in a
$14$-dimensional nilpotent orbit. We fix a maximal reductive subgroup $\PGL(2)$ in the centralizer
$Z_{\on{PSp}(6)}(e)$ of $e$. We choose an $\fsl_2$-triple $(e,h,f)$ in $\fsp(6)$.
The adjoint action of $h$ on $\fsp(6)$ equips it with a grading, and $\fsp(6)_{-1}=0$ ($e$ is an
even nilpotent element), while $\fu=\fsp(6)_{\leq-2}$ is the 7-dimensional nilpotent radical of
the parabolic subalgebra stabilizing an isotropic 2-plane in $\BC^6$. We denote by
$\on{U}\subset\on{PSp}(6)$ the unipotent subgroup with Lie algebra $\fu$. It is normalized by
$\PGL(2)\subset Z_{\on{PSp}(6)}(e)$. Finally, the Killing pairing with $e$ gives rise to a character
$\chi^0$ of $\fu$ and the same named character $\chi^0$ of $\on{U}$.

Extending the scalars to $\bF=\BC(\!(t)\!)$ we obtain a character $\chi^0\colon\on{U}(\bF)\to\bF$,
and we set $\chi:=\Res_{t=0}\chi^0\colon\on{U}(\bF)\to\BG_a$. Given
$\kappa\in\BC\setminus\overline\BQ$ we consider the dg-category
$\mathrm{D}\modu_\kappa(\Gr_{\on{PSp}(6)})$ of $D$-modules on the affine Grassmannian twisted
by the determinant line bundle to the power $\kappa$, and the category
$\mathrm{D}\modu_\kappa(\Gr_{\on{PSp}(6)})^{\PGL(2,\bO)\ltimes \on{U}(\bF),\chi,\lc}$ of locally compact
$(\PGL(2,\bO)\ltimes \on{U}(\bF),\chi)$-equivariant objects. (Recall that an equivariant
object is called locally compact if it becomes a compact object of $\mathrm{D}\modu^{\on{U}(\bF),\chi}_\kappa(\Gr_{\on{PSp}(6)})$ after forgetting the $\PGL(2,\bO)$ equivariant structure.)
Since we are going to relate our category
to representations of quantum {\em super}groups, all the categories of $D$-modules and constructible sheaves
we use in this paper will have super vector spaces as coefficients, as opposed to just vector spaces
(unless explicitly stated otherwise).

On the dual side, we consider the quantum group $U_q(\ff(4)),\ q=\exp(\pi\sqrt{-1}\kappa)$,
and we denote by $D^b(\Rep_q(\on{F}(4)))$ the dg-category of finite dimensional complexes of
$U_q(\ff(4))$-modules.

\begin{conj}
  \label{gaiotto f4}
   The categories
  \[\mathrm{D}\modu_\kappa(\Gr_{\on{PSp}(6)})^{\PGL(2,\bO)\ltimes \on{U}(\bF),\chi,\lc}\ \on{and}\
  D^b(\Rep_q(\on{F}(4)))\]
  are equivalent as braided monoidal categories, and this equivalence is compatible with the
  tautological $t$-structures.
\end{conj}

\subsection{Outline of the proof}
The strategy of our proof is the same as in~\cite{ty} (which follows a combination of the strategy of \cite{g} and \cite{bft1}), and we refer the reader to~\cite[\S1]{ty}
for the overview of the strategy. The only new ingredient is the proof of rigidity (dualizability)
of the Gaiotto monoidal category $\cC_{\kappa}^{\on{glob},\on{lc}}$ for transcendental level $\kappa$
in~\S\ref{Rigidity}. The proof is obtained by deformation from $\kappa=0$. For the general linear
quantum supergroups this part of the argument could be avoided and replaced by an embedding of
$\Rep_q(\GL(M|N))$ into the Gaiotto monoidal category due to the general results of~\cite{c}.

Four other cornerstones of the proof are:

\noindent (a) A choice of Borel subalgebra in $\ff(4)$ such that the highest weights of
finite dimensional irreducible representations of $U_q(\ff(4))$ are in natural bijection with
the set of relevant orbits in $\Gr_{\PSp(6)\times\PGL(2)}$, and the adjacency order on the set of
relevant orbits corresponds to the dominance order on the set of highest weights;

\noindent (b) dimension estimates for the intersections of certain ``semi-infinite'' orbits in
$\Gr_{\PSp(6)\times\PGL(2)}$;

\noindent (c) an isomorphism between two line bundles
on a certain zastava space; one pulled back from a colored configuration space of a curve $C$,
another one pulled back from $\Bun_{\PSp(6)\times\PGL(2)}(C)$;

\noindent (d) an equivalence of the local and global versions of the Gaiotto category along with
their fusion and factorization structures.

The bulk of the present paper is devoted
to the verification of (a,b,c). Unfortunately, this verification is absolutely independent of
similar verification for $\GL(M|N)$. In other words, our argument for various basic classical
quantum supergroups is case by case. This should be compared e.g.\ to the uniform proof of
Satake equivalence or Fundamental local equivalence for arbitrary reductive groups. Thus our
current pathetic level of understanding is similar to that of a bright student ignorant of the
structure theory of reductive Lie groups, discovering case by case numerous identities like
$\dim\fg=\on{rk}\fg\cdot(h+1)$.

\subsection{Acknowledgments}
The present note is a part of a joint project with Alexander Braverman.
We are honored to dedicate our note to his 50th birthday.

We are grateful to A.~Braverman, L.~Chen, I.~Entova-Aizenbud, Y.~Fu, D.~Gaitsgory, V.~Serganova
and H.~Yamane for the useful discussions.
The research of M.F.~was supported by the Israel Science Foundation (grant No.~994/24).

\section{Quantum supergroup $U_q(\ff(4))$ and factorizable sheaves}

\subsection{Quantum supergroup}
We fix a transcendental complex number $q\in\BC$.
For a definition of $U_q(\ff(4))$ see~\cite{y,chw}. Note that the definition
depends on a choice of a Borel subalgebra of $\ff(4)$, but the resulting quantum
algebras are all isomorphic according to~\cite[Proposition 7.4.1]{y2}.
We will use the Borel subalgebra with~3 odd isotropic simple roots, corresponding to the
simple root system number~4 in~\cite[\S18]{fss} or $\varSigma'''$ in~\cite[\S4.1]{m}.
More precisely, we fix a basis $\varepsilon_1,\varepsilon_2,\varepsilon_3,\delta$
of the dual Cartan subalgebra $\fh^*$ of $\ff(4)$ with the scalar product
$(\varepsilon_i,\varepsilon_j)=\delta_{ij},\ (\delta,\delta)=-3,\ (\varepsilon_i,\delta)=0$.
The positive simple roots of our Borel subalgebra of $\ff(4)$ are as follows:
\begin{equation}
  \label{b4}
  \alpha_1=\frac12(\varepsilon_1-\varepsilon_2-\varepsilon_3+\delta),\
  \alpha_2=\frac12(-\varepsilon_1+\varepsilon_2+\varepsilon_3+\delta),\
  \alpha_3=\frac12(\varepsilon_1-\varepsilon_2+\varepsilon_3-\delta),\
  \alpha_4=\varepsilon_2-\varepsilon_3.
\end{equation}

In this case the defining relations (quantum analogues of Serre relations) of $U_q(\ff(4))$
are explicitly written down in~\cite[Proposition 2.7(F4)]{chw}.

In the above basis
$\varepsilon_1,\varepsilon_2,\varepsilon_3,\delta$, the weight lattice 
$X$ is formed by all the collections
$\blambda=(\lambda_1,\lambda_2,\lambda_3,\lambda_4)\in(\frac12\BZ)^4$ of coordinates such that
$\lambda_1-\lambda_2\in\BZ\ni\lambda_2-\lambda_3$.
It contains the root sublattice $Q\subset X$ of index~2, formed by all the collections
$\blambda\in X$ such that $\lambda_3-\lambda_4\in\BZ$. In other words, for $\blambda\in Q$ the
coordinates are either all integers, or all half-integers. We also denote by $Q_{\on{pos}}$ (resp.\ $Q_{\on{neg}}$) the
monoid of nonnegative (resp.\ nonpositive) integral combinations of $\alpha_1,\alpha_2,\alpha_3,\alpha_4$.

We consider the category $\Rep_q(\on{F}(4))$ of finite dimensional representations of $U_q(\ff(4))$
equipped with an $X$-grading (compatible with the action of the Cartan part of the quantum supergroup).
According to~\cite[Theorem 48]{ge}, the highest weights of the irreducible finite dimensional
representations in $\Rep_q(\on{F}(4))$
(with respect to the above Borel subalgebra) are the same as the highest weights of the
irreducible finite dimensional representations of non-quantized supergroup $\on{F}(4)$.\footnote{Note
that if we do not require the existence of a compatible $X$-grading, there are finite dimensional
$U_q(\ff(4))$-modules whose limit as $q\to1$ does {\em not} exist, cf.~\cite[Remark 7.6]{ayy}
and~\cite[Remark 4.2 of the arXiv version]{y3}. We are
grateful to H.~Yamane for correcting our wrong definition of $\Rep_q(\on{F}(4))$ in an earlier version.}
Let us recall the classification of such highest weights. 
The highest weight of a finite
dimensional $\on{F}(4)$-module must be a dominant weight of the even part
$\on{F}(4)_{\bar0}=\on{Spin}(7)\times\SL(2)$, i.e.\ we must have
\begin{equation}
  \label{even dominant}
\lambda_1\geq\lambda_2\geq\lambda_3\geq0\leq\lambda_4,\ \lambda_1-\lambda_2\in\BZ\ni\lambda_2-\lambda_3.
\end{equation}

\begin{lem}
  \label{martir}
  A weight $\blambda$ is the highest weight of an irreducible finite dimensional representation
  of $\on{F}(4)$ (and of the quantum supergroup $U_q(\ff(4))$) if and only if
  the condition~\eqref{even dominant} along with the following
  auxiliary conditions are satisfied:
  \begin{equation}
    \label{odd dominant}
    \on{if}\ \lambda_4=0,\ \on{then}\ \lambda_1=\lambda_2+\lambda_3;\ \&\
    \on{if}\ \lambda_4>0,\ \on{then}\ \lambda_1>\lambda_2\geq\lambda_3>0.
    \end{equation}
\end{lem}

\begin{proof}
  The highest weights with respect to the Borel subalgebra with~3 even simple roots
  (the simple root system $\varSigma$ in~\cite[\S4.1]{m}) are classified
  in~\cite[Lemma~4.1]{m}. The positive root system $\varSigma'''$ is obtained from
  $\varSigma$ by composition of~3 odd reflections
  $\varSigma\leadsto\varSigma'\leadsto\varSigma''\leadsto\varSigma'''$: the first one
  in the root $\frac12(\varepsilon_1-\varepsilon_2-\varepsilon_3+\delta)$, the second one
  in the root $\frac12(-\varepsilon_1-\varepsilon_2+\varepsilon_3+\delta)$, the third one
  in the root $\frac12(\varepsilon_1+\varepsilon_2-\varepsilon_3+\delta)$.
  Now we apply~\cite[Corollary 8.6.2]{mu}.
\end{proof}

\subsection{Configuration spaces}
\label{config spaces}
From now on we will work with the adjoint group $\on{F}(4)^{\on{ad}}$ whose weight lattice
coincides with the root lattice $Q$. Its even part $\on{F}(4)^{\on{ad}}_{\bar0}$ is the
quotient $(\on{Spin}(7)\times\SL(2))/\Delta_{\pm1}$ modulo the diagonal subgroup in the center
$\{\pm1\}\times\{\pm1\}\subset\on{Spin}(7)\times\SL(2)$.

We fix a smooth projective curve $C$, we consider the configuration space $C^{\alpha}$ of $Q_{\on{neg}}$-colored
divisors on $C$ of
total degree $\alpha$. If $\alpha=-\sum_{i=1}^4 a_i\alpha_i$ for $a_i\in \BN$, we have $C^{\alpha}\simeq \prod_{i=1}^4 C^{(a_i)}$. We denote by $C^{\bullet}:= \bigsqcup_{\alpha\in Q_{\on{neg}}\setminus 0} C^{\alpha}$.

Similarly, we can consider the configuration space with a marked point. Namely, we choose  a marked point $c\in C$. We denote by $C_{\infty\cdot c}^\bullet:= \bigsqcup_{\bmu\in Q} C^{\bmu}_{\infty\cdot c}$ the space which classifies $Q$-colored divisors on $C$, such that the coefficients of any $x\in C\setminus c$ belongs to $Q_{\on{neg}}$. Here, $C^{\bmu}_{\infty\cdot c}$ is the connected component of total degree $\bmu$.

To be more precise, given a weight $\bmu\in Q$, the configuration space $C^{\bmu}_{\infty\cdot c}$ classifies $Q$-colored
divisors $D=-\varepsilon_1E_1-\varepsilon_2E_2-\varepsilon_3E_3-\delta\Delta$ on $C$ of
total degree $\bmu$ with the following positivity condition.
We have $D=\sum_{x\in C}(-{\mu}_{1,x}\varepsilon_1x-
{\mu}_{2,x}\varepsilon_2x-{\mu}_{3,x}\varepsilon_3x-{\mu}_{4,x}\delta x)$, and $\bmu_x\in Q$ for any
$x\in C$.

We set
\begin{equation}
  \label{D from E}
  D_1=E_1+\Delta,\ D_2=E_1+E_2+E_3+\Delta,\ D_3=2E_1+E_2+E_3,\ D_4=E_1+E_2.
\end{equation}
Equivalently,
\begin{multline}
  \label{E from D}
  E_1=\frac12(D_1-D_2+D_3),\ E_2=\frac12(-D_1+D_2-D_3)+D_4,\ E_3=\frac12(-D_1+D_2+D_3)-D_4,\\
  \Delta=\frac12(D_1+D_2-D_3).
  \end{multline}
These divisors are the coefficients of $D$ in the basis of negative simple roots:
\[D=-D_1\alpha_1-D_2\alpha_2-D_3\alpha_3-D_4\alpha_4.\]
Then we require
\begin{equation}
  \label{effective}
  D_i\ \on{is}\ \on{effective}\ \on{integral}\ \on{valued}\ \on{away}\
  \on{from}\ c\in C\ \on{for}\ 1\leq i\leq 4.
\end{equation}

The space $C^\bullet_{\infty\cdot c}$ has an ind-scheme structure. We have $C_{\infty\cdot c}^\bmu=\bigcup_{\blambda}C^{\bmu}_{\leq\blambda\cdot c}$, where $D=\bmu_c\cdot c- \sum_{C\setminus c}\alpha_x\cdot x \in C^{\bmu}_{\leq\blambda\cdot c}$ if the coefficient $\bmu_c$ of $c$ satisfies
$\bmu_c\leq \blambda$.
We have $C^{\bmu}_{\leq\blambda\cdot c}\simeq C^\alpha=\prod_{i=1}^4C^{(a_i)}$, where $\blambda-{\bmu}=\alpha=
\sum_{i=1}^4a_i\alpha_i$ for $a_i\in\BN$.

The configuration space $C^\bullet$ is a non-unital space with the monoidal structure given by the addition map. Furthermore, the addition map gives $C^\bullet_{\infty\cdot c}$ a $C^\bullet$-module space structure.

\subsection{A factorizable line bundle}
\label{rho}
We have $\rho_0=\frac12(5\varepsilon_1+3\varepsilon_2+\varepsilon_3+\delta),\
\rho_1=\frac12(3\varepsilon_1+\varepsilon_2+\varepsilon_3+\delta),\
\rho=\rho_0-\rho_1=\varepsilon_1+\varepsilon_2$.

We consider a line bundle $\CP$ on $C^{\bmu}_{\infty\cdot c}$ (or, $C^{\alpha}$) with fibers
\begin{multline}
  \label{global product}
  \CP_D=\det{}\!^{-8}\!R\Gamma(C,\bomega_C^{-1/2}(-E_1))
  \otimes\det{}\!^{-8}\!R\Gamma(C,\bomega_C^{-1/2}(-E_2))\\
  \otimes\det{}\!^{-8}\!R\Gamma(C,\bomega_C^{1/2}(-E_3))
  \otimes\det{}\!^{24}R\Gamma(C,\bomega_C^{1/2}(-\Delta))\\
  \otimes\det{}\!^{16}R\Gamma(C,\bomega_C^{-1/2})\otimes\det{}\!^{-16}R\Gamma(C,\bomega_C^{1/2}),
\end{multline}
where $\bomega_C$ is the canonical line bundle on $C$, and we fix a square root $\bomega_C^{1/2}$.
We add the factor
$\det{}\!^{16}R\Gamma(C,\bomega_C^{-1/2})\otimes\det{}\!^{-16}R\Gamma(C,\bomega_C^{1/2})$,
so that $\CP$ is canonically trivial on $C^{-\alpha_i}$ for any simple root $\alpha_i$.
Since the divisors $\Delta,E_i$ can have half-integral coefficients, an explanation is in
order.

Recall that in~\cite[\S 7]{d}, Deligne defined a symmetric bilinear pairing on
$\on{Pic}(C)\colon\langle\CM_1,\CM_2\rangle=\det\!R\Gamma(C,\CM_1\otimes\CM_2)\otimes
\det\!R\Gamma(C,\CO_C)\otimes\det\!{}^{-1}R\Gamma(C,\CM_1)\otimes\det\!{}^{-1}R\Gamma(C,\CM_2)$.
For a line bundle $\CN$ and a divisor $S$ on $C$, we have
\[\det\!R\Gamma(C,\CN(S))=\langle\CN,\CO_C(S)\rangle\otimes\det\!R\Gamma(C,\CN)
\otimes\det\!R\Gamma(C,\CO_C(S))\otimes\det{}\!^{-1}\!R\Gamma(C,\CO_C),\]
\[\det\!R\Gamma(C,\CN(-S))=\langle\CN,\CO_C(-S)\rangle\otimes\det\!R\Gamma(C,\CN)
\otimes\det\!R\Gamma(C,\CO_C(-S))\otimes\det{}\!^{-1}\!R\Gamma(C,\CO_C).\]
Taking product of these equalities, we obtain
\begin{multline*}
\det\!R\Gamma(C,\CN(S))\otimes\det\!R\Gamma(C,\CN(-S))\\
=\det{}\!^2\!R\Gamma(C,\CN)\otimes\det\!R\Gamma(C,\CO_C(S))\otimes\det\!R\Gamma(C,\CO_C(-S))
\otimes\det{}\!^{-2}\!R\Gamma(C,\CO_C).
\end{multline*}
Furthermore,
\[\det\!R\Gamma(C,\CO_C(S))\otimes\det\!R\Gamma(C,\CO_C(-S))=\det{}\!^2\!R\Gamma(C,\CO_C)
\otimes\langle\CO_C(S),\CO_C(S)\rangle.\]
Hence
\[\det{}\!^2\!R\Gamma(C,\CN)=\det\!R\Gamma(C,\CN(S))\otimes\det\!R\Gamma(C,\CN(-S))
\otimes\langle\CO_C(S),\CO_C(-S)\rangle,\]
\[\det{}\!^8\!R\Gamma(C,\CN)=\det{}\!^4\!R\Gamma(C,\CN(S))\otimes\det{}\!^4\!R\Gamma(C,\CN(-S))
\otimes\langle\CO_C(2S),\CO_C(-2S)\rangle.\]
We apply the latter equality to $S=\Delta$, and e.g.\ $\CN=\bomega_C^{-1/2}(-E_i)$.
Then the divisors $\pm2\Delta,\ -E_i\pm\Delta$ are all integral valued, so the RHS of the latter
equality is well defined, and we take it for the definition of the LHS.

\medskip

In other words, the fiber of $\CP$ at 
$D=\sum_{x\in C}(-{\mu}_{1,x}\varepsilon_1x-
{\mu}_{2,x}\varepsilon_2x-{\mu}_{3,x}\varepsilon_3x-{\mu}_{4,x}\delta x)$ is equal to
\begin{equation}
  \label{local product}
  \CP_D=\bigotimes_{x\in C}\bomega_x^{4({\mu}_{1,x}({\mu}_{1,x}-2)+{\mu}_{2,x}({\mu}_{2,x}-2)+
  {\mu}_{3,x}^2-3{\mu}_{4,x}^2)},
\end{equation}
where $\bomega_x$ is the fiber of the canonical line bundle $\bomega_C$ at $x\in C$.
Note that the exponent of $\bomega_x$ is integral for ${\mu}\in Q$.

In particular, the line bundle $\cP$ is a factorizable line bundle, i.e., 
$\on{add}^*(\cP)|_{(C^\alpha\times C^\beta)_{\on{disj}}}\cong 
(\cP\boxtimes \cP)|_{(C^\alpha\times C^\beta)_{\on{disj}}}$ and 
$\on{add}^*(\cP)|_{(C^\alpha\times C^{\bmu}_{\infty\cdot 
c})_{\on{disj}}}\cong (\cP\boxtimes \cP)|_{(C^\alpha\times 
C^{\bmu}_{\infty\cdot c})_{\on{disj}}}$. Here, $(C^\alpha\times 
C^\beta)_{\on{disj}}$ is an open subset of $C^\alpha\times 
C^\beta$ such that the supports of divisors are disjoint, 
$(C^\alpha\times C^{\bmu}_{\infty\cdot c})_{\on{disj}}$ is 
defined similarly but the divisor in $C^\alpha$ is required to 
miss the marked point $c$. 

Thus, it makes sense to consider the $\cP^{\kappa}$-twisted factorizable sheaves.

\subsection{A twisted factorization algebra}
We fix $\kappa\in\BC$ and set $q=\exp(\pi\sqrt{-1}\kappa)$. We assume that $q$ is
transcendental throughout the paper.

We consider the category $\on{D-mod}_{\CP^\kappa}(C^{\bmu})$ of $\CP^\kappa$-twisted $D$-modules
on $C^{\bmu}$. (Recall that all our $D$-modules have {\em super} vector spaces as coefficients.)
We describe the most important (for us) object of this category.
We assume $\alpha=-\sum_{i=1}^4a_i\alpha_i,\ a_i\in\BN$.
Then $C^\alpha\supset\oC^\alpha$: the open subset formed by all the multiplicity
free divisors (i.e.\ each point has multiplicity either~0 or a simple root).
Due to~\eqref{local product} and the fact that $(\alpha_i-2\rho,\alpha_i)=0$ for a simple root
$\alpha_i$, the restriction of the line bundle $\CP$ to $\oC^\alpha$ trivializes canonically.
Under this trivialization, we denote by $\oCI^\alpha$ the $\CP^\kappa$-twisted $D$-module
corresponding to the 1-dimensional local system with trivial monodromies, unless two points
colored by $-\alpha_4$ (the only {\em even} simple root) swap their positions, in which case
the monodromy is~$-1$ (the {\em sign} local system along $\alpha_4$). We assign to $\oCI^\alpha$ the
parity $a_1+a_2+a_3\pmod{2}$. Finally, we define
$\CI^\alpha\in\on{D-mod}_{\CP^\kappa}(C^\alpha)$ as the Goresky-MacPherson extension of $\oCI^\alpha$
to the whole of $C^\alpha$, cf.~\cite[Remark 2.4.2]{ty}. By construction, the
collection $(\CI^\alpha,\ \alpha\in Q_{\on{neg}}\setminus 0)$, forms a factorization algebra in
$\on{D-mod}_{\CP^\kappa}(C^\bullet)$.

\subsection{Factorizable sheaves}
\label{factor shv}
A {\em factorizable sheaf} $\CF$ is a collection of twisted $D$-modules
$\CF^{\bmu}\in\on{D-mod}_{\CP^\kappa}(C_{\infty\cdot c}^{\bmu})$ equipped with factorization isomorphisms
\begin{equation}
  \label{factor sheaves}
  \on{add}^*(\CF^{\alpha+\bmu})|_{\left(C^\alpha\times C^{\bmu}_{\infty\cdot c}
  \right)_{\on{disj}}}
  \cong\left(\CI^\alpha\boxtimes\CF^{\bmu}\right)|_{\left(C^\alpha\times C^{\bmu}_{\infty\cdot c}
  \right)_{\on{disj}}},
\end{equation}
satisfying the appropriate cocycle conditions. 

We also impose the following finiteness conditions:

(a) There is $\blambda\in Q$ such that the support of $\CF^{\bmu}$ lies in $C^{\bmu}_{\leq\blambda\cdot c}$.

(b) There are only finitely many ${\bmu}$ such that the singular support of $\CF^{\bmu}$
contains the conormal to the point ${\bmu}\cdot c$ in $C^{\bmu}_{\leq\blambda\cdot c}$.


The factorizable sheaves with the above finiteness conditions form an abelian category $\on{FS}$
(the morphisms are required to be compatible with the factorization isomorphisms).

One can also allow the marked point $c$ to vary in $C$; moreover, one can allow $n$ 
marked points  to vary in 
$C^n$. The resulting category $\on{FS}_n$~\cite[\S3]{g} is used to
construct a braided monoidal structure on $\on{FS}$ via the nearby cycles functor as the marked
points collide. The following theorem is proved similarly to the main result of~\cite{bfs}.
A conceptual proof is due to J.~Lurie, see the proof of~\cite[Theorem 29.2.3]{gl}
or~\cite[Theorem 1.2.1]{cf} in the classical setup.

\begin{thm}
  \label{bfsl}
  There is a braided monoidal equivalence
  $\Rep_q(\on{F}(4)^{\on{ad}})\simeq\on{FS}$. 
  In particular, for any $\blambda\in Q$ satisfying condition~\eqref{odd dominant},
  the corresponding irreducible $U_q(\ff(4))$-module $V_\blambda$ goes to
  the irreducible factorizable sheaf $\CF_\blambda$. Here, $\CF_\blambda$ is the unique irreducible factorizable sheaf whose support is $C^\bullet_{\leq \blambda\cdot c}$. \hfill $\Box$
\end{thm}

\section{The affine Grassmannian and SW zastava}

\subsection{Renormalized determinant line bundle}
\label{renormalized line bundle}
We will need the following groups and related moduli spaces.

We consider a basis $v_1,\ldots,v_6$ in a 6-dimensional symplectic vector space $V=\BC^6$
such that $\langle v_i,v_{7-i}\rangle=1=-\langle v_{7-i},v_i\rangle$ for $i=1,2,3$, and all
the other products vanish. We view $\varepsilon_1,\varepsilon_2,\varepsilon_3$ as a
natural basis in the space of coweights of $\Sp(6)$, so that the cone of dominant coweights is 
$\{\lambda_1\varepsilon_1+\lambda_2\varepsilon_2 +\lambda_3\varepsilon_3\ :\
\lambda_i\in\BZ,\ \lambda_1\geq\lambda_2\geq\lambda_3\geq0\}$.

The lattice of coweights of $\PSp(6)$ is
$\{\lambda_1\varepsilon_1+\lambda_2\varepsilon_2+ \lambda_3\varepsilon_3\ :\ \lambda_i\in\frac12\BZ,\
\lambda_1-\lambda_2\in\BZ\ni\lambda_2-\lambda_3\}$. The cone 
$\Lambda^+$ of dominant coweights of $\PSp(6)$ is
$\{\lambda_1\varepsilon_1+\lambda_2\varepsilon_2+ \lambda_3\varepsilon_3\ :\
\lambda_i\in\frac12\BN,\ \lambda_1-\lambda_2\in\BN\ni\lambda_2-\lambda_3\}$.

We also denote by $\delta$ the simple coroot of $\SL(2)$, so that the coweight lattice
of $\PGL(2)$ is $\frac12\BZ\delta$.

Let $\on{P}\subset\PSp(6)$ be the stabilizer of the subspace $V_{12}:=\BC v_1\oplus\BC v_2$.
Let $\on{P}_-\subset\PSp(6)$ be the stabilizer of the subspace $V_{56}:=\BC v_5\oplus\BC v_6$,
and let 
$\on{L}:=\on{P}\cap\on{P}_-$ (the common Levi of the opposite parabolic subgroups).
Then $\on{L}$ is the quotient of $\GL(V_{12})\times\Sp(V_{34})$ by the diagonal subgroup
$\pm\Id$ (here $V_{34}:=\BC v_3\oplus\BC v_4$ is a symplectic plane). The cone
$\Lambda^+_{\on{L}}$ of dominant coweights of $\on{L}$ is
$\{\lambda_1\varepsilon_1+\lambda_2\varepsilon_2+ \lambda_3\varepsilon_3\ :\
\lambda_i\in\frac12\BZ,\ \lambda_3\in\frac12\BN,\ \lambda_1-\lambda_2\in\BN,\ \lambda_2-\lambda_3\in\BZ\}$.

Let ${\on{T}}$ denote the diagonal (in the basis $v_1,\ldots,v_6$) Cartan torus of $\PSp(6)$
Recall the coweight $\rho=\varepsilon_1+\varepsilon_2\colon \BG_m\to {\on{T}}$ of~\S\ref{rho}.
Its image is the subgroup $\on{T}':=\diag(x,x,1,1,x^{-1}, x^{-1})\subset \PSp(6)$. 

Clearly, the quotient of $\Sp(V_{34})$ by $\pm\Id$ is $\PGL(2)$, so we obtain the diagonal
embedding $\PGL(2)\hookrightarrow\on{L}$. Recall that $\on{U}$ is the radical of $\on{P}$,
and $\PGL(2)\ltimes\on{U}\subset\on{L}\ltimes\on{U}=\on{P}$ is the subgroup considered
in~\cite[\S3.3]{bft2}. In particular, the quotient $\PSp(6)/\on{U}$ is a spherical
$\PSp(6)\times\PGL(2)$-variety. The $\PGL(2)$-invariant character $\chi^0$ of the Lie algebra
$\fu$ of $\on{U}$ is the sum $\chi^0(u)=u_{13}+u_{24}$ of matrix elements. We keep the same
notation for the corresponding character of the unipotent group $\on{U}$. Then we extend
scalars to $\bF$ and obtain the character $\chi^0_\bF\colon\on{U}(\bF)\to\bF$. It extends to the
same named character $\PGL(2,\bO)\ltimes\on{U}(\bF)\to\bF$. Finally, the character
$\chi\colon\PGL(2,\bO)\ltimes\on{U}(\bF)\to\BG_a$ is the residue of $\chi^0_\bF$.

\begin{defn}
We set $\widetilde{G}:=\SL(2)\times\Sp(6),\ \ol{G}:=\PGL(2)\times\PSp(6),\
  G:=(\SL(2)\times\Sp(6))/\{\pm1\}$ (quotient modulo the diagonal subgroup in the center).
  We set $\widetilde{H}:=\SL(2)\ltimes\on{U}\hookrightarrow\widetilde{G}$ (the diagonal embedding),
  it gives rise to the embeddings $H:=\PGL(2)\ltimes\on{U}\hookrightarrow G$ and
  $H\hookrightarrow\ol{G}$.
\end{defn}

The moduli space $\Bun_{\PSp(6)}$ of $\PSp(6)$-bundles on $C$ is the moduli space of rank~6
vector bundles $\CV$ on $C$ equipped with a non-degenerate skew-symmetric pairing
$\Lambda^2\CV\to\CM$ with values in a line bundle $\CM$, considered modulo the following
equivalence relation. For a line bundle $\CR$ on $C$, the twisted vector bundle $\CV\otimes\CR$
inherits a non-degenerate skew-symmetric pairing $\Lambda^2(\CV\otimes\CR)\to\CM\otimes\CR^2$,
and we declare
$(\CV,\ \Lambda^2\CV\to\CM)\sim(\CV\otimes\CR,\ \Lambda^2(\CV\otimes\CR)\to\CM\otimes\CR^2)$.
According to the following lemma, the renormalized determinant line bundle
$\CP_{\det,\on{GSp}(6)}$ whose fiber at $(\CV,\ \Lambda^2\CV\to\CM)$ is
\begin{equation*}
    \begin{split}
        \det\!{}^{-4}\!R\Gamma(C,\CV)\otimes\det\!{}^{6}\!R\Gamma(C,\CM\otimes\bomega_C^{-1/2})\\ \otimes \on{det}^{16}R\Gamma(C, \CO)\otimes \on{det}^8 R\Gamma(C, \bomega^{-1}_C)\otimes \on{det}^{-6}R\Gamma(C, \bomega_C^{-\frac{1}{2}}),
    \end{split}
\end{equation*}
descends from $\Bun_{\on{GSp}(6)}$ to $\Bun_{\PSp(6)}$.

\begin{lem}
  \label{independence}
  We have canonical isomorphisms
  \[\det\!{}^2\!R\Gamma(C,\CV)\otimes\det\!{}^{-3}\!R\Gamma(C,\CM\otimes\bomega_C^{-1/2})= 
  \det\!{}^2\!R\Gamma(C,\CV\otimes\CR)\otimes
  \det\!{}^{-3}R\Gamma(C,\CM\otimes\bomega_C^{-1/2}\otimes\CR^2),\]
  satisfying an appropriate cocycle condition.
\end{lem}

\begin{proof}  
We have 
\[\det\!R\Gamma(C,\CR^2)=\langle\CR,\CR\rangle\otimes\det\!{}^2\!R\Gamma(C,\CR)
\otimes\det\!{}^{-1}\!R\Gamma(C,\CO_C),\]
\[\det\!R\Gamma(C,\bomega_C)=\langle\CR,\bomega_C\otimes\CR^{-1}\rangle\otimes\det\!{}^2\!R\Gamma(C,\CR)
\otimes\det\!{}^{-1}\!R\Gamma(C,\CO_C),\]
since $\det\!R\Gamma(C,\CR)=\det\!R\Gamma(C,\bomega_C\otimes\CR^{-1})$. Plugging the latter formula
for $\det\!{}^2\!R\Gamma(C,\CR)$ into the former formula, we obtain
\[\det\!R\Gamma(C,\CR^2)=\langle\CR,\CR\rangle^2\otimes\langle\bomega_C^{-1},\CR\rangle\otimes
\det\!R\Gamma(C,\CO_C).\]
Replacing $\CR$ by a square root of $\CR$, we get
\[\det\!{}^2\!R\Gamma(C,\CR)=\langle\CR,\CR\rangle\otimes\langle\bomega_C^{-1},\CR\rangle\otimes
\det\!{}^2\!R\Gamma(C,\CO_C).\]
Furthermore, \[\det\!{}^2\!R\Gamma(C,\CV\otimes\CR)=\langle\det\CV,\CR\rangle^2\otimes
\det\!{}^2\!R\Gamma(C,\CV)\otimes\det\!{}^{12}\!R\Gamma(C,\CR)
\otimes\det\!{}^{-12}\!R\Gamma(C,\CO_C),\]
\begin{multline*}
  \det\!{}^{-3}\!R\Gamma(C,\CM\otimes\bomega_C^{-1/2}\otimes\CR^2)\\
=\langle\CM\otimes\bomega_C^{-1/2},\CR^2\rangle^{-3}\otimes
\det\!{}^{-3}\!R\Gamma(C,\CM\otimes\bomega_C^{-1/2})\otimes\det\!{}^{-3}\!R\Gamma(C,\CR^2)
\otimes\det\!{}^3\!R\Gamma(C,\CO_C).
\end{multline*}
Taking product of these two equalities, we obtain the desired result.
\end{proof}

The embedding $\PGL(2)\hookrightarrow \PSp(6)$ induces a map of stacks
$\Bun_{\PGL(2)}\rightarrow \Bun_{\PSp(6)}$. We denote the pullback of $\CP_{\det, \PSp(6)}$ to
$\Bun_{\PGL(2)}$ by $\CP_{\det, \PGL(2)}$. To be more precise, its fiber over a ${\PGL(2)}$-bundle, i.e.\
a rank 2 vector bundle $\CW$ modulo tensoring by a line bundle, is 
\begin{equation*}
    \begin{split}
      \on{det}^{-12}(C,\CW)\otimes \on{det}^{6}(C, \det \CW\otimes \bomega_C^{-\frac{1}{2}})\otimes
      \on{det}^{24}R\Gamma(C, \CO)\otimes \on{det}^{-6} R\Gamma(C, \bomega_C^{-\frac{1}{2}})
    \end{split}
\end{equation*}

 We define $\CP_{\det}$ on
$\Bun_{\overline{G}}=\Bun_{\PGL(2)}\times \Bun_{\PSp(6)}$ to be
$\CP_{\det, \PGL(2)}^{-1}\boxtimes \CP_{\det, \PSp(6)}$.


With some abuse of notation, we keep the same
notation for its further pullback to other spaces, such as the affine Grassmannian
${}^{\bomega}\Gr_{\overline{G}}$ in~\S\ref{gaiotto category}, and the global model $\scrM_{\infty\cdot c}$
in~\S\ref{global model}.

\subsection{Gaiotto category}
\label{gaiotto category}
In order to have a canonically defined character $\chi$, we need to consider the \textit{renormalized
affine Grassmannian} ${}^{\bomega}\Gr_{\overline{G}}$.

Recall from~\S\ref{renormalized line bundle} that ${\on{T}}$ denotes the diagonal Cartan torus of
$\PSp(6)$, and $\on{T}':=\diag(x,x,1,1,x^{-1}, x^{-1})\subset \PSp(6)$ is the image of the coweight
$\rho=\varepsilon_1+\varepsilon_2\colon \BG_m\to {\on{T}}$.
We denote by $\bomega_C^\rho$ the $\on{T}'$-bundle induced from the $\BG_m$-bundle $\bomega_C$.

We define the $\bomega_C^\rho$-renormalized affine Grassmannian ${}^\bomega\Gr_{\overline{G}}$
as the moduli space of pairs $(\CP_G, \alpha)$, where $\CP_G=(\CP_{\PGL(2)}, \CP_{\PSp(6)})$ is a
$\overline{G}$-bundle on $C$, and $\alpha$ is an isomorphism of $\CP_{\overline{G}}$ with
$\CP_{\overline G}^{\bomega}\colon= \bomega_C^\rho\overset{\on{T}'}{\times} \overline{G}=
(\CP_{\PGL(2)}^{\on{triv}}, \CP_{\PSp(6)}^{\bomega})$ on $C\setminus c$.
Here $\CP_{\PGL(2)}^{\on{triv}}$ is the trivial
$\PGL(2)$-bundle and $\CP_{\PSp(6)}^{\bomega}=\bomega_C^\rho\overset{\on{T}'}{\times} \PSp(6)$.
The renormalized affine Grassmannian ${}^\bomega\Gr_{\overline{G}}$ can also be defined as the affine
Grassmannian twisted by the $\on{T}'$-bundle $\bomega^{\rho}_C$, i.e.\
$\Gr_{\overline{G}}\overset{\on{T}'(\bO)}{\times} \bomega_C^{\rho}|_{\cD}$.

By definition, $^\bomega\Gr_{\overline{G}}$ is (the \'etale sheafitification
of) the stacky quotient $^\bomega\overline{G}(\bF)/{}^\bomega\overline{G}(\bO)$, where
$^\bomega\overline{G}(\bF)$ (resp.\ $^\bomega\overline{G}(\bO)$) denotes the automorphisms of
$\CP_{\overline{G}}^{\bomega}$ on the punctured formal disc $\overset{\circ}{\cD}$ (resp.\ formal disc
${\cD}$). Similarly, we can define the $\bomega^\rho_C$-renormalized groups $^\bomega\!H(\bF)$ and $^\bomega\!\on{U}(\bF)$. 

The character
\begin{equation*}
    \chi^0\colon \on{U}\longrightarrow \BG_a,\ 
        u\mapsto u_{1,3}+u_{2,4}
\end{equation*}
induces a character of $^\bomega\!\on{U}(\bF)$
\begin{equation}
    \begin{split}
         \chi_\bF\colon ^\bomega\!\on{U}(\bF)\longrightarrow \bomega_C|_{\overset{\circ}{\cD}}.
    \end{split}
\end{equation}
The composition with the residue map induces a character $\chi$ of $^\bomega\!\on{U}(\bF)$ that is
invariant under the conjugation action of $\PGL(2,\bF)$. In particular, it defines a character of
$^\bomega\!H(\bF)$. With some abuse of notations, we still denote the resulting character of
$^\bomega\!H(\bF)$ by $\chi$.

We define the twisted Gaiotto category
$\cC_{\kappa}^{\loc}:=\on{D-mod}^{^\bomega\!H(\bF),\chi}_{\CP_{\det}^{\kappa}}({}^{\bomega}\Gr_{\overline{G}})$ as the category of
$\CP_{\det}^{\kappa}$-twisted ($^\bomega\!H(\bF),\chi$)-equivariant $D$-modules on $^\bomega\Gr_{\overline{G}}$.
(Recall that all our $D$-modules have {\em super} vector spaces as coefficients.) It is equivalent to
$\on{D-mod}^{\PGL(2,\bO)\ltimes{}^\bomega\!\on{U}(\bF),\chi}_{\CP_{\det, \PSp(6)}^{\kappa}}({}^{\bomega}\Gr_{\PSp(6)})$,
i.e.\  the category of  $\CP_{\det, \PSp(6)}^{\kappa}$-twisted
($\PGL(2,\bO)\ltimes{}^\bomega\!\on{U}(\bF),\chi$)-equivariant $D$-modules on ${}^{\bomega}\Gr_{\PSp(6)}$. Here, the superscript ``$\loc$" means ``local", we will also introduce a global
version of the twisted Gaiotto category in~\S\ref{3.5.2}. Sometimes, in order to emphasize the chosen point $c$ in the definition of the affine Grassmannian, we will denote $\cC_{\kappa}^{\loc}$ by $\cC_{\kappa,c}^{\loc}$.

We denote by $\cC_{\kappa}^{\loc,\lc}:= \on{D-mod}^{\PGL(2,\bO)\ltimes{}^\bomega\!\on{U}(\bF),\chi,\lc}_{\CP_{\det, \PSp(6)}^{\kappa}}({}^{\bomega}\Gr_{\PSp(6)})$ the category of \textit{locally compact} Gaiotto equivariant twisted D-modules, i.e., it consists of objects which are compact after applying the forgetful functor $\on{D-mod}^{\PGL(2,\bO)\ltimes{}^\bomega\!\on{U}(\bF),\chi}_{\CP_{\det, \PSp(6)}^{\kappa}}({}^{\bomega}\Gr_{\PSp(6)})\rightarrow \on{D-mod}^{{}^\bomega\!\on{U}(\bF),\chi}_{\CP_{\det, \PSp(6)}^{\kappa}}({}^{\bomega}\Gr_{\PSp(6)})$.

\subsection{Classification of irreducible twisted $D$-modules}
\label{matrix g}
Recall the notation of~\cite[\S3.3]{bft2}. In particular, $q=\exp(\pi\sqrt{-1}\kappa)$ for
an irrational number $\kappa$ (denoted $c^{-1}$ in {\em loc.cit.}).

We first address the question when the closure of a $\PGL(2,\bO)\ltimes {}^{\bomega}\!\on{U}(\bF)$-orbit in ${}^\bomega \Gr_{\PSp(6)}$
supports an irreducible $\chi$-equivariant $\kappa$-twisted $D$-module. To avoid cumbersome notation,
we will consider the corresponding (and equivalent) question for the usual (non-renormalized) affine
Grassmannian $\Gr_{\PSp(6)}$. 

To this end we
consider a matrix $g\in\PSp(6)$ equal to (the image under the projection
$\Sp(6)\to\PSp(6)$ of) the sum $\Id+E_{34}$ (the identity matrix plus an elementary matrix).
Given a cocharacter $(\lambda_1,\lambda_2,\lambda_3)$ (resp.\ $\lambda_4\in\frac12\BZ$) of $\PSp(6)$
(resp.\ of $\PGL(2)$) we view $t^{\lambda_4}\in\PGL(2,\bF)$ as an element of $\PSp(6,\bF)$ under
the natural embedding $\PGL(2,\bF)\hookrightarrow\PSp(6,\bF)$, and consider the point
$t^{\lambda_4}gt^{(\lambda_1,\lambda_2,\lambda_3)}\in\Gr_{\PSp(6)}$.

\begin{prop}
\label{supports PSp}
\textup{(a)} The set of points $t^{\lambda_4}gt^{(\lambda_1,\lambda_2,\lambda_3)}$ where
$(\lambda_1,\lambda_2,\lambda_3)\in\Lambda^+_{\on{L}}$ and $\lambda_4\in\frac12\BN$, is the full set
of representatives of $\PGL(2,\bO)\ltimes\on{U}(\bF)$-orbits in $\Gr_{\PSp(6)}$.

\textup{(b)} The set of points $t^{\lambda_4}gt^{(\lambda_1,\lambda_2,\lambda_3)}$ where
$(\lambda_1,\lambda_2,\lambda_3)\in\Lambda^+$ and $\lambda_4\in\frac12\BN$, is the full set
of representatives of $\PGL(2,\bO)\ltimes\on{U}(\bF)$-orbits in $\Gr_{\PSp(6)}$ whose
closure supports a $\chi$-equivariant $D$-module.

\textup{(c)} The set of points $t^{\lambda_4}gt^{(\lambda_1,\lambda_2,\lambda_3)}$ where
$(\lambda_1,\lambda_2,\lambda_3)\in\Lambda^+$ and $\lambda_4\in\frac12\BN$ satisfy the
condition~\eqref{odd dominant}, is the full set
of representatives of $\PGL(2,\bO)\ltimes\on{U}(\bF)$-orbits in $\Gr_{\PSp(6)}$ whose
closure supports a $\chi$-equivariant $\kappa$-twisted $D$-module.
\end{prop}
\begin{proof}
    (a) Note that for geometric points, we have $|\on{P}(\bF)/\on{P}(\bO)|=|\PSp(6,\bF)/\PSp(6,\bO)|$, so it is enough to classify $\PGL(2,\bO)\ltimes \on{U}(\bF)$-orbits in $\on{P}(\bF)/\on{P}(\bO)$, which is further enough to classify $\PGL(2,\bO)$-orbits in $\Gr_{\on{L}}$. 
    
    We will classify $\SL_2(\bO)$-orbits in $\Gr_{\GL_2\times \SL_2}$. In this case, we should prove that the set of points $t^{\lambda_4} g t^{(\lambda_1,\lambda_2,\lambda_3)}$ is the full set of representatives of orbits, where we require one more condition $\lambda_3-\lambda_4\in \BZ$. The proof of classification for $\PGL(2,\bO)$-orbits in $\Gr_{\on{L}}$ is similar.

    Given a point of $\Gr_{\GL_2\times \SL_2}$, by the usual Bruhat decomposition, we can find a representative $(L_1, L_2)\in \GL_2(\bF)\times \SL_2(\bF)$ in its $\SL_2(\bO)$-orbit, such that $L_1=\begin{pmatrix}
        t^{n_1}&0\\ 0& t^{n_2}
    \end{pmatrix}$, and $n_1\geq n_2\in \BZ$. We denote by $k:= n_1-n_2\in \BN$. We claim that the  $\Stab_{\SL_2(\bO)}(L_1)$-orbits in $\Gr_{\SL_2}$ is indexed by the set \begin{equation}\label{condition n}
        \{n_3\in \BN, n_4\in \BZ\mid -n_3\leq n_4\leq n_3,\ n_4\leq -n_3+k\}.
    \end{equation}

    Indeed, there is $\Stab_{\SL_2(\bO)}(L_1)= L_1\GL_2(\bO)L_1^{-1}\cap \SL_2(\bO)=\{\begin{pmatrix}\bO & t^k \bO\\ \bO& \bO           \end{pmatrix}\subset \SL_2(\bO)\}$. Using the fact $|\SL_2(\bF)/\SL_2(\bO)|=|\on{B}_{\SL_2}(\bF)/\on{B}_{\SL_2}(\bO)|$, we can find an upper-triangular representative for any chosen  $\Stab_{\SL_2(\bO)}(L_1)$-orbit in $\Gr_{\SL_2}$. By left-multiplying an element of the form $\begin{pmatrix}
\bO^\times&0\\
\bO& \bO^\times
    \end{pmatrix}$ and right-multiplying an element of the form $\begin{pmatrix}
        \bO^\times&0\\
        \bO& \bO^\times
    \end{pmatrix}$, we can assume that the $(2,2)$-entry of the representative is non-zero, and the degree is the lowest. Furthermore, by left-multiplying an element of the form  $\begin{pmatrix}
\bO^\times& t^k \bO\\ 0& \bO^\times
    \end{pmatrix}$ and right-multiplying an element of the form $\begin{pmatrix}
\bO^\times&\bO\\ 0& \bO^\times
    \end{pmatrix}$, we can get a representative $L_2=\begin{pmatrix}
t^{n_3}&t^{n_4}\\ 0& t^{-n_3}
    \end{pmatrix}$, which satisfies the condition \eqref{condition n}.

    Then, we prove that for any two points $L_2=\begin{pmatrix}
t^{n_3}&t^{n_4}\\ 0& t^{-n_3}
    \end{pmatrix}, L_2'=\begin{pmatrix}
t^{n'_3}&t^{n'_4}\\ 0& t^{-n'_3}
    \end{pmatrix}$ which satisfy the condition \eqref{condition n}, they lie in different $\Stab_{\SL_2(\bO)}(L_1)$-orbits in $\Gr_{\SL_2}$. First, we have $n_3=n_3'$, because the lowest degree of entries in a representative of $\SL_2(\bO)$-orbit in $\Gr_{\SL_2}$ is invariant. Assume that $L_2$ and $L_2'$ are in the same orbit, that says there exists $\begin{pmatrix}
        e & t^k f\\
        f & h
    \end{pmatrix}\in \Stab_{\SL_2(\bO)}(L_1)$, and $\begin{pmatrix}
        a & b\\
        c & d
    \end{pmatrix}\in \SL_2(\bO)$, such that
    \begin{equation}
        \begin{pmatrix}
t^{n_3}&t^{n_4}\\
& t^{-n_3}
    \end{pmatrix} \begin{pmatrix}
        a & b\\
        c & d
    \end{pmatrix}=\begin{pmatrix}
        e & t^k f\\
        f & h
    \end{pmatrix}\begin{pmatrix}
t^{n'_3}&t^{n'_4}\\
& t^{-n'_3}
    \end{pmatrix}
    \end{equation}. 

In other words, we have
\begin{equation}\label{3.1.3}
    \begin{pmatrix}
    t^{n_3} a+t^{n_4} c& t^{n_3}b+t^{n_4}d\\
    t^{-n_3}c & t^{-n_3}d
\end{pmatrix}= \begin{pmatrix}
    t^{n_3}e& t^{n'_4} e+t^{k-n_3} f\\ t^{n_3} g& t^{n'_4}g+t^{-n_3}h
\end{pmatrix}
\end{equation}.

In the case $k=0$, the classification of $\Stab_{\SL_2(\bO)}(L_1)=\SL_2(\bO)$-orbits in $\Gr_{\SL_2}$ is given by the usual Bruhat decomposition, it matches the classification \eqref{condition n}. So, we can assume $k>0$. In this case, we have $e, h\in \bO^\times$. Also, we can assume $n_3>0$, otherwise $n_3=0$ implies $n_4=n_4'=0$.

First, we claim $d\in \bO^\times$. Otherwise, since $t^{-n_3}c=t^{n_3}g$, there is $c\in t^{2n_3}\bO$.
In particular, $t^{-n_3+n_4}cd\in t\bO$ and $t^{-n_3}c t^{n_3}b\in t^{2n_3}\bO\subset t\bO$. It forces
$t^{n_3}a t^{-n_3}d\in \bO^\times$, which contradicts the assumption $d\notin \bO^\times$.

Then, we assume $n'_4>n_4$. Since $k-n_3\geq n_4'$, we have $t^{n'_4}e+t^{k-n_3}f\in t^{n_4'}\bO$.
Since $d\in \bO^\times$, we have $t^{n_4}d\in t^{n_4}\bO^\times$, so that $n_3=n_4$. However,
$n_3\geq n'_4>n_4$. We obtain a contradiction.

In conclusion, we obtain that the $\SL_2(\bO)$-orbits in $\Gr_{\GL_2\times \SL_2}$ are indexed by the set $\{n_i\in \BZ\mid n_1\geq n_2, -n_3\leq n_4\leq n_3, n_4\leq -n_3+n_1-n_2\}$. Now, we claim it is the same as the index set $\{\lambda_i\in \frac{1}{2}\BZ\mid (\lambda_1,\lambda_2,\lambda_3)\in\Lambda^+_{\on{L}}, \lambda_4\in\frac12\BN, \lambda_3-\lambda_4\in \BZ\}$.

Indeed, the bijection maps are given by 
\begin{equation}
    \begin{split}
        (n_1, n_2, n_3, n_4)\mapsto (\lambda_1, \lambda_2, \lambda_3, \lambda_4)= (&n_1-\frac{n_3+n_4}{2}, n_2+\frac{n_3+n_4}{2}, \frac{n_3-n_4}{2}, \frac{n_3+n_4}{2})\\
        (n_1, n_2, n_3, n_4)=(\lambda_1+\lambda_4, &\lambda_2-\lambda_4, \lambda_3+\lambda_4, -\lambda_3+\lambda_4)\mapsfrom (\lambda_1, \lambda_2, \lambda_3,\lambda_4).
    \end{split}
\end{equation}

\medskip

    (b) Note that the character $\chi$ is given by taking the residue of the sum of the entries (1,3) and (2,4). We should prove that \begin{equation}
        \Stab_{\on{U}(\bF)}(t^{\lambda_4}g t^{(\lambda_1,\lambda_2,\lambda_3)})\subset \Ker \chi,
    \end{equation}
i.e., the $t^{-1}$-coefficient of the sum of the entries (1,3) and (2,4) is zero. In fact, it is
not hard to see that it is equivalent to the vanishing of the $t^{-1}$-coefficient of each of these
two entries. 

  Let $A\in  \Stab_{\on{U}(\bF)}(t^{\lambda_4}g t^{(\lambda_1,\lambda_2,\lambda_3)})$, we only care about the $2\times 2$-minor $A_{(1,2)\times (3,4)}= \begin{pmatrix}
      u_{1,3} & u_{1,4}\\ u_{2,3}& u_{2,4}
  \end{pmatrix}$. It belongs to 
  \[\begin{pmatrix}
      t^{\lambda_1+\lambda_4}&0\\ 0& t^{\lambda_2-\lambda_4}
  \end{pmatrix}\begin{pmatrix}
      \bO&\bO\\
      \bO&\bO
  \end{pmatrix}\begin{pmatrix}
      t^{-\lambda_3-\lambda_4}& -t^{-\lambda_3+\lambda_4}\\0& t^{\lambda_3+\lambda_4}
  \end{pmatrix}.\] The $t^{-1}$-coefficient of $u_{1,3}$ is always zero if and only if $\lambda_1-\lambda_3\geq 0$, and the $t^{-1}$-coefficient of $u_{2,4}$ is always zero if and only if $\lambda_2-\lambda_3\geq 0$.

  \medskip
  
  (c) We need to calculate the reductive part of $\Stab_{\PGL(2, \bO)}(t^{\lambda_4}g t^{(\lambda_1,\lambda_2,\lambda_3)})$. By definition, for each element in the stabilizer of $t^{\lambda_4}g t^{(\lambda_1,\lambda_2,\lambda_3)}$ in $\PGL(2, \bO)$, it belongs to 
  \begin{equation}\label{3.1.6}
      \begin{pmatrix}
          \bO& t^{\lambda_1-\lambda_2+2\lambda_4}\bO\\
          \bO& \bO
      \end{pmatrix} \cap \begin{pmatrix}
          1 & t^{2\lambda_4}\\
          0 & 1
      \end{pmatrix} \begin{pmatrix}
          \bO& t^{2(\lambda_3+\lambda_4)}\bO\\
         \bO & \bO
      \end{pmatrix}\begin{pmatrix}
          1 & -t^{2\lambda_4}\\
          0 & 1
      \end{pmatrix}.
  \end{equation}
 
  Assume $\lambda_4=0$, and $\lambda_1-\lambda_2, \lambda_3$ are non-zero, then the reductive part of the stabilizer
  group is $\begin{pmatrix}
      a&0 \\
      -a^{-1}+a& a^{-1}
  \end{pmatrix}/ \{\pm \Id\}$, where $a\in \BG_m$. If $\lambda_1-\lambda_2=0$, and $\lambda_3>0$, the reductive
  part of the stabilizer is 
  $\begin{pmatrix}
      a^{-1}&a-a^{-1} \\
      0& a
  \end{pmatrix}/ \{\pm \Id\}$. If $\lambda_1-\lambda_2>0$, and $\lambda_3=0$, the reductive part of the
  stabilizer is 
  $\begin{pmatrix}
      a&0\\
      0&a^{-1}
  \end{pmatrix}/ \{\pm \Id\}$. In these three cases, { the action on the fiber is trivial only if $\lambda_1-\lambda_2-\lambda_3=0$}. If  $\lambda_1-\lambda_2=0$, and $\lambda_3=0$, the reductive part of the stabilizer is 
  $\PGL(2)$, it acts trivially on the fiber.

Then, we assume $\lambda_4>0$.

In the case of $\lambda_1-\lambda_2>0$ and $\lambda_3>0$, \eqref{3.1.6} is a pro-unipotent subgroup of $\PGL(2)(\bO)$, which acts trivially on the fiber. If $\lambda_3=0$, then \begin{equation}\eqref{3.1.6}= 
      \begin{pmatrix}
          \bO^\times& t^{\lambda_1-\lambda_2+2\lambda_4}\bO\\
          \bO& \bO^\times
      \end{pmatrix},
\end{equation}
it has a $\BG_m$-factor
$\begin{pmatrix}
    a & 0\\
    0 & a^{-1}
\end{pmatrix}/\{\pm \Id\}$, the action on the fiber is trivial only if $\lambda_1-\lambda_2+3\lambda_4=0$,
which contradicts the assumption.

In the case of $\lambda_1=\lambda_2$, we have $ \begin{pmatrix}
         1 & t^{2\lambda_4}\\ 0& 1
    \end{pmatrix}t^{\lambda_4}g t^{(\lambda_1,\lambda_2,\lambda_3)}= t^{\lambda_4}t^{(\lambda_1,\lambda_2,\lambda_3)}$, and \eqref{3.1.6} has a $\BG_m$-factor
$\begin{pmatrix}
         1 & t^{2\lambda_4}\\ 0& 1
    \end{pmatrix}\begin{pmatrix}
    a & 0\\
    0 & a^{-1}
\end{pmatrix}\begin{pmatrix}
         1 & -t^{2\lambda_4}\\ 0& 1
    \end{pmatrix}/\{\pm \Id\}$. {The action on the fiber is trivial only if $3\lambda_4+\lambda_3=0$. However, $3\lambda_4+\lambda_3>0$ by assumption}.
\end{proof}

\begin{defn}
  \textup{(a)} Given $\blambda\in X$ satisfying the condition of~Proposition~\ref{supports PSp}(a),
  we denote by $\BO^\blambda\subset\Gr_{\PSp(6)}$ the corresponding $\PGL(2,\bO)\ltimes\on{U}(\bF)$-orbit.
  Under the natural bijection between the set of $\PGL(2,\bO)\ltimes\on{U}(\bF)$-orbits in $\Gr_{\PSp(6)}$
  and the set of $H(\bF)$-orbits in $\Gr_{\PGL(2)}\times\Gr_{\PSp(6)}$, it corresponds to the $H(\bF)$-orbit
  $\sO^\blambda\subset\Gr_{\PGL(2)}\times\Gr_{\PSp(6)}$.

  \textup{(b)} A weight  $\blambda\in X$ satisfying the condition of~Proposition~\ref{supports PSp}(a)
  (resp.\ (b,c)) is called relevant (resp.\ $\chi$-relevant, $(\chi,\kappa)$-relevant). The orbits
  $\BO^\blambda,\sO^\blambda$ for a $\chi$-relevant (resp.\ $(\chi,\kappa)$-relevant) weight $\blambda$ are called
  $\chi$-relevant (resp.\ $(\chi,\kappa)$-relevant).
\end{defn}

\begin{rem}
    We keep the same notation for the corresponding orbits in the $\bomega$-renormalized affine Grassmannian ${}^{\bomega}\Gr_{\PSp(6)}$ and ${}^{\bomega}\Gr_{\overline{G}}$.
\end{rem}

\begin{rem}
    An object in the twisted Gaiotto category is locally compact if and only if it is supported on finitely many relevant orbits, and the stalk on each orbit is finite dimensional. 
\end{rem}

\subsection{Adjacency of orbits}
The closure relation of $\PGL(2, \bO)\ltimes \on{U}(\bF)$-orbits in $\Gr_{\PSp(6)}$ is equivalent to the closure relation of $\PGL(2, \bO)\ltimes {}^{\bomega}\!\on{U}(\bF)$-orbits in ${}^{\bomega}\Gr_{\PSp(6)}$. It is enough to study the adjacency of orbits in the non-renormalized affine Grassmannian.
\begin{prop}
  \label{adjacency}
  If $\BO^{\blambda'}$ lies in the closure of $\BO^\blambda\subset\Gr_{\PSp(6)}$ (equivalently,
  $\sO^{\blambda'}\subset\ol\sO{}^\blambda\subset\Gr_{\ol{G}}$),
  then ${\blambda'}\leq\blambda$, i.e.\ $\blambda-{\blambda'}\in Q_{\on{pos}}$, i.e.\
  $\blambda-{\blambda'}=\sum_{i=1}^4a_i\alpha_i,\ a_i\in\BN$.
\end{prop}

\begin{proof}
We construct a family of lower semi-continuous $\wt{H}(\bF)$-invariant functions on $\Gr_{\wt{G}}$.
We denote by $\fL$ the line bundle $\CO_{\BP^1}(1)$ on the flag variety $\CB_{\SL(2)}=\BP(W)$.
It corresponds to the fundamental weight $\delta^*$ of $\SL(2)$. Its dual is the tautological
line subbundle $\fL^*\subset\CW=W\otimes\CO_{\CB_{\SL(2)}}$. We denote by $\fL_1,\fL_2,\fL_3$ the line
bundles on the flag variety $\CB_{\Sp(6)}$ corresponding to the weights
$\varepsilon_1^*,\varepsilon_2^*,\varepsilon_3^*$ of $\Sp(6)$. That is, if
$V\otimes\CO_{\CB_{\Sp(6)}}=:\CV\supset\CV^1\supset\CV^2\supset\CV^3\supset\CV^4\supset\CV^5\supset0$
is the tautological complete self-orthogonal flag in the tautological symplectic vector bundle on
$\CB_{\Sp(6)}$, then $\fL_1=\CV/\CV^1,\ \fL_2=\CV^1/\CV^2,\ \fL_3=\CV^2/\CV^3$.

Also recall that $\wt{G}/\wt{H}\cong\Sp(6)/\on{U}$ is the moduli space of the following data:

(a) an isotropic plane $\theta\colon W\hookrightarrow V$;

(b) isomorphisms $\zeta\colon W\iso W^\perp/W\iso V/W^\perp$.

\noindent Note that $\zeta$ gives rise to an isomorphism $(\det W)^{\otimes3}\iso\det V=\BC$, and
we choose a compatible trivialization of $\det W$.

We have the following $\wt{G}$-equivariant morphisms:
\begin{multline}
  \label{varphi}
  \wt{G}/\wt{H}\xrightarrow{\varphi_1}\Gamma(\CB_{\wt{G}},\fL\boxtimes\fL_1)
=V^{\delta^*}\boxtimes V^{\varepsilon_1^*},\\
\wt{G}/\wt{H}\xrightarrow{\varphi_2}\Gamma(\CB_{\wt{G}},\fL\boxtimes(\fL_1\otimes\fL_2\otimes\fL_3))
=V^{\delta^*}\boxtimes V^{\varepsilon_1^*+\varepsilon_2^*+\varepsilon_3^*},\\
\wt{G}/\wt{H}\xrightarrow{\varphi_3}\Gamma(\CB_{\Sp(6)},\fL_1^{\otimes2}\otimes\fL_2\otimes\fL_3)
=V^{2\varepsilon_1^*+\varepsilon_2^*+\varepsilon_3^*},\\
\wt{G}/\wt{H}\xrightarrow{\varphi_4}\Gamma(\CB_{\Sp(6)},\fL_1\otimes\fL_2)
=V^{\varepsilon_1^*+\varepsilon_2^*}
\end{multline}
(irreducible representations with given highest weights. E.g.\ $V^{\delta^*}$ is the tautological
representation of $\SL(2)$ with basis $v_{\on{high}},v_{\on{low}}$, and $V^{\varepsilon_1^*}$ is the tautological representation $V$ of $\Sp(6)$ with basis $v_1,\ldots,v_6$).
Since $\varphi_i$ is $\wt{G}$-equivariant, it is uniquely determined by $\varphi_i(1)$: an
$\wt{H}$-invariant vector in the corresponding irreducible $\wt{G}$-module.
Namely,
\begin{multline}
  \label{varphi(1)}
\varphi_1(1)=v_{\on{low}}\boxtimes v_6+v_{\on{high}}\boxtimes v_5,\
\varphi_2(1)=v_{\on{low}}\boxtimes v_4\wedge v_5\wedge v_6+v_{\on{high}}\boxtimes v_3\wedge v_5\wedge v_6,\\
\varphi_3(1)=v_5\otimes v_4\wedge v_5\wedge v_6+v_6\otimes v_3\wedge v_5\wedge v_6,\
\varphi_4(1)=v_5\wedge v_6
\end{multline}
(here we view $V^{\varepsilon_1^*+\varepsilon_2^*}$ (resp.\ $V^{\varepsilon_1^*+\varepsilon_2^*+\varepsilon_3^*},\
V^{2\varepsilon_1^*+\varepsilon_2^*+\varepsilon_3^*}$) as a quotient of $\Lambda^2V$ (resp.\ $\Lambda^3V,\
V\otimes\Lambda^3V$).

More geometrically, the value $\varphi_1(W,\theta,\zeta)$ is the composition
\(\fL^*\hookrightarrow\CW\xrightarrow{\theta}\CV\twoheadrightarrow\fL_1,\) and the value
$\varphi_4(W,\theta,\zeta)$ is the determinant of the composition
$\CW\xrightarrow{\theta}\CV\twoheadrightarrow\CV/\CV^2$. Recall that
$\det(\CV/\CV^2)=\fL_1\otimes\fL_2$, and $\det\CW$ is trivialized, so we obtain a morphism from
$\CO_{\CB_{\Sp(6)}}$ to $\fL_1\otimes\fL_2$.

Furthermore, the value of $\varphi_2(W,\theta,\zeta)$ is defined as follows. We consider the preimage
$\CU\subset\CW^\perp\subset\CV$ of $\zeta(\fL^*)\subset\CW^\perp/\CW$ (a rank~3 vector subbundle).
The composition of natural
morphisms $\CU\hookrightarrow\CW^\perp\hookrightarrow\CV\twoheadrightarrow\CV/\CV^3$ gives rise
to a morphism $\det\CW\otimes\fL^*=\det\CU\to\det(\CV/\CV^3)=\fL_1\otimes\fL_2\otimes\fL_3$.
However, $\det\CW$ is trivialized, so we obtain the desired morphism
$\fL^*\to\fL_1\otimes\fL_2\otimes\fL_3$.

Finally, the value of $\varphi_3(W,\theta,\zeta)$ is defined as follows. We have the composition 
$\CW^\perp\hookrightarrow\CV\twoheadrightarrow\CV/\CV^3$ that gives rise to a morphism
$\Lambda^3\CW^\perp\to\Lambda^3(\CV/\CV^3)=\fL_1\otimes\fL_2\otimes\fL_3$.
The filtration $0\subset\CW\subset\CW^\perp$ gives rise to an embedding
$(\CW^\perp/\CW)\otimes\det\CW\hookrightarrow\Lambda^3\CW^\perp$. Composing it with the previous
morphism we obtain a morphism
\[\CW\otimes\det\CW\xrightarrow{\zeta}(\CW^\perp/\CW)\otimes\det\CW\hookrightarrow\Lambda^3\CW^\perp
\to\fL_1\otimes\fL_2\otimes\fL_3.\] Since
$\det\CW$ is trivalized, we get a morphism $\CW\to\fL_1\otimes\fL_2\otimes\fL_3$.
Adding to it the composition \(\CW\xrightarrow{\theta}\CV\twoheadrightarrow\fL_1\), we get a morphism
$\CW\to\fL_1\oplus\fL_1\otimes\fL_2\otimes\fL_3$. Taking its determinant we obtain a morphism
$\det\CW\to\fL_1^{\otimes2}\otimes\fL_2\otimes\fL_3$. Recalling once again that $\det\CW$ is
trivialized, we obtain a morphism from $\CO_{\CB_{\Sp(6)}}$ to $\fL_1^{\otimes2}\otimes\fL_2\otimes\fL_3$.

We consider the formal disc $\widehat{C}_c$ (the completion of $C$ at $c$) and the punctured
formal disc $\widehat{C}{}_c^*$. We choose a local parameter $t$.
Then $\Gr_{\wt{G}}$ is the moduli space of $\wt{G}$-bundles 
on $\widehat{C}_c$ equipped with a section on $\widehat{C}{}_c^*$. Let $V_i$ denote the target
of $\varphi_i,\ 1\leq i\leq 4$. Given a $\wt{G}$-bundle with a section $(\CF,\sigma)\in\Gr_{\wt{G}}$,
we obtain the associated vector bundle $\CV_i^\CF$ on $\widehat{C}_c$ with a section $\sigma_i$ on
$\widehat{C}{}_c^*$. We denote by $\on{ord}(\sigma_i)$ the degree of pole of $\sigma_i$ at $c$,
i.e., the minimal power $a$ of the parameter $t$ such that $t^a\sigma_i$ is regular (and nonvanishing)
at $c$. Clearly, $f_i:=\on{ord}(\sigma_i)$ is a lower semi-continuous $\wt{H}(\bF)$-invariant
function on $\Gr_{\wt{G}}$.

Note that the highest weights $\omega^\svee_i$ of $V_i,\ 1\leq i\leq4$, form the dual basis of
$\{\alpha_i\}$, i.e.\ $\langle\omega^\svee_i,\alpha_j\rangle=\delta_{ij}$, where
\begin{equation}
  \label{omega}
  \omega^\svee_1=\varepsilon_1^*+\delta^*,\
\omega^\svee_2=\varepsilon_1^*+\varepsilon_2^*+\varepsilon_3^*+\delta^*,\
\omega^\svee_3=2\varepsilon_1^*+\varepsilon_2^*+\varepsilon_3^*,\
\omega^\svee_4=\varepsilon_1^*+\varepsilon_2^*.
\end{equation}
So the following lemma completes the proof of the proposition for
$\sO^{\blambda'}\subset\ol\sO{}^\blambda\subset\Gr_{\wt{G}}$.

\begin{lem}\label{3.7.2}
  We have $f_i|_{\BO^\blambda}=\langle\omega^\svee_i,\blambda\rangle$.
\end{lem}

\begin{proof}
  Choose a point in $\BO^\blambda$ and let $M\in \PSp(6,\bF)$ be a representative. Note that we have the
  following description of the values of the functions $f_i$ in terms of matrices:
\begin{enumerate}
    \item $-f_1$ is the minimal degree of the entry in rows $5$ and $6$;
    \item  $-f_2$ is the minimal degree of the $3\times 3$-minors in rows $(3,5,6)$ and $(4,5,6)$;
    \item $-f_3$ is the minimal degree of
      $M_{5,i}\cdot \det M_{(4,5,6),(j,k,l)}+M_{6,i}\cdot \det M_{(3,5,6),(j,k,l)}$, where
      $1\leq i,j,k,l\leq 6$ and $j,k,l$ are pairwise distinct, and $M_{5,i}$ stands for the matrix
      element of $M$ at row 5 and column $i$, while      
      $\det M_{(4,5,6),(j,k,l)}$ stands for the minor with rows $(4,5,6)$ and columns $(j,k,l)$;
    \item $-f_4$ is the minimal degree of the $2\times 2$-minors in rows $5$ and $6$.
\end{enumerate}
In terms of the above description, we have 
\[f_1|_{\BO^\blambda}= \lambda_1+\lambda_4, f_2|_{\BO^\blambda}= \lambda_1+\lambda_2+\lambda_3+\lambda_4, f_3|_{\BO^\blambda}= 2\lambda_1+\lambda_2+\lambda_3, f_4|_{\BO^\blambda}=\lambda_1+\lambda_2.\]

\end{proof}

Finally, to take care of the orbits $\sO^{\blambda'}\subset\ol\sO{}^\blambda\subset\Gr_{\ol{G}}$, we note
that they are obtained from the orbits in $\Gr_{\wt{G}}$ by a translation.
Proposition~\ref{adjacency} is proved.
\end{proof}

\subsubsection{Intersections of $H(\bF)$-orbits and semi-infinite orbits}
\label{intersections of}
Recall the element $g\in\PSp(6)$ introduced in~\S\ref{matrix g}. We denote by
$B=B_{\PGL(2)}\times g B^-_{\PSp(6)} g^{-1}$ the Borel subgroup of $\overline{G}$. We have $B\cap H=\{e\}$.
Let $U=U_{\PGL(2)}\times g U^-_{\PSp(6)} g^{-1}$ be the unipotent radical of $B$. Given a coweight
$\bmu=(-\mu_4, (\mu_1, \mu_2, \mu_3))$ of $\PGL(2)\times \PSp(6)$, we denote by $\sS^{\bmu}$ the
$U(\bF)$-orbit of the point $(t^{-\mu_4}, g t^{(\mu_1, \mu_2,\mu_3)})\in \Gr_{\overline{G}}$.
\begin{prop}\label{so intersect}
    If $\sO^{\blambda}\cap \sS^{\bmu}\neq \emptyset$, then $\bmu\leq \blambda$.
\end{prop}
\begin{proof}
    The non-empty intersection condition can be written explicitly as 
    \begin{equation*}
        \begin{split}
    (\PGL(2,\bF)\ltimes \on{U}(\bF))\cdot (t^{-\lambda_4}, g\cdot t^{(\lambda_1,\lambda_2,\lambda_3)}) \overline{G}(\bO)\cap \\ \cap (U_{\PGL(2)}(\bF)\times g U^-_{\PSp(6)}(\bF) g^{-1})\cdot (t^{-\mu_4}, g t^{(\mu_1, \mu_2,\mu_3)}) \overline{G}(\bO)\neq \emptyset.        
        \end{split}
    \end{equation*}
Here, the intersection is taken inside $\overline{G}(\bF)$.
    
    The above non-empty intersection condition is equivalent to 
    \begin{equation*}
        \begin{split}
            (\PGL(2,\bO)\ltimes \on{U}(\bF))\cdot t^{\lambda_4}g\cdot t^{(\lambda_1,\lambda_2,\lambda_3)}\cdot \PSp(6,\bO)\cap \\
            \cap t^{\mu_4} \cdot U_{\PGL(2)}(\bF)\cdot g\cdot U_{\PSp(6)}^-(\bF)\cdot t^{(\mu_1, \mu_2,\mu_3)}\neq \emptyset.
        \end{split}
    \end{equation*}
Here, the intersection is taken inside $\PSp(6,\bF)$.
    
Given any point $M\in t^{\mu_4} \cdot U_{\PGL(2)}(\bF)\cdot g\cdot U_{\PSp(6)}^-(\bF)\cdot
t^{(\mu_1, \mu_2,\mu_3)}$, if we regard $M$ as a $6\times 6$-matrix with coefficients in $\bF$, then 
\begin{itemize}
    \item $\deg M_{6,6}= {-\mu_1-\mu_4}$;
    \item $\deg \det M_{(5,6),(5,6)}= {-\mu_1-\mu_2}$;
    \item $\deg (M_{5,6}\cdot \det M_{(4,5,6),(4,5,6)}+ M_{6,6}\cdot \det M_{(3,5,6),(4,5,6)})=
      {-2\mu_1-\mu_2-\mu_3}$;
    \item $\deg \det M_{(4,5,6),(4,5,6)}= {-\mu_1-\mu_2-\mu_3-\mu_4}$.
\end{itemize}

However, if $M$ belongs to $ \PGL(2,\bO)\ltimes \on{U}(\bF)\cdot t^{\lambda_4}g\cdot t^{(\lambda_1,\lambda_2,\lambda_3)}\cdot \PSp(6,\bO)$, then according to the descriptions (1)-(4) in the proof of Lemma \ref{3.7.2}, we have
\begin{equation*}
    \begin{split}
        -\mu_1-\mu_4\geq -\lambda_1-\lambda_4,& -\mu_1-\mu_2\geq -\lambda_1-\lambda_2,\\
        -2\mu_1-\mu_2-\mu_3\geq -2\lambda_1-\lambda_2-\lambda_3,& -\mu_1-\mu_2-\mu_3-\mu_4\geq -\lambda_1-\lambda_2-\lambda_3-\lambda_4,
    \end{split}
\end{equation*}
which is equivalent to $\bmu\leq \blambda$.
\end{proof}

\subsection{Global Gaiotto category}
For technical reasons, a global definition is required for the twisted Gaiotto category. The desired
category is defined by considering the generic Hecke equivariant sheaves on the global model. 

\subsubsection{Global model}\label{global model}
Recall the notation $\overline{G}=\PGL(2)\times \PSp(6)$ and $H=\PGL(2)\ltimes\on{U}$. Let $\cX$ be
the affine closure $(H\backslash \overline{G})^{\aff}$, and let $\overset{\circ}{\cX}$ be the
homogeneous spherical variety $H\backslash \overline{G}$. The global model
$\Maps_{\on{gen}}(C, \cX/\overline{G}\supset \overset{\circ}{\cX}/\overline{G})$~\cite[\S 3.1]{sw}
is defined as the moduli stack of maps $C\to\cX/\overline{G}$ generically landing into
$\overset{\circ}{\cX}/\overline{G}$. It is known~\cite[Proposition 3.1.2]{sw} that
$\Maps_{\on{gen}}(C, \cX/\overline{G}\supset \overset{\circ}{\cX}/\overline{G})$ is an algebraic stack
of finite type. 

In order to have a canonically defined character $\chi$, we need to consider the
\textit{renormalized} global model.

The renormalized global model $\scrM$ is defined as the moduli stack of pairs
$(\CP_{\overline{G}},\sigma)$, where $\CP_{\overline{G}}$ is a $\overline{G}$-bundle on $C$ and
$\sigma\colon C\longrightarrow \on{T}'\backslash \cX\overset{\overline{G}}{\times}\CP_{\overline{G}}$
is a section such that the induced $\on{T}'$-bundle on $C$ is $\bomega_C^\rho$, and $\sigma$
generically lands into $\on{T}'\backslash \overset{\circ}{\cX}\overset{\overline{G}}{\times}\CP_{\overline{G}}$.

More precisely, using \cite[Theorem 1.1.2, Section 1.1.3]{bg}, we can write $\scrM$ as follows.
\begin{defn}
  The algebraic stack $\scrM$ classifies the data
  $(\CP_{\overline{G}}, \{\varkappa_V\}_{V\in \Rep(\on{L})})$, where
  $\CP_{\overline{G}}:=(\CP_{\PGL(2)}, \CP_{\PSp(6)})$ is a $\overline{G}$-bundle, and
    $\varkappa_V$ is an injective map of coherent sheaves
    \begin{equation}\label{3.6.1}
      \varkappa_V\colon ({\CP_{\PGL(2)}}\overset{\PGL(2)}{\times}\on{L}\otimes \bomega_{C\times S}^{\rho})\overset{\on{L}}{\times} V \longrightarrow
      \CP_{\PSp(6)}\stackrel{\PSp(6)}{\times}\Ind_{\on{P}}^{\PSp(6)}V,\
      \textnormal{for any } V\in \Rep(\on{L}),
    \end{equation}
    satisfying the Pl\"ucker relations. Here, ${\CP_{\PGL(2)}}\overset{\PGL(2)}{\times}\on{L}\otimes \bomega_{C\times S}^{\rho}$ denotes the $\on{L}$-bundle induced from the $\on{L}\times \on{T}'$-bundle $({\CP_{\PGL(2)}}\overset{\PGL(2)}{\times}\on{L}, \bomega_{C\times S}^{\rho})$ via the multiplication map $\on{L}\times \on{T}'\to \on{L}$, and $\Ind_{\on{P}}^{\PSp(6)}V$ denotes the
    $\PSp(6)$-module induced from the $\on{P}$-module $V$ (we regard $V$ as a $\on{P}$-module via
    the natural projection $\on{P}\to\on{L}$).     
\end{defn}

Similarly, we define the global model $\scrM_{\infty\cdot c}$ with a marked point $c$ that
classifies $(\CP_{\overline{G}},\sigma)$, where $\CP_{\overline{G}}$ is a $\overline{G}$-bundle on $C$ and
$\sigma\colon C\setminus c\longrightarrow \on{T}'\backslash \cX\overset{\overline{G}}{\times}\CP_{\overline{G}}$
is a section such that the induced $\on{T}'$-bundle on $C\setminus c$ is $\bomega_C^\rho$, and $\sigma$
generically lands into $\on{T}'\backslash \overset{\circ}{\cX}\overset{\overline{G}}{\times}\CP_{\overline{G}}$.
It is known that $\scrM_{\infty\cdot c}$ is an ind-algebraic stack of ind-locally finite type,
see~\cite[Lemma 3.5.2]{bft1}.

Recall the line bundle $\CP_{\det}$ defined in the Section \ref{renormalized line bundle}. We consider the category of $\CP^\kappa_{\det}$-twisted $D$-modules on $\scrM_{\infty\cdot c}$, denoted as $\on{D-mod}_{\CP^\kappa_{\det}}(\scrM_{\infty\cdot c})$.

\subsubsection{Stratification}
Given a geometric point of $\scrM_{\infty\cdot c}$, taking the restriction of $(\cP_{\oG},\sigma)$ to the formal disc $\cD_c$, we obtain a section $\sigma\colon \overset{\circ}{\cD}_c\rightarrow \on{T}'\backslash \overset{\circ}{\cX} \overset{\oG}{\times}\cP_{\oG}$ and the induced $\on{T}'$-bundle on $\overset{\circ}{\cD}_c$ is $\bomega_{\overset{\circ}{\cD}_c}^\rho$. 

Since any principal $\oG$-bundle on $\cD$ is (\'{e}tale) isomorphic to $\cP_{\oG}^\bomega$, the
above map gives rise to a section
$\overset{\circ}{\cD}_c\rightarrow \on{T}'\backslash \overset{\circ}{\cX} \overset{\oG}{\times}\cP_{\oG}^\bomega$.
Furthermore, two isomorphisms will differ by a multiplication by $^\bomega\!\oG(\bO)$, so we
obtain a point in \[|((\on{T}'\backslash (H\backslash G))(\bF)\underset{(\on{T}'\backslash \pt)(\bF)}{\times}\pt)/ ^{\bomega}\oG(\bO)|=| ^{\bomega}\!H(\bF)\backslash ^{\bomega}\Gr_{\oG}|.\]
Here, the map $\pt\rightarrow (\on{T}'\backslash \pt)(\bF)$ is given by $\bomega_{\overset{\circ}{\cD}_c}^\rho$.

For any $\blambda\in |^{\bomega}\!H(\bF)\backslash ^{\bomega}\Gr_{\oG}|$, we denote by
$\scrM_{=\blambda\cdot c}$ the locally closed substack of $\scrM_{\infty\cdot c}$ such that the above evaluation at $c$ is $\blambda$ and at other points is $0$. Note that if $\blambda=0$, then $\scrM_{=0\cdot c}$ is the $\bomega$-renormalized version of $\Bun_H$.

\subsubsection{Generic Hecke equivariant sheaves}\label{3.5.2}
Now, we mimic the definition of the generic Hecke equivariant sheaves in~\cite{gn,g} to give the definition of the twisted Gaiotto category. 

For any scheme $S$, using the same definition as in~\S\ref{global model}, we can also define the
$S$-stack $\scrM_{\infty\cdot c, S}$.
Now, given an $S$-point $\underline{x}$ in the Ran space, i.e.\ a non-empty set of
$\Maps(S,C\setminus c)$, one can
single out an open substack $\scrM_{\underline{x}, \infty\cdot c, S}\subset \scrM_{\infty\cdot c, S}$, where we
impose the condition that the maps $\{\varkappa_V\}$ are injective maps of vector bundles near
$\underline{x}$. Since $\varkappa_V$ is an injective vector bundle map near $\underline{x}$, the collection of maps
$\varkappa_V$ gives rise to an $^\bomega\!H$-bundle on the formal neighborhood $\cD_{\underline{x}}$ of $\underline{x}$,
i.e.\ a $\on{T}'\ltimes H$-bundle $\cP_{\on{T}'\ltimes H}$ with an identification of the induced $\on{T}'$-bundle with
$\bomega_{C\times S}^\rho$.

Let $\Hecke_{\underline{x}}$ be the stack that classifies the data
$(\CP_{\overline{G}}, \sigma, \cP_H, \varsigma)$, where $(\CP_{\overline{G}}, \sigma)$ belongs to
$\scrM_{\underline{x}, \infty\cdot c,S}$, $\cP_H$ is a principal
$^\bomega\!H$-bundle on the formal neighborhood of $\underline{x}$, and $\varsigma$ is an isomorphism
in the punctured formal neighborhood of $\underline{x}$ of
$\cP_H$ and the $^\bomega\!H$-bundle constructed from $(\CP_{\overline{G}}, \sigma)$ (since $\sigma$
generically lands in $\on{T}'\backslash \overset{\circ}{\cX}\overset{\overline{G}}{\times} \CP_{\overline{G}}$,
it induces a generic ${}^\bomega H$-reduction of $\CP_{\overline{G}}$ on $C$).
Given such a point, we can form another point $(\CP'_{\overline{G}},\sigma')$
in $\scrM_{\underline{x}, \infty\cdot c,S}$ by gluing $\sigma|_{(C\setminus c)\times S\setminus \underline{x}}$ and $\cP_H$ along the
punctured formal neighborhood of $\underline{x}$. 

We consider the following diagram,
\[\xymatrix{
&\Hecke_{\underline{x}}\ar[ld]_{\overset{\gets}{h}}\ar[rd]^{\overset{\to}{h}}&\\
\scrM_{\underline{x}, \infty\cdot c,S}&&\scrM_{\underline{x}, \infty\cdot c,S},
}\]
where $\overset{\gets}{h}$ sends $(\CP_{\overline{G}}, \sigma, \cP_H, \varsigma)$ to $(\CP_{\overline{G}}, \sigma)$ and ${\overset{\to}{h}}$ sends $(\CP_{\overline{G}}, \sigma, \cP_H, \varsigma)$ to $(\cP'_G, \sigma')$. 

For $\kappa\in\BC$, we define the category $\cC^{\on{glob}}_{\kappa,\underline{x}}$ as the category of
$\CP^\kappa_{\det}$-twisted $D$-modules $\cF$ on $\scrM_{\underline{x}, \infty\cdot c,S}$ equipped with an
isomorphism of $\overset{\gets}{h}{}^!(\cF)$ with the tensor product of $\overset{\to}{h}{}^!(\cF)$
and the local system $\chi$. 

Similarly to~\cite[\S 4.3.3]{ty}, if $\underline{x}\subset \underline{x}'\in \Ran_{C\setminus c}(S)$, then there exists a forgetful functor
\begin{equation}\label{trans 1}
\cC^{\on{glob}}_{\kappa,\underline{x}'}\longrightarrow \cC^{\on{glob}}_{\kappa,\underline{x}}.
\end{equation}
Additionally, for any $\underline{x}\in \Ran_{C\setminus c}(S_1)$ and $S_2\to S_1$, there exists a pullback functor
\begin{equation}\label{trans 2}
\cC^{\on{glob}}_{\kappa,\underline{x}}\longrightarrow \cC^{\on{glob}}_{\kappa,\underline{x}\underset{S_1}{\times}S_2}.
\end{equation}

\begin{defn}
 We define\[\cC_\kappa^{\on{glob}}:= \lim_{\underline{x}\in \Ran_{C\setminus c}(S),\ S\in \Sch} \cC^{\on{glob}}_{\kappa,\underline{x}}.\]
    Here, the transition functors are the  functors \eqref{trans 1} and \eqref{trans 2} defined above.
\end{defn}

Recall the renormalized affine Grassmannian $^\bomega\Gr_{\overline{G}}$ defined in~\S\ref{gaiotto category}. There is a projection map
\begin{equation}\label{3.5.4}
    \pi\colon{} ^\bomega\Gr_{\overline{G}}\longrightarrow \scrM_{\infty\cdot c},
\end{equation}
that sends $(\CP_{\overline{G}}, \varsigma)$ to $(\CP_{\overline{G}},\{\varkappa_V\}_{V\in \Rep(\on{L})})$,
where the map $\varkappa_V$ is defined as the composed map of coherent sheaves on $C\setminus c$
\begin{multline}
  \varkappa_V\colon ({\CP_{\PGL(2)}}\overset{\PGL(2)}{\times}\on{L}\otimes \bomega_{C\times S}^{\rho})\overset{\on{L}}{\times} V\to
  ({{\CP^{\on{triv}}_{\PGL(2)}}\overset{\PGL(2)}{\times}\on{L}\otimes \bomega_{C\times S}^{\rho}})\overset{\on{L}}{\times} V\\
  \to\CP^\bomega_{\PSp(6)}\stackrel{\PSp(6)}{\times}\Ind_{\on{P}}^{\PSp(6)}V\to
  \CP_{\PSp(6)}\stackrel{\PSp(6)}{\times}\Ind_{\on{P}}^{\PSp(6)}V.
\end{multline}
Here, the $\PSp(6)$-bundle $\CP^\bomega_{\PSp(6)}$ is defined as the  bundle induced from the $\on{T}'$-bundle $\bomega_{C\times S}^\rho$.



With the same proof as \cite[Theorem 4.4.1]{ty}, we have
\begin{thm}\label{thm 3.8.3}
    The $!$-pullback functor $\pi^!$ induces an equivalence of categories
    \[\pi^!\colon \cC_\kappa^{\on{glob}}\longrightarrow
   \cC_{\kappa}^{\loc}.\]
\end{thm}

  The above equivalence also induces a strata-wise equivalence
  $\cC_{\kappa}^{\on{glob}}|_{\scrM_{=\blambda\cdot c}}\simeq  \on{D-mod}_{\CP^\kappa_{\det}}^{^\bomega\!H(\bF),\chi}(\sO_c^\blambda)$. In particular,  for any relevant
  $\blambda\in | ^{\bomega}H(\bF)\backslash ^\bomega \Gr_{\oG}|$, there is a $D$-module in
  $\cC_{\kappa}^{\on{glob},\heartsuit}|_{\scrM_{=\blambda\cdot c}}$ (the heart of the natural $t$-structure)
  that corresponds to the rank 1 local system $\chi$
  on $\sO^\blambda\subset {}^\bomega\Gr_{\overline{G}}$.  We denote by
  $\IC^{\blambda}_{\on{glob}}\in \on{D-mod}_{\CP^\kappa_{\det}}(\scrM_{\infty\cdot c})$ its IC-extension.
  The object $\IC^{\blambda}_{\on{glob}}$ acquires a natural generic Hecke equivariant structure, and the set of all $D$-modules $\IC^{\blambda}_{\on{glob}}$ forms the collection of irreducible objects in $\cC_\kappa^{\on{glob},\heartsuit}$.

  We denote by $\cC_{\kappa}^{\on{glob},\lc}$ the full subcategory of $\cC_{\kappa}^{\on{glob}}$, such that $\cF\in \cC_{\kappa}^{\on{glob},\lc}$ if and only if $\supp(\cF)$ only contains finitely many relevant strata of $\scrM_{\infty\cdot c}$ and the $!$-stalk at any point is finite dimensional. It is the full subcategory generated by $\IC^{\blambda}_{\on{glob}}$ via extensions and shifts, and the local-global equivalence functor restricts to an equivalence $\pi^!\colon \cC_\kappa^{\on{glob},\lc}\rightarrow
   \cC_{\kappa}^{\loc,\lc}.$

\subsection{Zastava space}
To construct the functor from the twisted Gaiotto category to the category of representations (i.e., the category of twisted factorization modules), we need to use the SW (Sakellaridis-Wang) zastava space, cf., \cite[\S 3.3]{sw}. In our case, in order to match the twistings in Section \ref{renormalized line bundle} and have a canonically defined character $\chi$, we need to consider the $\bomega$-renormalized SW zastava space.

To begin, let us recall the definition of the SW zastava space without the marked point.

\begin{defn}
We define the zastava space ${\CY}$ as the fiber product
\[\pt\underset{\Bun_{\on{T}'}(C)}{\times} \Maps_{\on{gen}}(C, \on{T}'\backslash \cX/B \supset \on{T}'\backslash{\cX}'/B).\]
Here, $\cX'$ denotes the open $B$-orbit in $\cX$ and the map $\pt\to \Bun_{\on{T}'}(C)$ is given by $\bomega_C^\rho\in \Bun_{\on{T}'}(C)$.

Similarly, we define ${\CY}_{\infty\cdot c}$ as the fiber product
\[\pt\underset{\Bun_{\on{T}'}(C\setminus c)}{\times} \Maps_{\on{gen}}(C\setminus c, \on{T}'\backslash \cX/B
\supset \on{T}'\backslash{\cX}'/B)\underset{\Bun_B(C\setminus c)}{\times}\Bun_B(C).\]
\end{defn}

The projection $\on{T}'\backslash \cX/B\to \on{T}'\backslash \cX/\overline{G}$ induces a map of mapping stacks
$\CY\rightarrow \scrM$ (resp. $\CY_{\infty\cdot c}\to \scrM_{\infty\cdot c}$). Also, the map $\on{T}'\backslash \cX/B\to \pt/B$ induces the map $\CY\to \Bun_B$ (resp. $\CY_{\infty\cdot c}\to \Bun_B$).  In particular, we have maps 
\[\CY\rightarrow \scrM\underset{\Bun_{\overline{G}}}{\times}\Bun_B,\] and \[\CY_{\infty\cdot c}\rightarrow \scrM_{\infty\cdot c}\underset{\Bun_{\overline{G}}}{\times}\Bun_B.\] It is known that the above maps are open embedding and the images are the open locus whose generic ${}^\bomega\! H$-reduction and the $B$-reduction are transverse.

\subsubsection{}
We also need the compactified zastava spaces $\bar{\CY}$ and $\bar{\CY}_{\infty\cdot c}$. Let us denote
by $\overline{\Bun}_B$ the Drinfeld compactification of $\Bun_B$, i.e.\ the closure of $\Bun_B$ in
$\Maps_{\on{gen}}(C, \overline{G}\backslash (\overline{G}/U)^\aff/
\overline{T}\supset \overline{G}\backslash\overline{G}/B)$. Here $B\subset\ol{G}$ was introduced
in~\S\ref{intersections of}, and $\overline{T}$ is the maximal torus of $B$ equal to the product
of the diagonal torus of $\PGL(2)$ and $g{\on{T}}g^{-1}$: the diagonal torus of $\PSp(6)$ conjugated
by $g$.

\begin{defn}
  We define the compactified zastava space $\bar{\CY}$ as the intersection of
  $\pt\underset{\Bun_{\on{T}'}(C)}{\times} \Maps_{\on{gen}}(C, \on{T}'\backslash
  \cX\overset{\overline{G}}{\times}(\overline{G}/U)^\aff/ \overline{T} \supset \on{T}'\backslash \pt)$ and
  $\scrM\underset{\Bun_{{\overline{G}}}}{\times}\overline{\Bun}_B$ inside $\pt\underset{\Bun_{\on{T}'}(C)}{\times} \Maps_{\on{gen}}(C, \on{T}'\backslash \cX\overset{\overline{G}}{\times}(\overline{G}/U)^\aff/ \overline{T} \supset \on{T}'\backslash \cX/B)$. In other words, it is the open subset of $\scrM\underset{\Bun_{{\overline{G}}}}{\times}\overline{\Bun}_B$ where the generic ${}^\bomega H$-reduction and the generic $B$-reduction are transverse.

We define the compactified zastava space with poles $\bar{\CY}_{\infty\cdot c}$~as~the~\nohyphens{intersection} of
$\pt\underset{\Bun_{\on{T}'}(C\setminus c)}{\times} \Maps_{\on{gen}}(C\setminus c,
\on{T}'\backslash \cX\overset{\overline{G}}{\times}(\overline{G}/U)^\aff/ \overline{T} \supset
\on{T}'\backslash \pt)\underset{\Maps_{\on{gen}}(C\setminus c, \overline{G}\backslash
(\overline{G}/U)^\aff/ \overline{T})}{\times} \Maps_{\on{gen}}(C, \overline{G}\backslash
(\overline{G}/U)^\aff/ \overline{T}\supset \pt/B)$ and
$\scrM_{\infty\cdot c}\underset{\Bun_{{\overline{G}}}}{\times}\overline{\Bun}_B$ inside
$\pt\underset{\Bun_{\on{T}'}(C\setminus c)}{\times} \Maps_{\on{gen}}(C\setminus c, \on{T}'\backslash
\cX\overset{\overline{G}}{\times}(\overline{G}/U)^\aff/ \overline{T} \supset \on{T}'\backslash
\cX/B)\underset{\Maps_{\on{gen}}(C\setminus c, \overline{G}\backslash (\overline{G}/U)^\aff/ \overline{T})}{\times}
\Maps_{\on{gen}}(C, \overline{G}\backslash (\overline{G}/U)^\aff/ \overline{T}\supset \pt/B)$.
It is the open subset of $\scrM_{\infty\cdot c}\underset{\Bun_{{\overline{G}}}}{\times}\overline{\Bun}_B$
where the generic ${}^\bomega H$-reduction and the generic $B$-reduction are transverse.
\end{defn}

\subsubsection{}
Let $H\backslash\!\!\backslash \overline{G}\dslash U$ be $\Spec \BC[\overline{G}]^{H\times U}$.
Consider the projection
\begin{equation*}
\begin{split}
    \cX\overset{\overline{G}}{\times}(\overline{G}/U)^\aff\simeq \overline{G}\backslash (\cX{\times}(\overline{G}/U)^\aff)\simeq  \overline{G}\backslash ((\overline{G}/H)^\aff){\times}(\overline{G}/U)^\aff)\\ \longrightarrow (\overline{G}\times \overline{G})\dslash(\overline{G}\times H\times U)\simeq H\backslash\!\!\backslash \overline{G}\dslash U.
\end{split}
\end{equation*}
Applying $\pt\underset{\Bun_{\on{T}'}(C)}{\times}\Maps_{\on{gen}}(C, \on{T}'\backslash - / \overline{T}\supset \on{T}'\backslash \pt)$ and \[\pt\underset{\Bun_{\on{T}'}(C\setminus\!c)}{\times}\Maps_{\on{gen}}(C\setminus c, \on{T}'\backslash - / \overline{T}\supset \on{T}'\backslash \pt)\underset{\Bun_{\overline{T}}(C\setminus c)}{\times}\Bun_{\overline{T}}(C)\] to the above projection we obtain morphisms
\begin{equation}\label{eq 3.9.1}
    \begin{split}
        \CY&\longrightarrow \pt\underset{\Bun_{\on{T}'}(C)}{\times}\Maps_{\on{gen}}(C, \on{T}'\backslash (H\backslash\!\!\backslash \overline{G}\dslash U)/ \overline{T}\supset \on{T}'\backslash \pt),\\
        \CY_{\infty\cdot c}\longrightarrow \pt&\underset{\Bun_{\on{T}'}(C\setminus c)}{\times}\Maps_{\on{gen}}(C\setminus c, \on{T}'\backslash (H\backslash\!\!\backslash \overline{G}\dslash U)/ \overline{T}\supset \on{T}'\backslash \pt)\underset{\Bun_{\overline{T}}(C\setminus c)}{\times}\Bun_{\overline{T}}(C).
    \end{split}
\end{equation}

Note that any point in $\pt\underset{\Bun_{\on{T}'}(C)}{\times}\Maps_{\on{gen}}(C,
\on{T}'\backslash (H\backslash\!\!\backslash \overline{G}\dslash U)/ \overline{T}\supset \on{T}'\backslash \pt)$
determines a $\overline{T}$-bundle $\CP_{\overline{T}}$ on $C$ and a  rational map
$\CP_{\overline{T}}^\bomega:=\bomega_C^\rho\overset{\on{T}'}{\times}{\overline{T}}\dasharrow \CP_{\overline{T}}$.
Taking degeneracy locus of this map
induces a morphism (actually, an isomorphism)
\begin{equation}\label{eq 3.9.2'}
    \begin{split}
      \pt\underset{\Bun_{\on{T}'}(C)}{\times}\Maps_{\on{gen}}(C, \on{T}'\backslash (H\backslash\!\!\backslash \overline{G}\dslash U)/ \overline{T}\supset \on{T}'\backslash \pt)\longrightarrow
      \bigsqcup_{\alpha\in Q_{\on{\on{neg}}}}C^\alpha
    \end{split}
\end{equation}
(notation of~\S\ref{config spaces}). 

Similarly, we have 
\begin{equation}\label{eq 3.9.2''}
\pt\underset{\Bun_{\on{T}'}(C\setminus c)}{\times}\Maps_{\on{gen}}(C\setminus c, \on{T}'\backslash (H\backslash\!\!\backslash \overline{G}\dslash U)/ \overline{T}\supset \on{T}'\backslash \pt)\underset{\Bun_{\overline{T}}(C\setminus c)}{\times}\Bun_{\overline{T}}(C)\longrightarrow \bigsqcup_{{\bmu}\in Q}C^{\bmu}_{\infty\cdot c}.
\end{equation}

We denote by $\fv\colon \bar{\CY}\rightarrow C^\bullet$ (resp. $\fv\colon \bar{\CY}_{\infty\cdot c}\rightarrow C^\bullet_{\infty\cdot c}$) the composition of the maps \eqref{eq 3.9.1} and \eqref{eq 3.9.2'} (resp. \eqref{eq 3.9.2''}), and denote by $\fp$ (resp. $\fq$) the projection $\bar{\CY}_{\infty\cdot c}\rightarrow \scrM_{\infty\cdot c}$ (resp. $\bar{\CY}_{\infty\cdot c}\rightarrow \overline{\Bun}_B$). With some abuse of notations, we keep the same notation for the restrictions of the above maps.

It is obvious that the map $\fv$ is compatible with the grading of the zastava space (given by the
degree of $\CP_{\overline{T}}$) and the grading of the configuration space (given by the total degree of $D$). Furthermore,  the map $\fv$ is factorizable, cf.~\cite[Proposition 3.4.1]{sw}. Namely, there are natural isomorphisms
\begin{equation}\label{eq 3.9.4}
    \begin{split}
        (\CY\times \CY) \underset{C^\bullet\times C^\bullet}{\times} (C^\bullet\times C^\bullet)_{\on{disj}}&\simeq \CY \underset{C^\bullet}{\times} (C^\bullet\times C^\bullet)_{\on{disj}},\\
         (\bar{\CY}\times \bar{\CY}) \underset{C^\bullet\times C^\bullet}{\times} (C^\bullet\times C^\bullet)_{\on{disj}}&\simeq \bar{\CY} \underset{C^\bullet}{\times} (C^\bullet\times C^\bullet)_{\on{disj}},\\
         (\bar{\CY}\times \bar{\CY}_{\infty\cdot c}) \underset{C^\bullet\times C^\bullet_{\infty\cdot c}}{\times} (C^\bullet\times C^\bullet_{\infty\cdot c})_{\on{disj}}&\simeq \bar{\CY}_{\infty\cdot c} \underset{C^\bullet_{\infty\cdot c}}{\times} (C^\bullet\times C^\bullet_{\infty\cdot c})_{\on{disj}}.
    \end{split}
\end{equation}

The fiber of the zastava space ${\CY}$ over $D\in C^\bullet$ and ${\CY}_{\infty\cdot c}$ over
$D\in C^\bullet_{\infty\cdot c}$ is given by the product of the intersections of
${}^\bomega\! H(\bF)$-orbit and ${}^\bomega\! U(\bF)$-orbit in ${}^\bomega\Gr_{\overline{G}}$. To be more precise, if $D=\bmu_c\cdot c-\sum_{x\in C\setminus c} \alpha_x\cdot x$, we have
\begin{equation}
    \begin{split}
        \CY_{\infty\cdot c}|_{D}\simeq \CY_{\infty\cdot c}|_{{\bmu_c}\cdot c}\times \prod_{x} \CY|_{-\alpha_x\cdot x}\simeq \sS^{{\bmu_c}}\times \prod_x (\sO^{0}\cap \sS^{-\alpha_x}).
    \end{split}
\end{equation}
Here, the first isomorphism is given by the factorization property \eqref{eq 3.9.4} of the zastava
space, and the second isomorphism is given by \cite[Lemma 4.3.2]{sw}.

\subsection{Dimension estimate of $\CY_{\infty\cdot c}|_{D}$}
Let us denote by ${\CY}^{\bmu}_{=\blambda\cdot c}$ the preimage of $\scrM_{= \blambda\cdot c}\subset \scrM_{\infty\cdot c}$ in $\CY^{\bmu}_{\infty\cdot c}$. We denote by ${\oCY}{}^{\bmu}_{= \blambda\cdot c}$ the preimage of $\overset{\circ}{C}{}^{\bmu}_{= \blambda\cdot c}$ in ${\CY}^{\bmu}_{=\blambda\cdot c}$. Here, $\overset{\circ}{C}{}^{\bmu}_{= \blambda\cdot c}$ is the open subscheme of $C^{\bmu}_{\leq \blambda\cdot c}$ consisting of the divisors $D=\bmu_c\cdot c- \sum_{x\in C\setminus c} \alpha_x\cdot x$, whose coefficient $\bmu_c$ at $c$ is exactly $\blambda$ and the coefficient $-\alpha_x$ at any other point belongs to $\{-\alpha_1,-\alpha_2,-\alpha_3,-\alpha_4\}$. 

The zastava space $\overset{\circ}{\CY}{}^{\bmu}_{= \blambda\cdot c}$ is flat over $\overset{\circ}{C}{}^{\bmu}_{= \blambda\cdot c}$, and its fiber over $D=\blambda\cdot c-\sum \alpha_x\cdot x$ is isomorphic to the product of the fibers ${\CY}^{{\bmu}}_{= \blambda\cdot c}|_{\blambda\cdot c}\times \prod_{x} \CY^{-\alpha_x}|_{-\alpha_x\cdot x}\simeq (\sO^{{\blambda}}\cap \sS^{{\blambda}})\times \prod_x (\sO^{0}\cap \sS^{-\alpha_x})$. Note that if $\alpha_x=\alpha_4$, then $\sO^{0}\cap \sS^{-\alpha_x}\simeq \BA^1$, and $\sO^{0}\cap \sS^{-\alpha_x}\simeq \pt$ otherwise. So, if we assume that ${\blambda}=\bmu+\sum_{i=1}^4 a_i \alpha_i$, then $\dim \overset{\circ}{\CY}{}^{\bmu}_{= \blambda\cdot c}=a_1+a_2+a_3+2\cdot a_4$.

\begin{prop}\label{prop 3.10.1}
    Let $\sY$ be an irreducible component of $\sO^{{\blambda}}\cap \sS^{\bmu}$. We have
    \begin{equation}\label{eq 3.9.2}
        2\dim \sY\leq \dim {\oCY}{}^{\bmu}_{= \blambda\cdot c}.
    \end{equation}

Furthermore, if the above equality holds, then $\blambda-\bmu$ belongs to the positive coroot monoid
of $\PGL(2)\times \PSp(6)$. If $\dim \overset{\circ}{\CY}{}^{\bmu}_{= \blambda\cdot c}- 2\dim \sY=1$, then
$\blambda-\bmu$ is of the form $\alpha+\alpha'$, where $\alpha$ belongs to the positive coroot monoid of
$\PGL(2)\times \PSp(6)$ and $\alpha'$ belongs to $\{\alpha_1,\alpha_2,\alpha_3\}$.
\end{prop}
\begin{proof}
    According to Proposition \ref{so intersect}, we have 
    \begin{equation*}
        \begin{split}
            \dim \sY\leq n+\sum_{i=1}^{3} m_i,
        \end{split}
    \end{equation*}
    for any $\bmu\leq_{\overline{G}} \bmu'\leq \blambda$ and $\bmu'= \bmu+ n\delta+
    m_1(\varepsilon_1-\varepsilon_2)+m_2(\varepsilon_2-\varepsilon_3)+m_3\varepsilon_3$.

    If we assume $\blambda=\bmu+\sum_{i=1}^4 a_i \alpha_i$, then the condition $\bmu'\leq \blambda$ implies $n+m_1\leq a_1$, $n+m_3\leq a_2$, $2m_2\leq 2a_4$, and $m_1+m_3\leq a_3$. Taking their sum, we obtain the desired inequality \eqref{eq 3.9.2}: $\dim \oCY{}^{\bmu}_{= \blambda\cdot c}=a_1+a_2+a_3+2\cdot a_4 \geq 2 (n+m_1+m_2+m_3)\geq 2\dim \sY$.
\end{proof}

\begin{rem}
   The closure of $\oCY{}^{\bmu}_{= \blambda\cdot c}$ is a top dimensional irreducible component of ${\CY}^{\bmu}_{= \blambda\cdot c}$.   
\end{rem}

\subsubsection{}\label{section 3.10.1}
For any relevant orbit $\sO^\blambda\subset 
 {}^{\bomega}\Gr_{\overline{G}}$, there is a (up to a scalar) ${}^{\bomega}\!H(\bF)$-equivariant trivialization of $\CP_{\det}|_{\sO^\blambda}$. Also, there is a unique (up to a scalar) ${}^{\bomega}\!U(\bF)$-equivariant trivialization of $\cP_{\det}|_{\sS^\bmu}$. We denote by $\psi_{\blambda,\bmu}$ the ratio of these two trivializations on $\sO^\blambda\cap \sS^\bmu$, and denote by $\Psi_{\blambda,\bmu}$ the pullback of the shifted Kummer local system on $\BG_m$ corresponding to $q$ along $\psi_{\blambda,\bmu}$, such that its $!$-stalks are concentrated in degree $0$.

The following lemma is similar to~\cite[Lemma 9.1.7]{ty} and can be proved using the same method
as~\cite[Proposition 3.7.2]{bft1} and~\cite[Lemma 5.5]{g}.
\begin{lem}\label{lem 3.7.3}
 Let us  assume that $\blambda-\bmu \notin \{0, \alpha_1, \alpha_2, \alpha_3\}$. Then, the shifted local system $\Psi_{\blambda, \bmu}$ is not trivial on any irreducible component $\sY\subset \sO^\blambda\cap \sS^\bmu$, if $2\dim \sY= \dim \oCY{}^{\bmu}_{= \blambda\cdot c}$ or $2\dim \sY+1= \dim \oCY{}^{\bmu}_{= \blambda\cdot c}$.
\end{lem}

\subsection{Construction of the functor}
Let $\CP_{\det, \Bun_{\overline{T}}}$ be the relative determinant line bundle on $\Bun_{\overline{T}}(C)$ whose fiber over $\CP_{\overline{T}}\in \Bun_{\overline{T}}$ is
\begin{equation*}
    \begin{split}
      \on{det}^{-8}R\Gamma(C, \varepsilon_1^*(\CP_{\overline{T}})\otimes \bomega_C^{\frac{1}{2}}) \otimes
      \on{det}^{-8}R\Gamma(C, \varepsilon_2^*(\CP_{\overline{T}}) \otimes \bomega_C^{\frac{1}{2}}) \\
      \otimes\on{det}^{-8}R\Gamma(C, \varepsilon_3^*(\CP_{\overline{T}})\otimes \bomega_C^{\frac{1}{2}})
      \otimes \on{det}^{24}R\Gamma(C, \delta^*(\CP_{\overline{T}})\otimes \bomega_C^{\frac{1}{2}})\\
      \otimes\on{\det}^{16}R\Gamma(C, \bomega_C^{-\frac{1}{2}}) \otimes \on{det}^{-16} R\Gamma (C,
      \bomega_C^{\frac{1}{2}})
    \end{split}
\end{equation*}
(notation of the proof of~Proposition~\ref{adjacency}).
It is defined similarly to the line bundle $\CP$ of~\S\ref{rho} on the configuration space.
Furthermore, its pullback along the Abel-Jacobi map
\begin{equation}
    \begin{split}
     \on{AJ}\colon   C^\bullet_{\infty\cdot c}\longrightarrow \Bun_{\overline{T}}\\
        D\mapsto \CP_{\overline{T}}^\bomega(-D),
    \end{split}
\end{equation}
is canonically isomorphic to the line bundle $\CP$ on $C^\bullet_{\infty\cdot c}$. Let
$\CP_{\det, {\Bun_B}}$ be the pullback of $\CP_{\det, {\Bun_{\overline{T}}}}$.

\begin{lem}
The line bundles
  $\fp^*(\CP_{\det})$ and $\fv^*(\CP)$ are isomorphic on $\CY_{\infty\cdot c}$.
\end{lem}
\begin{proof}
Since the map $\on{AJ}\circ\fv$ is the composition $\CY_{\infty\cdot c}\rightarrow \Bun_B\rightarrow \Bun_{\overline{T}}$,
it is enough to check that the pullback of $\CP_{\det \Bun_{\overline{T}}}$ along $\Bun_B\rightarrow \Bun_{\overline{T}}$ is
isomorphic to the pullback of $\CP_{\det}$ along $\Bun_B\rightarrow \Bun_{\overline{G}}$. The latter claim
can be obtained as follows: if a $\overline{G}$-bundle admits a $B$-reduction,
then 
    \[\on{det}^{24}R\Gamma(C, \delta^*(\CP_{\overline{T}})\otimes \bomega_C^{\frac{1}{2}}) \otimes\on{\det}^{-24}R\Gamma(C, \bomega_C^{\frac{1}{2}})\simeq \] 
    \[   \on{det}^{-12}(C,\CU)\otimes \on{det}^{6}(C, \det \CU\otimes \bomega_C^{-\frac{1}{2}})\otimes \on{det}^{24}R\Gamma(C, \CO)\otimes \on{det}^{-6} R\Gamma(C, \bomega_C^{-\frac{1}{2}}),\]
    and
    \[\on{det}^{-8}R\Gamma(C, \epsilon_1^*(\CP_{\overline{T}})\otimes \bomega_C^{\frac{1}{2}}) \otimes \on{det}^{-8}R\Gamma(C, \varepsilon_2^*(\CP_{\overline{T}}) \otimes \bomega_C^{\frac{1}{2}}) \\ \otimes\on{det}^{-8}R\Gamma(C, \varepsilon_3^*(\CP_{\overline{T}})\otimes \bomega_C^{\frac{1}{2}})\]
    \[\otimes\on{\det}^{16}R\Gamma(C, \bomega_C^{-\frac{1}{2}})\otimes\on{\det}^{8}R\Gamma(C, \bomega_C^{\frac{1}{2}})\simeq  \det\!{}^{-4}\!R\Gamma(C,\CV)\otimes\det\!{}^{6}\!R\Gamma(C,\CM\otimes\bomega_C^{-1/2}) \]
    \[ \otimes \on{det}^{16}R\Gamma(C, \CO)\otimes \on{det}^8 R\Gamma(C, \bomega^{-1}_C)\otimes \on{det}^{-6}R\Gamma(C, \bomega_C^{-\frac{1}{2}}).\]
\end{proof}

We define the functor
\begin{equation} \label{functor F}
F\colon \on{D-mod}_{\CP^\kappa_{\det}}(\scrM_{\infty\cdot c})\longrightarrow       \on{D-mod}_{\CP^\kappa}(C^\bullet_{\infty\cdot c}),\quad
\CF\mapsto\fv_*(\fp^\bullet(\CF)),
\end{equation}
where $\fp^\bullet(\CF):=\fp^!(\CF)[-\on{rel.dim} (\Bun_B,\Bun_{\overline{G}})]$ with parity changed according to
$\on{rel.dim} (\Bun_B,\Bun_{\overline{G}})\pmod{2}$.
Here, $\on{rel.dim} (\Bun_B,\Bun_{\overline{G}})$ equals the dimension of the corresponding connected component of $\Bun_B$ minus the dimension of $\Bun_{\overline{G}}$.

The following lemma is proved the same way as~\cite[Theorem 7.3]{g}.
\begin{lem}
  For $\CF\in \cC_{\kappa}^{\on{glob}}$, the object
  $\fv_!(\fp^!(\CF))\in \on{D-mod}_{\CP^\kappa}(C^\bullet_{\infty\cdot c})$ is well-defined, and 
    \[\fv_!(\fp^!(\CF))\iso\fv_*(\fp^!(\CF)).\]
\end{lem}

\subsubsection{}
Let $\IC_{\on{glob}}^0$ be the rank $1$ generic $(^{\bomega}\!H,\chi)$-Hecke equivariant twisted
$D$-module on $\scrM_{\infty\cdot c}$. Up to a shift, it corresponds to the rank $1$ $(^{\bomega}\!H(\bF),\chi)$-equivariant twisted $D$-module on the unital  $^{\bomega}\!H(\bF)$-orbit in $^{\bomega}\Gr_{\overline{G}}$, under the local-global equivalence in Theorem \ref{thm 3.8.3}.

The following lemmas can be proved using the same analysis as \cite[\S 7.2-7.4, Proposition 9.1.1]{ty},

\begin{lem}\label{lem 3.8.3}
    The functor $F$ is $t$-exact, and it commutes with the Verdier duality.
\end{lem}

\begin{lem}\label{lemma 3.8.4}
    The twisted $D$-module $\Omega:= F(\IC^0_{\on{glob}}) \in \on{D-mod}_{\CP^\kappa}(C^\bullet)$ is a factorization algebra, and for any $\CF\in \cC_{\kappa}^{\on{glob}}$, the image $F(\cF)$ admits a canonical factorization module structure with respect to $\Omega$.
\end{lem}

The following theorem is one of the key steps in the proof of Gaiotto equivalence.
\begin{thm}\label{3.8.5}
    There is an isomorphism of factorization algebras
    \begin{equation}\label{eq 3.8.3}
        \Omega\simeq \CI.
    \end{equation}
\end{thm}


\subsection{Vanishing of the sub-bottom cohomology}
In this section, we will prove Theorem \ref{3.8.5}.

First of all, we note that the restrictions of $\Omega$ and $\CI$ to the open locus
$\overset{\circ}{C^\bullet}\subset C^\bullet$ are canonically isomorphic. Indeed, according to the
factorization structure of $\Omega$ and $\CI$, it is enough to check that they are canonically
isomorphic on $C^{-\alpha}$, for $\alpha\in \{\alpha_1,\alpha_2,\alpha_3,\alpha_4\}$. If $\alpha$
is even, i.e., $\alpha=\alpha_4$, then it follows as in~\cite[\S 5.1]{g}; and if $\alpha$ is odd,
i.e.\ $\alpha\in \{\alpha_1,\alpha_2,\alpha_3\}$, then the map from the zastava space $\CY^{-\alpha}$
to $C^{-\alpha}$ is an isomorphism, and $\Omega|_{C^{-\alpha}}$ is the (twisted) $D$-module corresponding
to the constant sheaf.

To establish an isomorphism~\eqref{eq 3.8.3}, we need to extend the isomorphism of
$\Omega|_{\overset{\circ}{C^\bullet}}$ and $\CI|_{\overset{\circ}{C^\bullet}}$ to the whole configuration space
$C^\bullet$. Due to the self-duality of $\Omega$~(Lemma~\ref{lem 3.8.3}) it is enough to prove that
the $!$-restriction of $\Omega$ to any stratum in $C^{-\alpha}\setminus \overset{\circ}{C}{}^{-\alpha}$
is concentrated in degrees no less than $1$, for any $\alpha\in Q_{\on{pos}\setminus 0}$ not equal to
a simple root. Furthermore, according to the factorization property of $\Omega$, we only need to
prove the following proposition.

\begin{prop}\label{prop 3.12.1}
For any $\alpha\in Q_{\on{pos}}\setminus 0$ which is not a simple root, the $!$-stalk of $\Omega$ at
$-\alpha\cdot x\in C^\bullet$ is concentrated in degrees no less than $2$.
\end{prop}
 By the base change theorem, we should prove that
  \begin{equation}\label{eq 3.12.1}
H^i(\CY|_{-\alpha\cdot x}, \fp^!(\IC^0_{\on{glob}})|_{-\alpha\cdot x}[-\on{rel.dim}(\Bun_B,\Bun_G)])=0,       
    \end{equation}
    if $\alpha$ is not simple, and $i\leq 1$.



\subsubsection{}
Note that the $!$-restriction of $\fp^!(\IC^0_{\on{glob}})[-\on{rel.dim}(\Bun_B,\Bun_G)]$ to any irreducible component $\sY$ of $\CY|_{-\alpha\cdot x}\simeq \sO^0\cap \sS^{-\alpha}$ is a lisse D-module with $!$-stalks $\BC[-\dim \oCY{}^{-\alpha}]$ everywhere. In particular, the $!$-restriction of $\fp^!(\IC^0_{\on{glob}})[-\on{rel.dim}(\Bun_B,\Bun_G)]$ to $\sY$ is concentrated in degrees no less than $\dim \overset{\circ}{\CY}{}^{-\alpha}-\dim \sY$, and its cohomology is concentrated in degrees no less than $\dim \oCY{}^{-\alpha}-2\dim \sY$.

To prove the vanishing of the cohomology in degrees no more than $1$, we only need to consider the case  $2\dim \sY= \dim \oCY{}^{-\alpha}$ or $2\dim \sY= \dim \oCY{}^{-\alpha}-1$. 

The latter case is easier. Namely, in this case, the cohomology is concentrated in degrees no less than
$1$, and $H^1$ is non-zero if and only if the $!$-restriction of
$\fp^!(\IC^0_{\on{glob}})[-\on{rel.dim}(\Bun_B,\Bun_G)]$ to all $\sY$ (or, a dense open subset of $\sY$)
is non-constant. However, up to a shift, this D-module is the tensor product of the local system
corresponding to $\chi$ and the tame local system $\Psi_{0, -\alpha}$. If $\alpha$ is not simple,
then $\Psi_{0, -\alpha}$ is non-constant and the local system on $\sY$ corresponding to $\chi$ is
non-tame, so their tensor product is non-constant.

Now, we only need to focus on the case $2\dim \sY=\oCY{}^{-\alpha}$. According
to~Proposition~\ref{prop 3.10.1}, we have $-\alpha= -n\cdot\delta-m_1(\varepsilon_1-\varepsilon_2)-
m_2(\varepsilon_2-\varepsilon_3)-m_3\varepsilon_3$ (negative combination of simple coroots of $\PGL(2)\times \PSp(6)$).
In this case, the cohomology of \eqref{eq 3.12.1} is concentrated in degrees no less than $0$.
Applying the same analysis as in the previous paragraph, we deduce that the $0$-th cohomology
of~\eqref{eq 3.12.1} vanishes.

So, we only need to check $H^1=0$ and we are interested in the irreducible components which are of the critical dimension (i.e., $n+m_1+m_2+m_3$). 

In order to describe the intersection $\sO^0\cap \sS^{-\alpha}\subset \Gr_{\PGL(2)}\times \Gr_{\PSp(6)}$,
we introduce the following notation.

\begin{defn}
  \textup{(1)} $\BS_{\on{U}}^0$ is defined as the unital $\on{U}(\bF)$-orbit in $\Gr_{\PSp(6)}$.

  \textup{(2)} Given ${\gamma}$ a coweight of $\PSp(6)$, we define
  $\BS_-^{{\gamma}}$ as the $U^-_{\PSp(6)}(\bF)$-orbit of $t^{{\gamma}}$ in $\Gr_{\PSp(6)}$. Similarly, we denote by  $\BS^{{\gamma}}$ the $U_{\PSp(6)}(\bF)$-orbit of $t^{{\gamma}}$ in $\Gr_{\PSp(6)}$.

  \textup{(3)} $U^-_{\on{L}}$ is defined as the unipotent radical of the negative Borel subgroup
  of the Levi $\on{L}\subset\PSp(6)$ introduced in~\S\ref{renormalized line bundle}.

  \textup{(4)} $U^\pm\subset\PSp(6)$ is defined as the semi-direct product $U^-_{\on{L}}\ltimes \on{U}$.

  \textup{(5)} Given a coweight $\lambda$ of $\PSp(6)$, we denote by $\fS^\lambda$ the
  $U^\pm(\bF)$-orbit of $t^\lambda\in \Gr_{\PSp(6)}$.
\end{defn}

By definition, the intersection $\sO^0\cap \sS^{-\alpha}\subset \Gr_{\PGL(2)}\times \Gr_{\PSp(6)}$ is isomorphic to the following intersection inside $\Gr_{\PSp(6)}$:
\begin{equation}\label{eq 3.12.2}
g\cdot \begin{pmatrix}
        t^{n}&*\\ 0& t^{-n}
    \end{pmatrix}\cdot \BS_{\on{U}}^0\cap \BS_-^{-{\gamma}}\subset \Gr_{\PSp(6)}.
    \end{equation}
Here, $g$ is the element in $\PSp(6,\bF)$ corresponding to $\Id+E_{3,4}$, introduced
in~\S\ref{matrix g}, and $\gamma=m_1(\varepsilon_1-\varepsilon_2)+
m_2(\varepsilon_2-\varepsilon_3)+m_3\varepsilon_3$. Furthermore, $\begin{pmatrix}
     t^{n}&*\\ 0& t^{-n}
\end{pmatrix}\subset \PGL(2,\bF)$ acts on $\Gr_{\PSp(6)}$ via the group embedding $\PGL(2,\bF)\hookrightarrow \PSp(6,\bF)$ introduced in~\S\ref{renormalized line bundle}.

\subsubsection{Stratification of $\sO^0\cap \sS^{-\alpha}$}
 The following lemma is easy to check.

\begin{lem}
  The intersection of \eqref{eq 3.12.2} and $\fS^\lambda$ is non-empty only if $\lambda$ is of the form $\lambda=-k_1(\varepsilon_1-\varepsilon_2) -k_2\varepsilon_3$.
\end{lem}

For any $k_1, k_2$, we denote by $\BY_{k_1, k_2}$ the intersection of \eqref{eq 3.12.2} and
$\fS^{-k_1(\varepsilon_1-\varepsilon_2)-k_2\varepsilon_3}$.

We claim
\begin{prop}
    The dimension of $\BY_{k_1,k_2}$ is of dimension $n+m_1+m_2+m_3$, i.e., the critical dimension, only if $n=0$ and $(k_1, k_2)=(0,0)$.

    In particular, in this case, $\BY_{k_1,k_2}= \BS_{\on{U}}^0\cap \BS_-^{-\gamma}\subset \Gr_{\PSp(6)}$.
\end{prop}
\begin{proof}
  Note that the stratum $\BY_{k_1,k_2}$ has a natural projection to $\Gr_{\on{L}}$, which is given by
  the composition $\BY_{k_1,k_2}\hookrightarrow g\cdot \begin{pmatrix}
        t^{n}&*\\ 0& t^{-n}
  \end{pmatrix}\cdot \BS_{\on{U}}^0\cap \fS^{-k_1(\varepsilon_1-\varepsilon_2)-k_2\varepsilon_3}\rightarrow
  \Gr_{\on{P}}\to \Gr_{\on{L}}$. 

  We consider the MV-cycles $\sE_{\on{big}}:=\BS^{n(\varepsilon_1-\varepsilon_2)+n\varepsilon_3}\cap
  \BS_-^{-{\gamma}}\subset \Gr_{\PSp(6)}$. The intersection $\BY_{k_1,k_2}$ is a subspace of $\sE_{\on{big}}\cap \fS^{-k_1(\varepsilon_1-\varepsilon_2)-k_2\varepsilon_3}$, and
  the projection $\BY_{k_1,k_2}\rightarrow \Gr_{\on{L}}$ is the restriction of a similarly defined projection 
    \begin{equation}\label{3.9.1}
       \sE_{\on{big}}{\cap \fS^{-k_1(\varepsilon_1-\varepsilon_2)-k_2\varepsilon_3}}\to  \BS^{n(\varepsilon_1-\varepsilon_2)+n\varepsilon_3}\cap \fS^{-k_1(\varepsilon_1-\varepsilon_2)-k_2\varepsilon_3}\rightarrow \Gr_{\on{P}}\rightarrow \Gr_{\on{L}}.
    \end{equation}
  The image $\sB_{\on{big}}$ of the above projection \eqref{3.9.1} is contained in the following MV-cycles 
  \begin{equation}\label{3.9.2}
      \begin{pmatrix}
        t^{n}&*&&\\ 0& t^{-n}&&\\&&t^n&*\\&&0&t^{-n}
    \end{pmatrix}\on{L}(\bF)/\on{L}(\bO)\cap \begin{pmatrix}
        t^{-k_1}&0&&\\ *& t^{k_1}&&\\&&t^{-k_2}&0\\&&*&t^{k_2}
    \end{pmatrix}\on{L}(\bF)/\on{L}(\bO)\subset \Gr_{\on{L}}.
  \end{equation}

  We note that ${\on{T}}(\bO)$ acts on $\sE_{\on{big}}\cap\fS^{-k_1(\varepsilon_1-\varepsilon_2)-k_2\varepsilon_3}$ and on $\Gr_{\on{L}}$. Furthermore,
  the map~\eqref{3.9.1} is ${\on{T}}(\bO)$-invariant. Also, the action of ${\on{T}}(\bO)$ on the
  MV-cycles~\eqref{3.9.2} is transitive. In particular, it implies that:

  (a) if $ \BS^{n(\varepsilon_1-\varepsilon_2)+n\varepsilon_3}\cap \BS_-^{-{\gamma}}\cap
  \fS^{-k_1(\varepsilon_1-\varepsilon_2)-k_2\varepsilon_3}$ is non-empty, then $\sB_{\on{big}}=$\eqref{3.9.2};
  
  (b) the fibers of the projection~\eqref{3.9.1} are isomorphic and their dimensions are no more than
  $\dim \sE_{\on{big}}- \dim \sB_{\on{big}}$. 

  Since $\sE_{\on{big}}\subset   \BS^{n(\varepsilon_1-\varepsilon_2)+n\varepsilon_3}\cap \BS_-^{-{\gamma}}$ and the
  latter is of pure dimension $2\cdot n_1+m_1+m_2+m_3$, we have
  $\dim \sE_{\on{big}}\leq 2\cdot n_1+m_1+m_2+m_3$. Note that $\sB_{\on{big}}$ is of pure dimension
  $2\cdot n_1+k_1+k_2$. Consequently, the dimension of the fibers in~\eqref{3.9.1} is no more than
  $m_1+m_2+m_3-k_1-k_2$. 

  Let $\sB_{\on{sm}}$ be the image obtained by taking the restriction of~\eqref{3.9.1} to $\BY_{k_1,k_2}$.
  We claim that if the dimension of $\sB_{\on{sm}}$ is no less than $n_1+k_1+k_2$, then
  $(k_1,k_2)=(0,0)$ and $n=0$. Once we proved this claim, we will obtain that the dimension of the
  stratum $\BY_{k_1,k_2}$ is strictly smaller than the critical dimension $n+m_1+m_2+m_3$ unless
  $n=0$ and $(k_1,k_2)=(0,0)$ by the dimension estimate of fibers above.

Now, we prove the above claim.

We check a similar statement for $\on{L}'={\SL_2\times \GL_2}$, the case of $\on{L}=\on{L}'/\{\pm1\}$
is similar. We have 
\begin{equation}\label{3.12.5}
    \sB_{\on{sm}}\subset g\cdot \begin{pmatrix}
        t^{n}&f&&\\ 0& t^{-n}&&\\&&t^n&f\\&&0&t^{-n}
    \end{pmatrix}\on{L}'(\bF)/\on{L}'(\bO)\cap \begin{pmatrix}
        t^{-k_1}&0&&\\ *& t^{k_1}&&\\&&t^{-k_2}&0\\&&*&t^{k_2}
    \end{pmatrix}\on{L}'(\bF)/\on{L}'(\bO)\subset \Gr_{\on{L}'}.
\end{equation}

Here, $f$ varies in $\bF$.

 The intersection above is the intersection of the following two spaces: 
\begin{equation}\label{3.9.4}
\begin{pmatrix}
        t^{n}&*\\ 0& t^{-n}
    \end{pmatrix}\SL(2,\bF)/\SL(2,\bO)\cap \begin{pmatrix}
        t^{-k_1}&0\\ *& t^{k_1}
    \end{pmatrix}\SL(2,\bF)/\SL(2,\bO)    
\end{equation}
 and 
 \begin{equation}\label{3.9.4'}
 \begin{pmatrix}
        t^{n}&*\\ 0& t^{-n}
    \end{pmatrix}\SL(2,\bF)/\SL(2,\bO)\cap  \begin{pmatrix}
       1&-1\\ 0& 1
    \end{pmatrix} \begin{pmatrix}
        t^{-k_2}&0\\ *& t^{k_2}
    \end{pmatrix}\SL(2,\bF)/\SL(2,\bO).    
 \end{equation}

 In particular, $\sB_{\on{sm}}$ embeds into the above two spaces.
 The dimension of~\eqref{3.9.4} is $n+k_1$ and the dimension of~\eqref{3.9.4'} is $n+k_2$. It implies
 that the dimension of $\sB_{\on{sm}}$ is no more than $\min\{n+k_1, n+k_2\}$. Without loss of
 generality, we only need to consider the case when $k_1\leq 0$, $k_2\leq 0$, $n\geq -k_1$ and
 $n\geq -k_2$. 

 If $n>0$, we consider the upper triangular representatives of~\eqref{3.9.4}. They are of the form
 $\begin{pmatrix}
        t^{n}&t^{-k_1}\bO\\ 0& t^{-n}
 \end{pmatrix}$. The condition $\sB_{\on{sm}}\neq \emptyset$ implies that
 $\min\{\deg(t^{-k_1}o+t^{-n}), n\}=-k_2$ for a certain $o\in \bO$. Since $n>0$ and $-k_1\geq 0$,
 $\deg(t^{-k_1}o+t^{-n})=-n$ and $\min\{\deg(t^{-k_1}o+t^{-n}), n\}=-n<0$. However, we require that that
 $-k_2\geq 0$, which is a contradiction. 

 On the other hand, if $n=0$, the conditions $n\geq -k_1$ and $n\geq -k_2$ and $k_1, k_2\leq 0$ imply
 that $k_1=k_2=0$.
    \end{proof}

\subsubsection{Vanishing of the sub-bottom cohomology I}
According to the fact that $\BY_{0,0}$ is the only stratum which is possible of the critical dimension
$n+m_1+m_2+m_3$, so in order to prove \eqref{eq 3.12.1}, it remains to prove that the $H^1$
of the restriction to $\BY_{0,0}$ vanishes. And we can assume $n=0$.

Note that it is enough to check that for any irreducible component (or, its open dense subset)
$\sY_{0,0}$ of $\BY_{0,0}$, the first cohomology $H^1$ of the restriction of
$\fp^!(\IC^0_{\on{glob}})|_{-\alpha\cdot x}[-\on{rel.dim}(\Bun_B,\Bun_G)]$ vanishes. We rewrite
$\fp^!(\IC^0_{\on{glob}})|_{-\alpha\cdot x}[-\on{rel.dim}(\Bun_B,\Bun_G)]$ as
$\chi\overset{!}{\otimes} \Psi_{0,-\alpha}[-\dim \overset{\circ}{\CY}{}^{-\alpha}]$, where $\chi$
is the local system corresponding to the character $\chi$,
and $\Psi_{0,-\alpha}$ is the local system corresponds to the ratio of two trivializations of the determinant line bundle, we require that the $!$-stalks of these two local systems are concentrated in degree $0$.


Let us denote by $d:= m_1+m_2+m_3$ the critical dimension, and we assume that $\sY_{0,0}$ is of dimension $d$.
We only need to prove
\begin{equation}\label{3.9.6}
H^1(\sY_{0,0}, \chi\overset{!}{\otimes} \Psi_{0,-\alpha}[-2d])=0    
\end{equation}
We are going to use the method in~\cite{g} to reduce the proof of~\eqref{3.9.6} to the case $m_1=0$ or $m_2=0$.

Consider the projection 
\begin{equation}
\chi_{\on{univ}}\colon \sY_{0,0}\longrightarrow \BG_a^2,    
\end{equation} 
that sends a point in $\sY_{0,0}$ to the residue of the entries $u_{1,3}$ and $u_{2,4}$ of an element in $\on{U}(\bF)g$ (warning: not in $\on{U}(\bF)$!).

We denote by $l\colon \BG_a^2\to \BG_a$ the map which sends $(x, y)$ to $x+y$. Then
\begin{equation}
    \begin{split}
      H^1(\sY_{0,0}, \chi\overset{!}{\otimes} \Psi_{0,-\alpha}[-2d])=
      H^1(\sY_{0,0}, (l\circ \chi_{\on{univ}})^!(\exp)\overset{!}{\otimes} \Psi_{0,-\alpha}[-2d])\\
      = H^1(\BG_a^2, l^!(\exp)\overset{!}{\otimes} \chi_*(\Psi_{0,-\alpha})[-2d]). 
    \end{split}
\end{equation}

So, $H(\sY_{0,0}, \chi\overset{!}{\otimes} \Psi_{0,-\alpha}[-2d])$, shifted by $2$, is the $!$-stalk of
the Fourier transform of $\chi_*(\Psi_{0,-\alpha})[-2d]$ at $l$. We should show that
$\chi_*(\Psi_{0,-\alpha})[-2d]$ is concentrated in perverse degrees no less than~$2$. 

Let $\on{T}''\subset {\on{T}}$ be the subgroup consisting of the diagonal matrices $\diag(z,w,1,1,w^{-1}, z^{-1})$.
The group $\on{T}''(\bO)$ acts on $\sY_{0,0}$ by left multiplication, and acts on $\BG_a^2$ by sending
$(x,y)$ to $(zx, wy)$. Note that $\chi_{\on{univ}}$ is $\on{T}''(\bO)$-equivariant and $\Psi_{0,-\alpha}$ is $\on{T}''(\bO)$-equivariant against a certain character (see~\cite[Lemma 5.5]{g}), so the Fourier
transform of  $\chi_*(\Psi_{0,-\alpha})[-2d]$  is $\on{T}''(\bO)$-equivariant against a certain character.
In particular, it is lisse on the open stratum of $\BG_a^2$. 

If $m_1$ and
$m_2$ are not zero, then  $\chi_*(\Psi_{0,-\alpha})[-2d]$ is the clean extension from its restriction
to the open stratum (removing diagonals). Now, the claim follows from the fact that the dimension
of the fiber of $\chi_{\on{univ}}$ on the open stratum is no more than $d-2$. So, $\chi_*(\Psi_{0,-\alpha})[-2d]$ is
concentrated in degrees no less than $2$.

\subsubsection{Vanishing of sub-bottom cohomology II}
It remains to check the cases $m_1=0$ or $m_2=0$. Note that for any matrix representative of any
point in $\BY_{0,0}$, the minimal degree of its entries in the 5-th row is~$0$. So, if $\BY_{0,0}$
is non-empty, $t^{m_2-m_1}\in \bO$, i.e.\ $m_2\geq m_1$.

Thus we can assume $m_1=0$. In this case, the intersection $\BY_{0,0}$ is isomorphic to the
intersection of the unital $\on{U}'(\bF)$-orbit in $\Gr_{\PSp(4)}$ and the $U_{\PSp(4)}^-(\bF)$-orbit
of $t^{-m_2(\varepsilon_2-\varepsilon_3)-m_3\varepsilon_3}\in \Gr_{\PSp(4)}$. Here, $\on{U}'$ denotes the
intersection of $\PSp(4)$  (the stabilizer of $v_1,v_6$ in $\PSp(6)$) and $\on{U}$ in $\PSp(6)$.

Note that the projection $\Sp(4)\to \PSp(4)$ induces an isomorphism of $\Gr_{\Sp(4)}$ and the neutral
connected component of $\Gr_{\PSp(4)}$. So the intersection $\BY_{0,0}$ is isomorphic to the following
intersection in $\Gr_{\Sp(4)}$
\begin{equation}\label{3.9.9}
       \begin{pmatrix}
        1&a&b&c\\ & 1&&b\\&&1&-a\\&&&1
    \end{pmatrix}\Sp(4,\bO)/\Sp(4,\bO)\cap \begin{pmatrix}
       t^{-m_2}&&&\\ *& t^{m_2-m_3}&&\\ *&*&t^{m_3-m_2}&\\ *&*&*&t^{m_{2}}
  \end{pmatrix}\Sp(4,\bO)/\Sp(4,\bO),
\end{equation}
where, $a, b, c\in \bF$ and the latter space is the $U_{\Sp(4)}^-(\bF)$-orbit of
$-m_2(\varepsilon_2-\varepsilon_3)-m_3\varepsilon_3$ in $\Gr_{\Sp(4)}$. We denote~\eqref{3.9.9} by
$\BY_{0,0}'$.

If the intersection $\BY_{0,0}'$ is non-empty, then $c\in t^{-m_2}\bO$ and $ab\in t^{-m_3}\bO+c$.
For any point in $\BY_{0,0}'$, we can find a unique representative such that
$a,b,c\in t^{-1}\BC[t^{-1}]$. We denote the corresponding point by $(a,b,c)$. This assignment gives
rise to an embedding 
\begin{equation}
    \begin{split}
        \BY_{0,0}'\hookrightarrow \bF/\bO\times \bF/\bO\times \bF/\bO.
    \end{split}
\end{equation}

The restriction of the shifted local system $\chi$ to $\BY_{0,0}'$ is the pullback of the exponential
$D$-module on $\BG_a$ along the map $\BY_{0,0}'\rightarrow \BG_a$ sending a point $(a,b,c)$ to the
residue of the sum of $a+b$.

Now, we are going to prove 
\begin{prop}
  Unless $m_2=1$ and $m_3=0$, the sub-bottom cohomology
  $H^1(\BY_{0,0}', \chi\overset{!}{\otimes} \Psi_{0,-\alpha}[-\dim \overset{\circ}{\CY}{}^{-\alpha}])$
  vanishes for any $m_2, m_3\geq 0$.
\end{prop}

\begin{proof}
We can assume the intersection $\BY_{0,0}'$ is non-empty, which implies $2m_2-m_3\geq 0$.

Case I) $2m_2-m_3\geq 3$.

We claim that in this case, the intersection $\BY_{0,0}'$ is $\BG_a$-invariant. Here, 
\begin{equation}\label{3.9.11}
    \BG_a:= 
       \left\{\begin{pmatrix}
        1&t^{-1}x&&\\ & 1&&\\&&1&-t^{-1}x\\&&&1
    \end{pmatrix},\ x\in \BC\right\},
\end{equation}
and it acts on $\Gr_{\Sp(4)}$ via the left multiplication.

It is enough to check that $\min\{\deg(t^{-1}x), \deg(t^{-1}xb)\}\geq -m_2+1$, for any $x\in \BC$ and
$b$ appearing in~\eqref{3.9.9}. Since we assume $2m_2\geq m_3+3\geq 3$, we have
$\deg(t^{-1}x)\geq -1\geq -m_2+1$. Furthermore, the minimal degree of the $2\times 2$-minors in the
first and second rows of the point $(a,b,c)\in\BY_{0,0}'$ is $-m_3$, so $\deg(b^2)=2\deg(b)\geq -m_3$.
The assumption implies that $\deg(b)\geq 2-m_2$. In particular, $\deg(t^{-1}xb)\geq -m_2+1$.

Now, the intersection $\BY_{0,0}'$ is isomorphic to $Y'\times \BG_a$, where $Y'$ is the subset of
$\BY_{0,0}'$ consisting of the points $(a,b,c)$ such that the residue of $a$ is~$0$. Since there
is no non-constant map $\BG_a\rightarrow \BG_m$, the shifted local system
$\chi\overset{!}{\otimes}\Psi_{0,-\alpha}$ on $\BY_{0,0}'$ is equivariant with respect to $\BG_a$
against the exponential $D$-module. As a result, its cohomology vanishes in all degrees.

\medskip

Case II) $2m_2-m_3=2$.
 
We decompose $\BY_{0,0}'$ into two pieces: $P_1$ and $P_2$. Here, $P_1$ consists of the points
$(a,b,c)$ such that $\deg(b)=-\frac{m_3}{2}=-m_2+1$; $P_2$ is the complement
(i.e.\ $\deg(b)>-\frac{m_3}{2}$). We study the cohomologies of $P_1$ and $P_2$.

Using the same analysis as in Case I), we can prove that $P_2$ is $\BG_a$-invariant with respect to the left multiplication of the $\BG_a$ in~\eqref{3.9.11}. Indeed, $(a,b,c)$ belongs to the intersection if and only if $\on{min}\{\deg(ab-c),0, \deg(b^2),\deg(b)\}=-m_3$ and $\on{min}\{0,\deg(a), \deg(b),\deg(c)\}=-m_2$. However, if $(a,b,c)\in P_2$, the action of $\BG_a$ does not change the corresponding degrees. Also, the shifted local system
$\chi\overset{!}{\otimes}\Psi_{0,-\alpha}$ is $\BG_a$-equivariant against the exponential $D$-module
and its cohomologies vanish.

Let us focus on $P_1$ and prove that the cohomologies on $P_1$ vanish. The condition $2m_2-m_3=2$ forces
$m_2\geq 1$. If $m_2=1$, we have $m_3=0$, this is the case which we do not need to consider
(the simple coroot case). So we assume $m_2\geq 2$.

We claim that for any point $(a,b,c)$, we have $a\in t^{1-m_2}\bO$. 

Indeed, note that the minimal degree of the entries of the first row is $-m_2$, so $a\in t^{-m_2}\bO$.
Also, the minimal degree of the $2\times 2$-minors in the first and second rows is $-m_3$, so
$\deg(ab-c)\geq -m_3=-2m_2+2$. If $\deg(a)=-m_2$, then $\deg(c)=\deg(ab)=-2m_2+1$. The non-empty
intersection condition implies $-2m_2+1\geq -m_2$, so $1\geq m_2$, which contradicts the assumption. 

Now we can describe $P_1$ explicitly. Namely, a point $(a,b,c)\in (\bF/\bO)^3$ belongs to $P_1$ if and only if $a\in t^{1-m_2}\bO/\bO$, $b\in t^{1-m_2}\bO^\times/\bO$, and $c\in t^{-m_2}\bO^\times/\bO$. That is to say
\begin{equation}\label{3.9.12}
    P_1\simeq \BG_a^{m_2-1}\times \BG_m\times \BG_a^{m_2-2}\times \BG_m\times \BG_a^{m_2-1}.
\end{equation}

The shifted local system $\Psi_{0,-\alpha}$ comes from a Kummer local system on $\BG_m\times \BG_m$.

If $m_3\geq 3$, the shifted local system $\chi$ on $P_1$ is the pullback of the exponential $D$-module
on $\BG_a\times \BG_a$ along the map
\begin{equation}
P_1 \longrightarrow     \BG_a^{m_2-1}\times \BG_a^{m_2-2}\longrightarrow \BG_a\times \BG_a.
\end{equation}
The first map is given by projecting to the first and the third factors of~\eqref{3.9.12}, and the
second map is given by projecting to the $\BG_a$-factors of $\BG_a^{m_2-1}$ and $\BG_a^{m_2-2}$ corresponding to the coefficients of $t^{-1}$. In particular, the shifted local system
$\chi\overset{!}{\otimes}\Psi_{0,-\alpha}$ is $\BG_a\times \BG_a$-equivariant against the exponential
$D$-module, so its cohomology vanishes in all degrees.

If $m_2=2$, then the shifted local system $\chi$ is the pullback of exponential $D$-module on $\BG_a\times \BG_a$ along the following map
\begin{equation}
    P_1\longrightarrow \BG_a\times \BG_m\hookrightarrow \BG_a\times \BG_a.
\end{equation}
Here, the first map is given by projecting to the first and second factors of~\eqref{3.9.12}
(note that $m_2-1=1$), and the second map is the open embedding $\BG_m\hookrightarrow\BG_a$.

So, the shifted local system  $\chi\overset{!}{\otimes}\Psi_{0,-\alpha}$ is $\BG_a$-equivariant
against the exponential $D$-module. Here, the action (different from the former action on $P_2$!) of $x\in \BG_a$ on $\BY_{0,0}'$ is given by
sending $(a,b,c)\mapsto (a+t^{-1}x, b,c)$. In particular, the cohomology of
$\chi\overset{!}{\otimes}\Psi_{0,-\alpha}$ vanishes.

\medskip

Case III)  $2m_2-m_3=1$.

For any point $(a,b,c)$ in the intersection, we have $\min\{\deg(ab-c), 0, \deg(b^2),\deg(b)\}= -m_3$.
Since $m_3$ is odd, $\deg(b^2)>-m_3$, and in particular, $\deg(b)>-m_3$. It implies that
$ab-c\in t^{-m_3}\bO^\times$. Also, note that $\min\{0, \deg(a), \deg(b), \deg(c)\}=-m_2$. 

If $m_2\geq 2$, we claim that  $\deg(a)=-m_2$. 

Otherwise, $\deg(c)=-m_2$ and $\deg(ab)>-m_3$. The latter implies $\deg(c)=-m_3$, so $m_2=m_3$, which
implies $m_2=1$, a contradiction.

So for any point in $\BY_{0,0}'$, we have $a\in t^{-m_2}\bO^\times$ and $c\in t^{-m_2}\bO$. As
$m_2\geq 2$, we have $m_3>m_2$. So $\deg(ab)= -m_3$ and $b\in t^{-m_2+1}\bO^\times$. On the other hand,
any $a,b,c$ satisfying the above conditions give rise to a point in $\BY_{0,0}'$. We have
\begin{equation}
    \BY_{0,0}'\simeq \BG_m\times \BG_a^{m_2-1}\times \BG_m\times \BG_a^{m_2-2}\times \BG_a^{m_2}.
\end{equation}

Consider the action of $\BG_a$ on $\BY_{0,0}'$ that sends $(a,b,c)$ to $(a+t^{-1}x, b,c)$. The $D$-module
$\chi\overset{!}{\otimes}\Psi_{0,-\alpha}$ is $\BG_a$-equivariant against the exponential $D$-module,
so its cohomology vanishes.

If $m_2=m_3=1$, we have $\BY_{0,0}'\simeq \BG_m\times \BG_a$, and $\chi$ is the pullback of the
exponential $D$-module on $\BG_a$ and $\Psi_{0,-\alpha}$ is the pullback of a Kummer local system on
$\BG_m$. The cohomology vanishes in all degrees as well.

\medskip

Case IV) $2m_2-m_3=0$.

Assume $m_3\neq 0$ (the case $m_3=m_2=0$ is trivial). Note that
$\min\{\deg(ab-c), 0, \deg(b^2),\deg(b)\}= -m_3$, so we must have $\deg(b)=-m_2$. So
$(a,b,c)\in \bF/\bO\times \bF/\bO\times \bF/\bO$ belongs to $\BY_{0,0}'$ if and only if
$\deg(a)\geq -m_2$, $\deg(b)=-m_2$, and $\deg(c)\geq -m_2$. So we have
\begin{equation}
    \BY_{0,0}'\simeq \BG_a^{m_2}\times \BG_m\times \BG_a^{m_2-1}\times \BG_a^{m_2}.
\end{equation}
Simlarly, the $D$-module $\chi\overset{!}{\otimes}\Psi_{0,-\alpha}$ is $\BG_a$-equivariant against the exponential
$D$-module, so its cohomologies vanish.
\end{proof}

Applying the same method as in \cite[\S 6.7]{g}, we obtain the following corollary
from~Theorem~\ref{thm 3.8.3}. We sketch the proof for the sake of completeness.

\begin{cor}\label{cor irr}
  The functor $F$ of~\eqref{functor F} takes irreducible generic Hecke equivariant $D$-modules
  to irreducible factorizable sheaves:
    \[F(\IC^\blambda_{\on{glob}})= \CF_\blambda.\]
\end{cor}

\begin{proof}
  Note that both of $F(\IC^\blambda_{\on{glob}})$ and $\CF_\blambda$ are supported on
  $C^\bullet_{\leq \blambda \cdot c}$. 
We prove by induction on the length of $\blambda-\bmu$ that their restrictions to
$C^{\bmu}_{\leq \blambda \cdot c}$ are isomorphic for any $\bmu\leq \blambda$.

For $\blambda-\bmu=0$, it just follows from $\sO^\blambda\cap \sS^\blambda=\pt$. Now, assume the
result for any $\blambda'-\bmu'< \blambda-\bmu$. According to Theorem~\ref{thm 3.8.3} and the factorization property of
$F(\IC^\blambda_{\on{glob}})$ and $\CF_\blambda$, we obtain that they are isomorphic on the complement of
the closed stratum $\{\bmu\cdot c\} {\hookrightarrow} C^{\bmu}_{\leq \blambda\cdot c}$.

Due to Lemma \ref{lemma 3.8.4}, to prove that $F(\IC_{\on{glob}}^\blambda)|_{C^{\bmu}_{\leq \blambda\cdot c}}$ is the IC-extension, we should prove that the $!$-stalk of $F(\IC_{\on{glob}}^\blambda)|_{\bmu\cdot c}$ is concentrated in degrees no less than $1$. In other words,
 \begin{equation}\label{eq 3.9.19-}
H^i(\CY_{\infty\cdot c}|_{\bmu \cdot c}, \fp^!(\IC^\blambda_{\on{glob}})|_{\bmu \cdot c}[-\on{rel.dim}(\Bun_B,\Bun_G)])=0,       \end{equation}
if $i\leq 0$ and $\bmu<\blambda$. 

Since $F$ is $t$-exact, the above cohomology vanishes if $i<0$, and the $0$-th cohomology equals the $!$-stalk of the image of $\nabla_{\on{glob}}^{\blambda}$, the $*$-extension of the rank 1 local system on $\scrM_{=\blambda\cdot c}$. The $!$-stalks of $\fp^!(\nabla_{\on{glob}}^{\blambda})[-\on{rel.dim}(\Bun_B,\Bun_G)]$ are concentrated in degrees $\dim \overset{\circ}{\CY}{}^{\bmu}_{=\blambda\cdot c}$.

Note that \begin{equation}\label{eq3.9.19}
\CY_{\infty\cdot c}|_{\bmu \cdot c}\simeq  \sS^{\bmu}.    
\end{equation}
According to Proposition \ref{prop 3.10.1},  $\dim \overset{\circ}{\CY}{}^{\bmu}_{=\blambda\cdot c}$ is no less
than twice the dimension of (any irreducible component of) $\sO^\blambda\cap \sS^\bmu$, so $0$-th cohomology is the bottom cohomology. Thus $H^0$ is non-zero only if $\fp^!(\IC^\blambda_{\on{glob}})|_{\bmu \cdot c}[-\on{rel.dim}(\Bun_B,\Bun_G)]$ is constant on some irreducible component of the critical dimension $\frac{1}{2}\cdot \dim \overset{\circ}{\CY}{}^{\bmu}_{=\blambda\cdot c}$. 


Under the above isomorphism, $\fp^!(\IC^\blambda_{\on{glob}})|_{\bmu \cdot c}$ is isomorphic to the tensor product of $\pi^!(\IC^\blambda_{\on{glob}})\subset \on{D-mod}_{\cP_{\det}^\kappa}( ^{\bomega}\Gr_{\oG})$ and the Kummer local system $\Psi_{\blambda, \bmu}$, up to a shift. Since the Kummer local system is tame and $\pi^!(\IC^\blambda_{\on{glob}})$ is the pullback of the exponential D-module, we only need to prove that $\Psi_{\blambda, \bmu}$ is non-constant. According to Lemma \ref{lem 3.7.3} and the assumption $\blambda>\bmu$, it is constant only if $\blambda-\bmu$ is a simple odd root. However, if $\blambda-\bmu$ is a simple odd root, $\frac{1}{2}\cdot \dim \oCY^{\bmu}_{=\blambda}$ is not an integer. In particular, there is no irreducible component of the critical dimension $\frac{1}{2}\cdot \dim \oCY^{\bmu}_{=\blambda}$, so \eqref{eq 3.9.19-} is zero for $i=0$.
\end{proof}

We sketch the proof of the following corollary, see \cite[Proposition 4.2.3, 4.2.4]{bft1} or \cite[Corollary 9.5.2]{ty} for details.
\begin{cor}\label{cor 4.1.1}
  \label{Fglob}
  The functor $F$ induces a conservative, faithful functor between abelian categories $F^{\on{glob}}\colon \cC_{\kappa}^{\on{glob},\lc,\heartsuit}\to \on{FS}$. 
\end{cor}
\begin{proof}
According to Lemma \ref{lem 3.8.3}, the functor $F$ restricts to $\cC_{\kappa}^{\on{glob}}$ is $t$-exact. Furthermore, we note that $\cC_{\kappa}^{\on{glob},\lc}$ is generated by irreducible objects $\IC^{\blambda}_{\on{glob}}$, and the image of $\IC^{\blambda}_{\on{glob}}$ under $F$ is $\cF_{\blambda}$ (see, Corollary \ref{cor irr}). Thus, the image of $\cC_{\kappa}^{\on{glob},\lc}$ satisfies the finiteness conditions in Section \ref{factor shv} as each $\cF_{\blambda}$ satisfies the finiteness conditions. We have a functor between abelian categories  $F^{\on{glob}}\colon \cC_{\kappa}^{\on{glob},\lc,\heartsuit}\to \on{FS}$.

The conservativity and faithfulness just follow from~Corollary~\ref{cor irr}.
\end{proof}

\section{Rigidity of $\cC_{\kappa}^{\on{glob},\on{lc}}$ and equivalence}
\label{Rigidity}

\subsection{Monoidal structure of $\cC_{\kappa}^{\on{glob}}$}
First let us recall the monoidal structure of the twisted global Gaiotto category. Recall that when we define the global Gaiotto category, we need to take $C=\BP^1$.

In \S \ref{global model}, we defined the global model $\scrM_{\infty\cdot c}$, for a given point $c\in C$. However, a similar definition applies to multiple points
$\underline{c}=\{c_1, c_2,\ldots, c_n\}\subset C$, the resulting stack is denoted by $\scrM_{\infty\cdot \underline{c}}$. 

Consider the global model over $\Ran_C$. To be more precise, $\scrM_{\infty\cdot \Ran}$ is defined as the ind-stack locally of finite type which classifies the data $(\underline{c}, \CP_{\oG}, \sigma)$, where $\underline{c}=\{c_1, c_2,\cdots, c_n\}\in C$, $\cP_{\oG}\in \Bun_{\oG}$, and $\sigma$ is a section $C\setminus \underline{c}\rightarrow \cX\overset{G}{\times}\CP_{\oG}$ which generically lands in $\overset{\circ}{\cX}\overset{G}{\times}\CP_{\oG}$. Its fiber over $\underline{c}\in \Ran_C$ is exactly  $\scrM_{\infty\cdot \underline{c}}$. We denote by $\scrM_{\infty\cdot C^n}$ the relative fiber product of $C^n$ and $\scrM_{\infty\cdot \Ran}$ over $\Ran$.

Using the same definition of the global Gaiotto category on $\scrM_{\infty\cdot c}$,
we can define the global Gaiotto categories of $\scrM_{\infty\cdot \Ran}$, $\scrM_{\infty\cdot \underline{c}}$ and $\scrM_{\infty\cdot C^n}$, which are denoted by $\cC_{\kappa,\Ran}^{\on{glob}}$, $\cC_{\kappa,\underline{c}}^{\on{glob}}$ and $\cC_{\kappa,C^n}^{\on{glob}}$, respectively. Then $\cC_{\kappa,\Ran}^{\on{glob}}$ is a category over $\Ran$, and its fiber over any point $\underline{c}\in \Ran_C$ is canonically equivalent to $\cC_{\kappa,\underline{c}}^{\on{glob}}$.

The map $\pi$ in~\S\ref{3.5.4} extends to a map over $\Ran_C$, i.e., we have a projection $ ^\bomega\Gr_{\overline{G}, \Ran}\longrightarrow \scrM_{\infty\cdot \Ran}$. By the same method for the proof
of~\cite[Theorem 4.4.1]{ty}, the $!$-pullback induces an equivalence of categories 
\begin{equation}\label{eq 4.1.1}
     \cC_{\kappa,\Ran}^{\on{glob}}\simeq  \on{D-mod}_{\CP^\kappa_{\det}}^{^\bomega\!H(\bF)_{\Ran},\chi}({}^\bomega\Gr_{\overline{G}, \Ran})=: \cC_{\kappa, \Ran}^{\loc}.
\end{equation}

Namely, for any $S$-points $\underline{c}\in \Ran_C$ and $S\in \Sch^{\aff}$, the same proof as in \textit{loc.cit.} implies that the $!$-pullback along $^\bomega\Gr_{\overline{G}, \underline{c}}\longrightarrow \scrM_{\infty\cdot \underline{c}}$ induces an equivalence $\cC_{\kappa,\underline{c}}^{\on{glob}}\simeq  \on{D-mod}_{\CP^\kappa_{\det}}^{^\bomega\!H(\bF)_{\underline{c}},\chi_{\underline{c}}}({}^\bomega\Gr_{\overline{G}, \underline{c}})$.

In particular, the factorization structure of $\cC_{\kappa, \Ran}^{\loc}$ gives a factorization category structure of $\cC_{\kappa,\Ran}^{\on{glob}}$: 
\begin{equation}\label{4.1.2}
     (\cC_{\kappa,C^2}^{\on{glob}})\underset{C^2}{\times} (C\times C)_{\on{disj}} \simeq (\cC_{\kappa,C}^{\on{glob}}\otimes \cC_{\kappa,C}^{\on{glob}})\underset{C\times C}{\times} (C\times C)_{\on{disj}}.
\end{equation}


\subsubsection{}
Due to the $\BG_a$-equivariance of $\cC_{\kappa,C}^{\on{glob}}|_{\BA^1}$, for any generic Hecke equivariant $D$-module $\cF\in \cC_{\kappa}^{\on{glob}}$, we can form an object in $\cC_{\kappa, C}^{\on{glob}}|_{\BA^1}$. 

Then, according to \eqref{4.1.2}, we can form an object $\cF_1\boxtimes\cF_2|_{\BA^2\setminus \Delta}$ on $ \scrM_{\infty\cdot C^2}|_{\BA^2\setminus \Delta}$ for any two generic Hecke equivariant $D$-modules $\cF_1,\cF_2\in \cC_{\kappa}^{\on{glob}}$. We can make an appropriate shift such that the functor 
\begin{equation}
\begin{split}
     \cC_{\kappa}^{\on{glob}}\otimes \cC_{\kappa}^{\on{glob}}\longrightarrow \cC_{\kappa, C^2}^{\on{glob}}|_{\BA^2\setminus \Delta}\\
     \cF_1, \cF_2\mapsto \cF_1\boxtimes\cF_2|_{\BA^2\setminus \Delta}
\end{split}
\end{equation}
is $t$-exact.

Now, we define the fusion $\cF_1\star \cF_2\in \cC_{\kappa}^{\on{glob}}$ to be
\begin{equation}\label{eq 4.13}
    \cF_1\star \cF_2\colon = \pr_{\Delta,*}\Psi_{x-y}(\cF_1\boxtimes\cF_2|_{\BA^2\setminus \Delta})[-1].
\end{equation}
Here, $\Psi_{x-y}$ is the nearby cycles functor along the diagonal, and $\pr_{\Delta}$ is the projection $\scrM_{\infty\cdot C}|_{\BA^1}\simeq \BA^1\times \scrM_{\infty\cdot \{0\}}\to \scrM_{\infty\cdot \{0\}}$. 

It is not hard to see that fusion preserves local compactness. 
\begin{lem}\label{preserve locally compactness}
    If $\cF_1, \cF_2\in \cC_{\kappa}^{\on{glob}, \lc}$, then $\cF_1\star \cF_2\in \cC_{\kappa}^{\on{glob}, \lc}$.
\end{lem}

\begin{proof}
  The finiteness of $!$-stalks of $\cF_1\star \cF_2$ follows from the definition. It is
  enough to check that $\supp(\IC_{\on{glob}}^{\blambda}\star \IC^{{\blambda'}}_{\on{glob}})$ only contains finitely many relevant strata of $\scrM_{\infty\cdot c}$, for any relevant $\blambda, \blambda'$. 

  Take an algebraic substack of $\scrM_{\infty\cdot C^2\setminus\Delta}$ locally of finite type, whose fiber  over any point $\underline{c}=(c_1, c_2)$ in  $C^2\setminus\Delta$ contains the closure of $\scrM_{=(\blambda, \blambda')\cdot \underline{c}}$, and consider its closure inside $\scrM_{\infty\cdot C^2}$. The fiber of the resulting stack over any point $c$
    in the diagonal contains finitely many relevant strata, which implies that $\supp(\IC_{\on{glob}}^{\blambda}\star \IC^{{\blambda'}}_{\on{glob}})$ only contains finitely many relevant strata of $\scrM_{\infty\cdot c}$. 
    \end{proof}

We will see later (ref., Corollary \ref{cor 4.1.8}) that the fusion gives a braided monoidal structure of $\cC_{\kappa}^{\on{glob}, \lc}$.

\subsubsection{Topological factorization category}\label{braided monoidal}
Consider the diagonal stratification of $C^n$, we denote by $\on{D}_{\on{constr}}(C^n)$ the category of constructible sheaves on $C^n$. For a stratum $C^{\frP}\subset C^n$, we denote by $\on{Loc}_{\on{rh}}(C^{\frP})$ the category of regular holonomic local systems on $C^{\frP}$. 

\begin{defn} (cf.~\cite[Corollary 3.2.2]{[Nick]})
  A topological factorization category is a collection of $\on{D}_{\on{constr}}(C^n)$-module categories $\cC_n$
  equipped with factorization structure that satisfies the following Ran condition:
\begin{itemize}
    \item For any surjection $\Delta_{I/J}\colon I\twoheadrightarrow J$ between two non-empty finite sets, the $!$-pullback induces an equivalence $\Delta_{I/J}^!\cC_{{|J|}}\simeq \cC_{{|I|}}$.
\end{itemize}    
\end{defn}
\begin{rem}
  This definition is different from the notion of topological factorization category
  in~\cite[Theorem 5.5.4.10]{[Lur]}, where the author considers the constructible sheaf of categories
  while we consider the module category over constructible sheaves. In particular, even when $C=\BA^1$,
  it is not obvious at all that one can get a braided monoidal category structure for the category
  $\cC_1$. However, see~Lemma~\ref{lin}.
\end{rem}

\subsubsection{}
We are going to first define a topological factorization category whose $\cC_1$ is given by $\cC_{\kappa}^{\on{glob},\lc}$, and then we prove that it satisfies the conditions in~Lemma~\ref{lin}.

According to the (Ran version) local-global equivalence, it is enough to define a topological factorization category for the local Gaiotto category. 

We define $\cC_n:=\cC_{\kappa, C^n,\on{rh}}^{\loc,\lc}$ as the full subcategory of $\cC_{\kappa, C^n}^{\loc,\lc}:=\on{D-mod}_{\CP^\kappa_{\det}}^{^\bomega\!H(\bF)_{C^n},\chi}({}^\bomega\Gr_{\overline{G}, C^n})$ consisting of the objects which are locally constant along each stratum $C^{\frP}\subset C^n$. To be more precise, an object belongs to $\cC_{\kappa, C^n,\on{rh}}^{\loc,\lc}$ if and only if its $!$-restriction to $C^{\frP}\subset C^n$ belongs to $\on{Loc}_{\on{rh}}(C^{\frP})\otimes \cC_{\kappa, \underline{c}}^{\loc,\lc}$, for any stratum $C^{\frP}$. According to the Riemann-Hilbert correspondence, the category $\cC_{\kappa, C^n,\on{rh}}^{\loc,\lc}$ is a category with an action of $\on{D}_{\on{constr}}(C^n)$. Furthermore, $\cC_{\kappa, C^n,\on{rh}}^{\loc,\lc}$ satisfies the Ran condition and acquires the factorization structure. 

An important feature  of $\cC_{\kappa, C^n,\on{rh}}^{\loc,\lc}$ is that six operations between categories for
various strata (closures) of $C^n$ are well-defined.
\begin{prop}\label{exist !}
    For $C^{\frP}\subset C^n$, the $!$-pushforward functor (i.e., the left adjoint functor of $!$-restriction) $\cC_{\kappa, C^{\frP},\on{rh}}^{\loc,\lc}\rightarrow \cC_{\kappa, C^n,\on{rh}}^{\loc,\lc}$ is well-defined. Similarly, the $*$-pullback functor is also well-defined.
\end{prop}
\begin{proof}
We only prove the first claim.    We denote $\tilde{\cC}_{\kappa, C^n,\on{rh}}^{\loc,\lc}$ the full subcategory of compact objects in $\on{D-mod}_{\CP^\kappa_{\det}}^{ ^\bomega\! \on{U}(\bF)_{C^n},\chi}({}^\bomega\Gr_{\PSp(6), C^n})$ which are locally constant along each stratum $C^{\frP}\subset C^n$. So $\cC_{\kappa, C^n,\on{rh}}^{\loc,\lc}$ is the $\PGL(2, \bO)_{C^n}$-invariants of $\tilde{\cC}_{\kappa, C^n,\on{rh}}^{\loc,\lc}$. 

First, according to the locally constant assumption of $\tilde{\cC}_{\kappa, C^n,\on{rh}}^{\loc,\lc}  $, the $!$-pushforward functor 
\begin{equation}\label{exist ! 1}
\tilde{\cC}_{\kappa, C^{\frP},\on{rh}}^{\loc,\lc}\rightarrow \tilde{\cC}_{\kappa, C^n,\on{rh}}^{\loc,\lc}    
\end{equation}
 is well-defined.

The functor \eqref{exist ! 1} induces a functor $\cC_{\kappa, C^{\frP},\on{rh}}^{\loc,\lc}\rightarrow \cC_{\kappa, C^n,\on{rh}}^{\loc,\lc}$, such that the following diagram commutes
\[\xymatrix{
\cC_{\kappa, C^{\frP},\on{rh}}^{\loc,\lc}\ar[r]\ar[d] &\cC_{\kappa, C^n,\on{rh}}^{\loc,\lc}\ar[d]\\
\tilde{\cC}_{\kappa, C^{\frP},\on{rh}}^{\loc,\lc}\ar[r]& \tilde{\cC}_{\kappa, C^n,\on{rh}}^{\loc,\lc}.
}\]

Here, the vertical arrows are forgetful functors

Since $\PGL(2,\bO)_{C^n}$ is formally smooth pro-finite group over $C^n$, the corresponding functor $\cC_{\kappa, C^{\frP},\on{rh}}^{\loc,\lc}\rightarrow \cC_{\kappa, C^n,\on{rh}}^{\loc,\lc}$ is the left adjoint functor of the $!$-restriction $\cC_{\kappa, C^n,\on{rh}}^{\loc,\lc}\rightarrow \cC_{\kappa, C^{\frP},\on{rh}}^{\loc,\lc}$. 
\end{proof}


\begin{rem}
  One can also interpret our categories and the $!$-pushfoward and $*$-pullback functors entirely
  in the constructible setting. For example, the $!$-pushforward functor in \eqref{exist ! 1} can
  be defined as follows by the Kirillov model.

 According to the equivalence between the Whittaker model and Kirillov model, cf.~\cite[Appendix A]{gl}, we can replace $\tilde{\cC}_{\kappa, C^n}^{\loc,\lc}$ with its corresponding Kirillov model. The latter is defined as follows: denote by $\tilde{\on{U}}$ the kernel of $\chi\colon ^\bomega\! \on{U}(\bF)_{C^n}\rightarrow \BG_a$, then the category  $\on{D-mod}_{\CP^\kappa_{\det}}(\tilde{\on{U}}\backslash {}^\bomega\Gr_{\PSp(6), C^n}):= \on{D-mod}^{\tilde{\on{U}}}_{\CP^\kappa_{\det}}( {}^\bomega\Gr_{\PSp(6), C^n})$ admits an action of $\BG_m\ltimes \BG_a$, the Kirillov model is defined as the full subcategory of compact objects belonging to the kernel of $\Av_*^{\BG_a}$ in
$\on{D-mod}_{\CP^\kappa_{\det}}^{\BG_m}(\tilde{\on{U}}\backslash {}^\bomega\Gr_{\PSp(6), C^n})$. It is equivalent to $\tilde{\cC}_{\kappa, C^n}^{\loc,\lc}$, and makes sense in the constructible sheaf setting. In the constructible setting, it is equivalent to $\tilde{\cC}_{\kappa, C^n, \on{rh}}^{\loc,\lc}$ due to the Riemann-Hilbert correspondence and the constructibility with respect to the stratification of $C^n$.

The $!$-pushforward functor
\begin{equation}\label{4.1.8}
    \on{D}_{\on{constr},\CP^\kappa_{\det}}({}^\bomega\Gr_{\PSp(6), C^{\frP}})\rightarrow  \on{D}_{\on{constr},\CP^\kappa_{\det}}({}^\bomega\Gr_{\PSp(6), C^n})
\end{equation}
 induces the $!$-pushforward in the equivariant setting
\begin{equation}\label{4.1.7}
    \on{D}_{\on{constr},\CP^\kappa_{\det}}(\tilde{\on{U}}\backslash {}^\bomega\Gr_{\PSp(6), C^{\frP}})\rightarrow  \on{D}_{\on{constr},\CP^\kappa_{\det}}(\tilde{\on{U}}\backslash {}^\bomega\Gr_{\PSp(6), C^n}).
\end{equation}

Now, the corresponding functor of \eqref{exist ! 1} in the constructible setting for the Kirillov model is induced from \eqref{4.1.7}.
\end{rem}

\begin{rem}\label{rem 4.1.6}
By local-global comparison, we also have a topological factorization category whose corresponding $\on{D}_{\on{constr}}(C^n)$-module category is $\cC_{\kappa, C^n, \on{rh}}^{\on{glob},\lc}$.    
\end{rem}

\subsubsection{Construction of the global functor over $C^n$}
Similarly to the functor $F$ in \eqref{functor F}, we can also construct a factorizable functor between categories over $C^n$ (or, over the Ran space). Namely, we consider the SW-zastava space with $n$ marked points $\CY_{\infty\cdot C^n}:= \scrM_{\infty\cdot C^n}\underset{\Bun_{\oG}}{\times'} \Bun_B$ which has the projection $\fv_{C^n}$ to the configuration space with $n$ marked points $C^\bullet_{\infty\cdot C^n}$ and the projection $\fp_{C^n}$ to $\scrM_{\infty\cdot C^n}$. Here, $\scrM_{\infty\cdot C^n}\underset{\Bun_{\oG}}{\times'} \Bun_B$ denotes the open subspace of $\scrM_{\infty\cdot C^n}\underset{\Bun_{\oG}}{\times} \Bun_B$, such that the generic $H$-reduction and $B$-reduction are generically transverse. Then, we define the functor
\begin{equation} \label{functor F Cn}
  \on{D-mod}_{\CP^\kappa_{\det}}(\scrM_{\infty\cdot C^n})\longrightarrow
  \on{D-mod}_{\CP^\kappa}(C^\bullet_{\infty\cdot C^n}),\quad
\CF\mapsto\fv_{C^n,*}(\fp_{C^n}^\bullet(\CF)),
    \end{equation}
where similarly to~\eqref{functor F}, $\fp^\bullet_{C^n}(\CF):=\fp^!_{C^n}(\CF)[-\on{rel.dim} (\Bun_B,\Bun_{\overline{G}})]$ with parity changed according to
$\on{rel.dim} (\Bun_B,\Bun_{\overline{G}})\pmod{2}$.

Similarly to Lemma \ref{lemma 3.8.4}, \eqref{functor F Cn} induces a $t$-exact functor $F_{C^n}^{\on{glob}}\colon \cC_{\kappa, C^n}^{\on{glob},\lc}\rightarrow D^b(\FS_n)$. Furthermore, the resulting functor is factorizable. Indeed, by the same constructions and arguments as in \cite[\S 7.2, 7.3]{ty}, we can construct a factorizable functor from $\on{D-mod}_{\CP^\kappa_{\det}}^{^\bomega\!H(\bF)_{\Ran},\chi,\lc}({}^\bomega\Gr_{\overline{G}, C^n})$ to $D^b(\FS_n)$, which is compatible with the functor $F_{C^n}^{\on{glob}}$ with respect to the local-global equivalence functor \eqref{eq 4.1.1}. Furthermore, it restricts to a factorizable functor between topological factorization categories, $F^{\on{glob}}_{C^n, \on{rh}}\colon \cC_{\kappa, C^n, \on{rh}}^{\on{glob},\lc}\rightarrow D^b_{\on{rh}}(\FS_n)$, where $D^b_{\on{rh}}(\FS_n)$ denotes the bounded derived category of factorization modules with $n$-marked points, which is also constructible with respect to the diagonal stratification of $C^n$ and satisfies the finiteness conditions in \cite{bfs}.

\subsubsection{Braided monoidal category structure}
We learned the following facts about topological factorization category from Lin Chen.

\begin{lem}\label{lin}
  For a topological factorization category $\CalC$ over $\BA^1$, if we require the existence of the
  $*$-pullback functor from $\CalC_n$ to $\CalC_{\frP}:=\CalC_n\otimes_{{\on{D}}_{\on{constr}}(C^n)}\on{Loc}_{\on{rh}}(C^\frP)$,
  then the nearby cycles functor (along the diagonal $\BA^1\subset \BA^1\times \BA^1$) provides the fiber of the factorization category over one point with a lax braided monoidal structure. 

Furthermore, a factorizable functor between topological factorization categories preserving $*$-pushforward functor gives rise to a lax braided monoidal functor among the corresponding braided monoidal categories. 
\end{lem}

As a direct corollary, $\cC_{\kappa}^{\on{glob,lc}}$ is a lax braided monoidal category, and if the factorizable functor is conservative and the target lax braided monoidal category is strict, then the source category is also braided monoidal and the functor is braided monoidal.

\subsubsection{Comparison of lax monoidal structures}
    In this section, we want to show that the lax braided monoidal functor defined via Lemma \ref{lin} coincides with the fusion product in \eqref{eq 4.13}. Both of them are given by taking certain nearby cycles, but in different categores. Namely, the former nearby cycles functor is in the equivariant setting, which by definition uses the composition of $*$-pushforward and $*$-pullback functor (or the $!$-pushforward functor and $!$-pullback) for the equivariant D-module category, i.e., the one defined in Proposition \ref{exist !}; while the definition of \eqref{eq 4.13} uses the nearby cycles functor in the plain category (without equivariant structure) of ind-holonomic D-modules, which can be rewritten using the corresponding pushforward and pullback functors for the plain category of D-modules. Since $*$-pushforward and $!$-pullback functors in the equivariant setting and non-equivariant setting coincide, we need to check that the $*$-pullback functor (or, the $!$-pushforward functor) defined in Proposition \ref{exist !} and Remark \ref{rem 4.1.6} coincides with the $*$-pullback functor (resp.\ the $!$-pushforward functor) for the plain category of ind-holonomic D-modules on the global model. To prove that the lax braided monoidal functor defined via Lemma \ref{lin} coincides with the functor in \eqref{eq 4.13}, we should prove that the following diagram commutes
  \[ \xymatrix{
     \cC^{\on{glob},\lc}_{\kappa, C^2\setminus \Delta,\on{rh}}|_{\BA^2\setminus \Delta}\ar[r]\ar[d] & \cC^{\on{glob},\lc}_{\kappa, C^2,\on{rh}}|_{\BA^2}\ar[d]\\
    \on{D-mod}_{\cP_{\det}^{\kappa},\on{h}} (\scrM_{\infty\cdot C^2\setminus \Delta})|_{\BA^2\setminus \Delta}\ar[r] &\on{D-mod}_{\cP_{\det}^{\kappa},\on{h}}{(\scrM_{\infty\cdot C^2}|_{\BA^2}),
    }}\]
  where the horizontal arrows are $!$-pushforwards, vertical arrows are forgetful functors,
  and the subscript `h' means ind-holonomicity. The problem of $*$-pullback is dealt similarly.

Let us denote by  $\scrM_{\underline{x},\infty\cdot C^2}|_{\BA^2}$ the open substack of  $\scrM_{\infty\cdot C^2}|_{\BA^2}\underset{\BA^2}{\times} (\BA\setminus \underline{x})^2$, such that the Pl\"{u}cker maps are injective maps of vector bundles near $\underline{x}$.

The global Gaiotto category $\cC^{\on{glob},\lc}_{\kappa, C^2,\on{rh}}|_{\BA^2}$ is equivalent to the
limit of categories $\cC^{\on{glob},\lc}_{\kappa, \underline{x}, C^2,\on{rh}}|_{\BA^2}$~\cite[Section 5.4]{ty}.
The latter category is the category of Hecke equivariant (only at $\underline{x}$)
$\cP_{\det}^{\kappa}$-twisted D-modules on the open substack $\scrM_{\underline{x},\infty\cdot C^2}|_{\BA^2}$
with locally compact and constructible conditions. Thus, it is enough to show that the following
diagram commutes.
 \[ \xymatrix{
     \cC^{\on{glob},\lc}_{\kappa, \underline{x}, C^2\setminus \Delta,\on{rh}}|_{\BA^2\setminus \Delta}\ar[r]^{j_!^{eq}}\ar[d]^{\on{oblv}} & \cC^{\on{glob},\lc}_{\kappa,  \underline{x}, C^2,\on{rh}}|_{\BA^2}\ar[d]^{\on{oblv}}\\
    \on{D-mod}_{\cP_{\det}^{\kappa},\on{h}} (\scrM_{\underline{x},\infty\cdot C^2\setminus \Delta}|_{\BA^2\setminus \Delta})\ar[r]^{j_!} &\on{D-mod}_{\cP_{\det}^{\kappa},\on{h}}{(\scrM_{\underline{x},\infty\cdot C^2}|_{\BA^2}).
    }}\]
Here, $j_!^{eq}$ is the $!$-pushforward in the equivaraint category setting.
    
We should check that the natural transform 
\begin{equation}\label{natural trans}
    j_!\circ \on{oblv}\rightarrow \on{oblv}\circ j_!^{eq}
\end{equation} is an isomorphism. 

Note that for any point in $\scrM_{\underline{x},\infty\cdot C^2}$, we have a well-defined
$^\bomega\!H$-bundle on the formal neighborhood of $\underline{x}$ (see~\S\ref{3.5.2} for notation).
By adding an isomorphism of this $^\bomega\!H$-bundle and the one induced from $\bomega^\rho$,
we obtain an $^{\bomega}H(\bO)_{\underline{x}}$-torsor on $\scrM_{\underline{x},\infty\cdot C^2}$, denoted by
$\tilde{\scrM}_{\underline{x},\infty\cdot C^2}$. We consider the $!$-pullback from
$\scrM_{\underline{x},\infty\cdot C^2}|_{\BA^2}$ to $\tilde{\scrM}_{\underline{x},\infty\cdot C^2}|_{\BA^2}$. 

Since both $^{\bomega}H(\bF)_{\underline{x}}$ and $^{\bomega}H(\bO)_{\underline{x}}$ are formally smooth, the natural transforms of the following outer square and the lower square are isomorphisms.
 \[ \xymatrix{
     \cC^{\on{glob},\lc}_{\kappa, \underline{x}, C^2\setminus \Delta,\on{rh}}|_{\BA^2\setminus \Delta}\ar[r]^{j_!^{eq}}\ar[d]^{\on{oblv}} & \cC^{\on{glob},\lc}_{\kappa,  \underline{x}, C^2,\on{rh}}|_{\BA^2}\ar[d]^{\on{oblv}}\\
    \on{D-mod}_{\cP_{\det}^{\kappa},\on{h}} (\scrM_{\underline{x},\infty\cdot C^2\setminus \Delta}|_{\BA^2\setminus \Delta})\ar[r]^{j_!}\ar[d]^{\on{!-pullback}} &\on{D-mod}_{\cP_{\det}^{\kappa},\on{h}}{(\scrM_{\underline{x},\infty\cdot C^2}|_{\BA^2})}\ar[d]^{\on{!-pullback}}\\
    \on{D-mod}_{\cP_{\det}^{\kappa},\on{h}} (\tilde{\scrM}_{\underline{x},\infty\cdot C^2\setminus \Delta}|_{\BA^2\setminus \Delta})\ar[r]^{\tilde{j}_!} &\on{D-mod}_{\cP_{\det}^{\kappa},\on{h}}{(\tilde{\scrM}_{\underline{x},\infty\cdot C^2}|_{\BA^2})}
    }\]

 Now, the desired property following from the fact that $!$-pullback functor to the bundle is
 conservative, which implies the natural transform~\eqref{natural trans} of the upper square is
 also an isomorphism.

Thus, the two $!$-pushforward functors  match, one is defined in the twisted Gaiotto category,
and the other is in the plain category of twisted D-modules on the global model. For a similar
reason, the $*$-pullback functor shares the same property. In particular, the lax braided monoidal
functor in Lemma \ref{lin} indeed coincides with the fusion constructed in \eqref{eq 4.13}.

\subsubsection{}
Note that we have
\begin{lem}
    The functor $F^{\on{glob}}_{C^n,\on{rh}}$ is compatible with $*$-pullback and $!$-pushforward.
\end{lem}
\begin{proof}
    If $\bmu$ is antidominant enough,
  $\Bun_B^\bmu\rightarrow \Bun_G$ is smooth, and it follows that $*$-pullback and $!$-pushforward are compatible with
  $\fp^!$ on the corresponding components. Then by the factorization property of $\fp^!(\cF)$ for
  any $\cF\in \cC_{\kappa}^{\on{glob}}$, one can extend the compatibility to the whole zastava space.
  Furthermore, since  $\fv_*$ can be written as the composition of a clean extension and a proper
  pushforward, and $*$-pullback and $!$-pushforward commute with proper pushforward, it follows that $F^{\on{glob}}_{C^n,\on{rh}}$ preserves the $*$-pullback and $!$-pushforward.
\end{proof}
  
  Thus, since the category of representations of the quantum supergroup is braided monoidal, we have
  \begin{cor}\label{cor 4.1.8}
 $\cC_{\kappa}^{\on{glob},\lc}$ acquires a braided monoidal category structure $(\cC, \star, {b})$ and $F$ is braided monoidal.      
  \end{cor}
  \begin{proof}
    It remains to check that  $F^{\on{glob}}_{C^n,\on{rh}}$ is conservative, which follows
    immediately from the stratawise conservativity. The latter follows
    from~Corollary~\ref{cor 4.1.1} and the definition of $F^{\on{glob}}_{C^n,\on{rh}}|_{C^{\frP}}$. 
  \end{proof}

\subsection{Rigidity of $\cC^{\on{glob}}$}
Now, let us temporarily restrict our attention to the untwisted case, i.e., $\kappa=0$. The goal of this section is to prove
\begin{prop}\label{unt rigid}
    The untwisted Gaiotto category $\cC^{\on{glob},\lc}$ is rigid.
\end{prop}
\subsubsection{Satake action on $\cC^{\on{glob}}$}
The Satake category is defined as the category of $\oG(\bO)$-equivariant $D$-modules (with coefficients in
super vector spaces) on $\Gr_{\oG}$. Using the local-global comparison~\cite[Theorem 4.4.1]{ty}, we can realize it as a global category, i.e., the category of generic Hecke equivariant $D$-modules on the global Hecke stack $\Hecke^{\on{glob}}_{\oG,c}$, we denote this category by $\mathcal{H}^{\on{glob}}_{\oG}$. Here the global Hecke stack is the algebraic stack that classifies the data $(\cP_{\oG}, \tilde{\cP}_{\oG}, \alpha)$, where $\cP_{\oG}, \tilde{\cP}_{\oG}\in \Bun_{\oG}$, and $\alpha$ is an isomorphism of ${\cP}_{\oG}$ and $\tilde{\cP}_{\oG}$ on $C\setminus c$.

Denote by ${\scrM}^{\on{conv}}_{\infty\cdot c}$ the convolution stack that classifies the data $({\cP}_{\oG}, \tilde{\cP}_{\oG}, \sigma, \alpha)$, where $({\cP}_{\oG},\sigma)$ is a point in $\scrM_{\infty\cdot c}$ and $({\cP}_{\oG},\tilde{\cP}_{\oG},\alpha)$ is a point in $\Hecke^{\on{glob}}_{\oG,c}$.

Consider the following diagram
\begin{equation}
\xymatrix{
&{\scrM}^{\on{conv}}_{\infty\cdot c}\ar[ld]^{p_1}\ar[d]^{p_2}\ar[rd]^{p_3}&\\
\scrM_{\infty\cdot c}& \scrM_{\infty\cdot c}& \Hecke_{\oG, c}^{\on{glob}}.
}    
\end{equation}
Here $({\cP}_{\oG}, \tilde{\cP}_{\oG}, \sigma, \alpha)\in {\scrM}^{\on{conv}}_{\infty\cdot c}$
goes to $({\cP}_{\oG},\sigma)$ under the morphism $p_1$, and goes to
$(\tilde{\cP}_{\oG}, \tilde{\sigma}\colon C\setminus c\overset{\sigma}{\rightarrow} \cX\overset{\oG}{\times}\cP_{\oG}|_{C\setminus c}\overset{\alpha}{\simeq} \cX\overset{\oG}{\times}\tilde{\cP}_{\oG}|_{C\setminus c})$
under the morphism $p_2$, and goes to $({\cP}_{\oG},\tilde{\cP}_{\oG},\alpha)$ under the morphism $p_3$. 

Now, we define the convolution action of the global Satake category on the global Gaiotto category as
\begin{equation}
    \cF *  \CS:= p_{2,*} (p_1^!(\cF)\overset{!}{\otimes} p_3^!(\CS)),
\end{equation}
for $\CS\in \mathcal{H}^{\on{glob}}_{\oG}$ and $\cF\in \cC^{\on{glob}}$.

The following lemma immediately follows from the definition:
\begin{lem}\label{lem 4.2.2}
    The action of $\on{D-mod}^{{}^\bomega\! \oG(\bO)}({}^\bomega\Gr_{\oG})$ on $\cC^{\loc}=\on{D-mod}{}^{^\bomega\! H(\bF),\chi}({}^\bomega\! \Gr_{\oG})$ is compatible with the action of $\mathcal{H}^{\on{glob}}_{\oG}$ on $\cC^{\on{glob}}$.
\end{lem}

We claim that that the fusion product of the global Gaiotto category is compatible with
the convolution action of the global Satake category:
\begin{prop}\label{conv comp}
    For any $\cF_1, \cF_2\in \cC^{\on{glob}}$, and $\CS_1, \CS_2\in \mathcal{H}^{\on{glob}}_{\overline{G}}$, we have a canonical isomorphism
    \begin{equation}
        (\cF_1 * \CS_1)\star (\cF_2 *  \CS_2)\simeq (\cF_1\star \cF_2) *  (\CS_1\star \CS_2).
    \end{equation}
    Here, $\CS_1\star \CS_2$ denotes the convolution product (equivalently, fusion product) of $\CS_1$ and $\CS_2$.
\end{prop}
\begin{proof}
  Let us consider the algebraic stack $\scrM_{\infty C^2}^{\on{conv}}$ that classifies the data
  $((c_1, c_2)\in C^2, \cP_{\oG}, \tilde{\cP}_{\oG}, \sigma, \alpha)$, that is defined
  similarly to $\scrM_{\infty C^2}^{\on{conv}}$, but here we require that $\alpha$ and $\sigma$
  are defined on $C\setminus \{c_1, c_2\}$. Also, we have an algebraic stack
  $\Hecke_{\oG, C^2}^{\on{glob}}$.

Consider the following modified convolution diagram 
\begin{equation}
\xymatrix{
&{\scrM}^{\on{conv}}_{\infty\cdot C^2}\ar[ld]^{p_{1, C^2}}\ar[d]^{p_{2, C^2}}\ar[rd]^{p_{3,C^2}}&\\
\scrM_{\infty\cdot C^2}& \scrM_{\infty\cdot C^2}& \Hecke_{\oG, C^2}^{\on{glob}}.
}    
\end{equation}

We $!$-pullback $\cF_1\boxtimes \cF_2|_{\scrM_{\infty C^2}|_{\BA^2\setminus \Delta}}$ along $p_{1, C^2}$ to ${\scrM}^{\on{conv}}_{\infty\cdot C^2}$, and $!$-tensor with the $!$-pullback of $\CS_1\boxtimes \CS_2|_{\Hecke_{\oG, C^2}^{\on{glob}}|_{\BA^2\setminus \Delta}}$ along $p_{3, C^2}$. The resulting sheaf is denoted by $\cF_1\tilde{\boxtimes}\cF_2\tilde{\boxtimes}\CS_1\tilde{\boxtimes} \CS_2$.

Note that the map $p_{2,C^2}$ is ind-proper, so pushforward along $p_{2,C^2}$ commutes with taking nearby cycles.

On the one hand, the nearby cycles of $\cF_1\tilde{\boxtimes}\cF_2\tilde{\boxtimes}\CS_1\tilde{\boxtimes} \CS_2$ along $\scrM_{\infty\cdot C^2}|_{\Delta=\BA^1}$ is constant along $\BA^1$ with the stalk $(\cF_1\star \cF_2) \tilde{\boxtimes} (\CS_1\star \CS_2)$, and its pushforward to $\scrM_{\infty\cdot C^2}|_{\Delta}$ along $p_{2, C^2}$ is constant along $\BA^1$ with the stalk $(\cF_1\star \cF_2)  *  (\CS_1\star \CS_2)$.

On the other hand, the pushforward of $\cF_1\tilde{\boxtimes}\cF_2\tilde{\boxtimes}\CS_1\tilde{\boxtimes} \CS_2$ along $p_{2, C^2}$ is $(\cF_1 *  \CS_1)\tilde{\boxtimes} (\cF_2 *  \CS_2)$, and its nearby cycles along $\scrM_{\infty\cdot C^2}|_{\Delta}$ is constant along $\BA^1$ with the stalk $(\cF_1 *  \CS_1)\star (\cF_2 * \CS_2)$.

So the proposition is proved.
\end{proof}

\begin{prop}\label{generate unt}
    The untwisted Gaiotto category $\cC^{\on{glob},\lc}$ is generated by the action of $\Hecke^{\on{glob},\lc}$ on the unit object $\IC^0_{\on{glob}}$.
\end{prop}
\begin{proof}
  For convenience, we fix some notation. We call a coweight $\blambda$ of $\oG$ dominant
  if it is anti-dominant for $\PGL(2)$ part and is dominant for $\PSp(6)$ part. Also,
  we set $H'=g^{-1}H g$.

Note that according to the decomposition theorem, for any dominant $\blambda$, the convolution $\IC^0_{\on{glob}} *  \IC^{\blambda, \Hecke}$ is semi-simple. So, we should prove that $\IC^0_{\on{glob}} *  \IC^{\blambda, \Hecke}= \IC^{\blambda}_{\on{glob}}\bigoplus_{
\blambda'<\blambda} (\IC^{\blambda'}_{\on{glob}})^{\oplus n_\blambda'}$.

    According to Lemma~\ref{lem 4.2.2}, it is enough to do the calculation for the local action of the locally compact Satake category on $\on{D-mod}^{^\bomega\!H(\bF),\chi,\lc}(^\bomega\Gr_{\oG})$. We can ignore the $\bomega$-renormalization, and study the action of $\on{D-mod}^{\oG(\bO),\lc}(\Gr_{\oG})$ on the unit object in $\on{D-mod}^{H(\bF),\chi,\lc}(\Gr_{\oG})$.

We should prove that 
\begin{equation}\label{4.2.5}
    H'(\bF)\oG(\bO)\overset{\oG(\bO)}{\times} \overline{\Gr}_{\oG}^\blambda\longrightarrow \Gr_{\oG}
\end{equation}
is contained in the closure of the $H'(\bF)$-orbit passing through $t^\blambda$, and its fiber over $t^\blambda$ is exactly one point.

For the first claim, we first note that the closure of $H'(\bF)\oG(\bO)$ equals the closure of $H'(\bF)I$, where $I:=I_{\PGL(2)}^+\times I_{\PSp(6)}^-$. Furthermore, for any dominant $\blambda$, we always have $H'(\bF)I t^\blambda I= H(\bF)t^\blambda I$. In particular, we obtain that $H'(\bF)\oG(\bO)t^\blambda \oG(\bO)$ is contained in the closure of $H'(\bF)t^\blambda \oG(\bO)$.

For the second claim, note that the fiber of $t^\blambda\in G(\bF)$ along~\eqref{4.2.5} is the intersection of the unital $H'(\bF)$-orbit and the closure of the unital $t^{\blambda} \oG(\bO) t^{-\blambda}$-orbit in $\Gr_{\oG}$. The latter is contained in the closure of the unital $U_{\PGL(2)}(\bF)\times U_{\PSp(6)}(\bF)^-$-orbit in $\Gr_{\oG}$, which intersects with the unital $H'(\bF)$-orbit just at one point $1\in \Gr_{\oG}$ by~Proposition~\ref{so intersect}. 
\end{proof}

With the help of the above proposition, we can finish the proof of~Proposition~\ref{unt rigid}.
\begin{proof}(of Proposition~\ref{unt rigid})
The unit object $\IC^0_{\on{glob}}$ is rigid in $\cC^{\on{glob},\lc}$. Indeed, it is the unital object for the fusion monoidal structure.

To see this, we notice that there is a closed embedding of algebraic stacks
\begin{equation}\label{4.2.6}
    \scrM_{\infty\cdot C}\times C\hookrightarrow \scrM_{\infty\cdot C^2},
\end{equation}
which sends a point $(c, \CP_{\oG}, \sigma), \tilde{c}$ to $((c,\tilde{c})\in C^2, \CP_{\oG}, \sigma)$.

For any $\CF\in \cC^{\on{glob}}$, the sheaf $\cF\boxtimes \IC^0_{\on{glob}}|_{\scrM_{\infty C^2}|_{\BA^2\setminus \Delta}}$ is supported on the closed substack \eqref{4.2.6}. Furthermore, its restriction to \eqref{4.2.6} is constant along $\BA^2\setminus \Delta$ with fiber $\CF$, so its nearby cycles along $\BA^1$ is constant along $\BA^1$ with fiber $\cF$. It proves that $\CF\star \IC^0_{\on{glob}}\simeq \cF$. Similarly, one can prove that $\IC^0_{\on{glob}}\star \CF\simeq \CF$.

Now, the proposition follows from Propositions~\ref{conv comp},~\ref{generate unt}.
\end{proof}

\subsection{Rigidity of $\cC_{\kappa}^{\on{glob},\lc}$}
In order to prove that the twisted Gaiotto category is equivalent to the category of representations of the quantum supergroup, we should prove that the generic twisted Gaiotto category $\cC_{\kappa}^{\on{glob},\lc}$ is rigid. However, since the twisted Satake category is too small (when $\kappa$ is generic, it only contains one irreducible object), $\cC_{\kappa}^{\on{glob},\lc}$ is not generated by the action of the twisted Satake category on the unit object $\IC^0_{\on{glob}}\in \cC_{\kappa}^{\on{glob},\lc}$. We should use a different approach.

\subsubsection{}
Let $\cC^{\on{glob},\lc}(\CP_{\det})$ denote the category of locally compact generic Hecke equivariant
$D$-modules on the punctured determinant line bundle $\CP_{\det}$ on $\scrM_{\infty\cdot c}$. Since
$\BG_m$ acts on the punctured determinant line bundle, there is a strong action of $\BG_m$ on
$\cC^{\on{glob},\lc}(\CP_{\det})$, i.e., $\cC^{\on{glob},\lc}(\CP_{\det})$ is a module category over the
monoidal category $\on{D-mod}^c(\BG_m)$ (here we use $D$-modules with coefficients in vector spaces, not in
super vector spaces). The Fourier-Mellin transform provides an equivalence
$\on{D-mod}^c(\BG_m)\simeq \Perf(\BA^1/\BZ)$. Here $\Perf(\BA^1/\BZ)$ is the category of
$\BZ$-equivariant quasi-coherent sheaves on $\BA^1$ finitely generated over the ring generated by
the functions on $\BA^1$ and the translations by $\BZ$.
For example, $\oplus_{n\in \BZ} \delta_{n}\in \Perf(\BA^1/\BZ)$.

Now, for any $n\in \BN$, we define the category $\cC^{\on{glob},\lc}(\CP_{\det})_n$ as
\begin{equation}
    \cC^{\on{glob},\lc}(\CP_{\det})_n:= \cC^{\on{glob},\lc}(\CP_{\det})\otimes_{\on{D-mod}^c(\BG_m)} \Perf(\BA^1/\BZ)_{n}. 
\end{equation}
Here $\Perf(\BA^1/\BZ)_{n}$ is the category of perfect complexes on the $n$-th thickening of
$0\in \BA^1/\BZ$, it is equivalent to the category of perfect complexes on the $n$-th thickening of
$0\in \BA^1$, the functor $\on{D-mod}^c(\BG_m)\rightarrow  \Perf(\BA^1/\BZ)_{n}$ is given by taking restriction $\on{D-mod}^c(\BG_m)\simeq  \Perf(\BA^1/\BZ)\rightarrow \Perf(\BA^1/\BZ)_{n}$.

The fusion product of $\cC_{\kappa}^{\on{glob},\lc}$ naturally extends to a fusion product of $\cC^{\on{glob},\lc}(\CP_{\det})$ and $\cC^{\on{glob},\lc}(\CP_{\det})_n$. Furthermore, since the image of $\cC^{\on{glob},\lc}$ under the compsoition $\cC^{\on{glob},\lc}\rightarrow \cC^{\on{glob},\lc}(\CP_{\det})\rightarrow \cC^{\on{glob},\lc}(\CP_{\det})_n$ generates $\cC^{\on{glob},\lc}(\CP_{\det})_n$, each $ \cC^{\on{glob},\lc}(\CP_{\det})_n$ is rigid as well.

As a corollary, the limit of  $\cC^{\on{glob},\lc}(\CP_{\det})_n$ with respect to the restriction transition functors is also rigid, since taking limit preserves rigidity if transition functors are monoidal.

\subsubsection{}
Let us denote by $ \cC^{\on{glob},\lc}(\CP_{\det})_{\hat{0}}$ the relative tensor product 
\begin{equation}
     \cC^{\on{glob},\lc}(\CP_{\det})\otimes_{\on{D-mod}^c(\BG_m)} \Perf(\Spec(\BC[\![\kappa]\!])).
\end{equation}

The key observation is that 
\begin{prop}
  \label{global limit}
    The restriction functor induces an equivalence of categories
    \begin{equation}\label{4.3.3}
        \cC^{\on{glob},\lc}(\CP_{\det})_{\hat{0}}\simeq \lim_n \cC^{\on{glob},\lc}(\CP_{\det})_{n}.
    \end{equation}
\end{prop}

\begin{proof}
    First, we prove that~\eqref{4.3.3} is fully faithful.

    For any two objects $\cF_1,\CF_2\in  \cC^{\on{glob},\lc}(\CP_{\det})$,
    $\Hom(\cF_1\otimes \BC[\![\kappa]\!],\CF_2\otimes \BC[\![\kappa]\!])$ is finitely generated over
    $\BC[\![\kappa]\!]$, so as  $ \cC^{\on{glob},\lc}(\CP_{\det})_{\hat{0}}$ is generated by objects of the
    form $\cF\otimes \BC[\![\kappa]\!]$, we obtain that for any two objects in
    $\cC^{\on{glob},\lc}(\CP_{\det})_{\hat{0}}$, their Hom is finitely generated over $\BC[\![\kappa]\!]$.
    Now the fully faithfulness follows from the fact that for any finitely generated
    $\BC[\![\kappa]\!]$-module $M$, the map
    \begin{equation}
        M\longrightarrow R\lim_{n} M\otimes_{\BC[\![\kappa]\!]} \BC[\kappa]/\kappa^n
    \end{equation}
is an isomorphism. Here, ${R\lim}$ denotes \textit{derived limit}, that is the limit in the derived category of complexes.
    
   It remains to prove that \eqref{4.3.3} is essentially surjective.

   We note that both sides of \eqref{4.3.3} admit filtrations given by relevant orbits, and the
   functor in~\eqref{4.3.3} is compatible with respect to the orbits-filtration. So we should prove
   that for any $(\PGL(2,\bO)\ltimes U(\bF),\chi)$-relevant orbit $\BO^\blambda$ in $\Gr_{\PSp(6)}$
   (equivalently, any locally closed substack of $\scrM_{\infty\cdot c}$ which admits non-zero generic
   Hecke equivariant $D$-modules), the restriction of \eqref{4.3.3} to $\BO^\blambda$ is an equivalence,
   i.e.,
   \begin{equation}\label{4.3.5}
   \begin{split}
       \on{D-mod}^{\PGL(\bO)\ltimes U(\bF),\chi,\lc}(\CP_{\det}|_{\BO^\blambda})\otimes_{\on{D-mod}(\BG_m)^c} \Perf(\Spec \BC[\![\kappa]\!])\\
       \simeq\\
       \lim_n  \on{D-mod}^{\PGL(\bO)\ltimes U(\bF),\chi,\lc}(\CP_{\det}|_{\BO^\blambda})\otimes_{\on{D-mod}(\BG_m)^c} \Perf(\Spec \BC[\kappa]/\kappa^n).
   \end{split}
   \end{equation}

We divide the question into two cases:

I). Assume that the stabilizer in $\PGL(\bO)\ltimes U(\bF)$ of a point $L_{\blambda}$ in
$\BO^{\blambda}\subset \Gr_{\PSp(6)}$ (equivalently, of any point in $\BO^{\blambda}$) acts on the fiber
$\CP_{\det}|_{L_{\blambda}}$ non-trivially. Actually, in this case, the stabilizer must be $\BG_m$, and it
acts on the fiber via $x\mapsto x^m$, for a non-zero integer $m$.
   
Using the Riemann-Hilbert correspondence, we can consider the corresponding constructible sheaf category. The LHS of~\eqref{4.3.5} is isomorphic to 
\begin{equation}
\begin{split}
   & \on{D-mod}^{\BG_m,\lc}(\BG_m)\otimes_{\on{D-mod}^c(\BG_m)} \Perf(\Spec \BC[\![\kappa]\!])   \\
    \simeq& \BC[q,q^{-1}]/q^m-1 \otimes^L_{\BC[q,q^{-1}]} \BC[\![q-1]\!]\mPerf\\
    \simeq &\BC\mPerf.
\end{split}
\end{equation}

The RHS is isomorphic to
\begin{equation}\label{4.3.7}
\begin{split}
  \lim_{n} \BC[q,q^{-1}]/q^m-1 \otimes^L_{\BC[q,q^{-1}]} \BC[q-1]/(q-1)^m\mPerf.
\end{split}
\end{equation}

For any $m\geq 1$, we have an isomorphism of dg-algebras
\begin{equation}
     \BC[q,q^{-1}]/q^m-1 \otimes^L_{\BC[q,q^{-1}]} \BC[q-1]/(q-1)^m \simeq \BC[\epsilon].
\end{equation}
Here, $\epsilon$ is a generator in the degree $-1$.

Using the above isomorphism, we have
\begin{equation}
    \eqref{4.3.7}\simeq \lim_n \BC[\epsilon]\mPerf
\end{equation}
here, the transition map in the RHS is given by sending $\epsilon$ to $0$. In particular, the transition map between $\BC[\epsilon]\mPerf$ factors through $\BC\mPerf$. So, $\lim_n \BC[\epsilon]\mPerf\simeq \BC\mPerf$.

\medskip

II). Assume that the stabilizer in $\PGL(\bO)\ltimes U(\bF)$ of a point $L_{\blambda}$ in
$\BO^{\blambda}\subset \Gr_{\PSp(6)}$ (equivalently, of any point in $\BO^{\blambda}$) acts on the fiber
$\CP_{\det}|_{L_{\blambda}}$ trivially.

In this case, let us denote by $S$ the stabilizer group of $L_{\blambda}\in \Gr_{\PSp(6)}$.
The LHS of~\eqref{4.3.5} is equivalent to 
\begin{equation}
    H_*(S)\modu^{f.d}\otimes \BC[\![\kappa]\!]\mPerf,
\end{equation}
and the RHS is equivalent to 
\begin{equation}
    \lim_n  H_*(S)\modu^{f.d}\otimes \BC[\kappa]/\kappa^n\mPerf,
\end{equation}
where $ H_*(S)\modu^{f.d}$ denotes the category of complexes over $H_*(S)$ whose underlying complex is bounded and each cohomology is finite dimensional.

Now, the desired equivalence follows from the following lemma.
\end{proof}
\begin{lem}
    If $A$ is finitely generated non-negatively graded dg-algebra over $\BC$, then we have the following equivalence of categories
    \begin{equation}\label{4.3.11}
        A\modu^{f.d}\otimes \BC[\![\kappa]\!]\mPerf\simeq \lim_n A\modu^{f.d}\otimes \BC[\kappa]/\kappa^n\mPerf.
    \end{equation}
\end{lem}
\begin{proof}
  Under our assumption, $ A\modu^{f.d}\otimes \BC[\![\kappa]\!]\mPerf$ can be identified with the
  category of bounded complexes over the dg-algebra $A\otimes \BC[\![\kappa]\!]$, such that each
  cohomology is finitely generated over $\BC[\![\kappa]\!]$. A similar description holds for $A\modu^{f.d}\otimes \BC[\kappa]/\kappa^n$ for any $n$.

    Note that $\lim_n A\modu\otimes \BC[\kappa]/\kappa^n\modu$ admits a fully faithful functor into the category of complexes over $A\otimes \BC[\![\kappa]\!]$, and the image is identified with the full subcategory of bounded complete complexes over $A\otimes \BC[\![\kappa]\!]$. Under this identification, the right-hand-side of~\eqref{4.3.11} is identified with the category of bounded complete complexes over $A\otimes \BC[\![\kappa]\!]$ that are finite dimensional modulo $\kappa$.

    It remains to prove that any such complex is finitely generated over $\BC[\![\kappa]\!]$. Using the filtration given by truncation functors, we only need to prove that the bottom cohomology is finitely generated. However, it directly follows from the Nakayama lemma. 
\end{proof}

Now from Proposition~\ref{global limit} we obtain
\begin{cor}
    $\cC^{\on{glob},\lc}(\cP_{\det})_{\hat{0}}$ is rigid.
\end{cor}

The category 
$ \cC^{\on{glob},\lc}(\CP_{\det})\otimes_{\on{D-mod}^c(\BG_m)} \Perf(\Spec(\BC(\!(\kappa)\!)))$ is also rigid as
a localization.
Since $\BC(\kappa)\subset \BC(\!(\kappa)\!)$ is a sub-field, we see that the category $\cC^{\on{glob},\lc}(\CP_{\det})\otimes_{\on{D-mod}^c(\BG_m)} \Perf(\Spec(\BC(\kappa)))$ is also rigid.

Now we can finish the proof of the main result of this section.
\begin{prop}
  The $\kappa$-twisted Gaiotto category
  $\cC_{\kappa}^{\on{glob},\on{lc}}$ (equivalently,$\cC_{\kappa}^{\on{loc},\on{lc}}$) is rigid for a Weil-generic $\kappa$.
\end{prop}
\begin{proof}
  First, if we consider $D$-modules with coefficients in $\BQ$ and we regard the affine
  Grassmannian and its determinant line bundle as schemes defined over $\BQ$, then we can also
  define the category of locally compact Gaiotto equivariant $D$-modules on the punctured
  determinant line bundle. The resulting category is denoted by $\cC_{\BQ}^{\loc,\lc}(\CP_{\det})$. By taking relative tensor product, we obtain
    $\cC_{\BQ(\kappa)}^{\loc,\lc}:= \cC_{\BQ}^{\loc,\lc}(\CP_{\det})\otimes_{\Perf_{\BQ}(\BA^1/\BZ)} \BQ(\kappa)\mPerf$.

  We claim that $\cC_{\BQ(\kappa)}^{\loc,\lc}$ is rigid. It is enough to show that it is rigid after
  the extension of field of scalars. However, its scalar extension by
  $\BQ(\kappa)\subset \BC(\!(\kappa)\!)$ is the category $\cC^{\on{loc},\lc}(\CP_{\det})\otimes_{\on{D-mod}^c(\BG_m)} \Perf(\Spec(\BC(\kappa)))\simeq \cC^{\on{glob},\lc}(\CP_{\det})\otimes_{\on{D-mod}^c(\BG_m)} \Perf(\Spec(\BC(\kappa)))$, which is rigid.

    Take any transcendental number $\kappa_0$, and consider the map $\BQ(\kappa)\rightarrow \BC$ given by $\kappa=\kappa_0$. We obtain that
$\cC_{\kappa_0}^{\on{loc},\lc}= \cC_{\BQ(\kappa)}^{\loc,\lc}\otimes_{\BQ(\kappa)\mPerf} \BC\mPerf$ is rigid.
\end{proof}

\subsection{Equivalence}\label{equiv}
Now we are able to prove that the functor
$F^{\on{glob}}\colon \cC_{\kappa}^{\on{glob},\lc,\heartsuit}\to \on{FS}$
of~Corollary~\ref{Fglob} is an equivalence for transcendental $\kappa$. The proof is the same
as in~\cite[\S\S4.3--4.5]{bft1} or~\cite[\S10]{ty} with one minor modification.
Recall that a weight $\blambda\in Q$ satisfying condition~\eqref{odd dominant} is called
{\em typical} if the scalar product $(\blambda+\rho,\alpha)\ne0$ for odd roots $\alpha$.
For example, $\blambda_0=(\frac72,\frac12,\frac12,\frac12)$ is typical. For a typical $\blambda$,
the corresponding irreducible module $V_\blambda\in\Rep_q(\on{F}(4)^{\on{ad}})$ forms a block by
itself, i.e.\ for any weight ${\blambda'}\ne\blambda$ in $Q$ satisfying condition~\eqref{odd dominant},
we have $\Ext^\bullet(V_\blambda,V_{\blambda'})=0=\Ext^\bullet(V_{\blambda'},V_\blambda)$. Moreover, if $M_\blambda,M_{\blambda'}$
denote the corresponding Verma modules over $U_q(\ff(4))$, then
$\Ext^\bullet(M_\blambda,V_{\blambda'})=0=\Ext^\bullet(M_{\blambda'},V_\blambda)$.

Now the argument of~\cite[\S5]{bft1} proves that for a relevant orbit $\sO^{\blambda'}\subset\Gr_{\overline{G}}$,
the costalk of $\IC_{\loc}^{\blambda'}$ at $\sO^\blambda$ vanishes (under the assumption that $\blambda$
is typical, and ${\blambda'}\ne\blambda$). The same argument proves that the costalk of
$\IC_\loc^\blambda$ at $\sO^{\blambda'}$ vanishes under the same assumption.
Then the argument of~\cite[Corollary 4.4.2]{bft1} proves that $\IC_\loc^\blambda$ is a projective-injective object of the abelian category
$\cC_{\kappa}^{\on{loc},\lc,\heartsuit}$, where the $t$-structure of $\cC_{\kappa}^{\on{loc},\lc}$ is given by the pullback of the naive $t$-structure of $\cC_{\kappa}^{\on{glob},\lc}$. Thus, by~Theorem~\ref{thm 3.8.3} we conclude that $\IC_{\on{glob}}^\blambda$ a projective-injective object of the abelian category $\cC_{\kappa}^{\on{glob},\lc,\heartsuit}$. Due to the rigidity property of
$\cC_{\kappa}^{\on{glob},\lc,\heartsuit}$, the argument of~\cite[Corollary 4.4.3]{bft1} goes through, and
we obtain enough projective-injective objects of $\cC_{\kappa}^{\on{glob},\lc,\heartsuit}$.
The rest of the argument in~\cite[\S4.5]{bft1} goes through as well, and we derive our main
results.

\begin{thm}
  \label{derived}
  The natural functor $D^b(\cC_{\kappa}^{\on{glob},\lc,\heartsuit})\to\cC_{\kappa}^{\on{glob},\on{lc}}$ is an
  equivalence.
\end{thm}

\begin{thm}
  \label{gaiotto}
  The functor $F^{\on{glob}}\colon \cC_{\kappa}^{\on{glob},\lc,\heartsuit}\to \on{FS}$ is a braied monoidal
  equivalence.
\end{thm}
\begin{proof}
  For the full faithfulness, we note that $\IC_{\on{glob}}^{\blambda_0}\star \IC_{\on{glob}}^{\blambda}$ forms a collection of projective-injective generators of $\cC_{\kappa}^{\on{glob},\lc,\heartsuit}$, it is enough to check that the functor $F^{\on{glob}}$ induces an isomorphism $\Ext^i_{\cC_{\kappa}^{\on{glob},\lc,\heartsuit}}(\IC_{\on{glob}}^{\blambda_0}\star \IC_{\on{glob}}^{\blambda}, \IC_{\on{glob}}^{\blambda_0}\star \IC_{\on{glob}}^{\blambda'})=\Ext^i_{\on{FS}}(F(\IC_{\on{glob}}^{\blambda_0})\star F(\IC_{\on{glob}}^{\blambda}), F(\IC_{\on{glob}}^{\blambda_0})\star F(\IC_{\on{glob}}^{\blambda'}))=\Ext^i_{\on{FS}}(\cF_{\blambda_0}\star \cF_{\blambda}, \cF_{\blambda_0}\star \cF_{\blambda'})$. Due to the rigidity of irreducible objects, it is enough to show that $\Ext^i_{\cC_{\kappa}^{\on{glob},\lc,\heartsuit}}(\IC_{\on{glob}}^{\blambda_0}, \cF)= \Ext^i_{\on{FS}}(\cF_{\lambda_0}, F(\cF))$, for any $\cF\in \cC_{\kappa}^{\on{glob},\lc,\heartsuit}$. However, since $\IC_{\on{glob}}^{\blambda_0}$ and $\cF_{\lambda_0}$ are projective, both of $\Ext^i$ vanish if $i\neq 0$, and both of $\Ext^0$ equal the multiplicities of $\IC_{\blambda_0}^{\on{glob}}$ in $\cF$.

  The essential surjectivity follows from the fact that $\on{FS}$ is generated by $\cF_\blambda$ for all $\blambda$ relevant.

\end{proof}

\begin{thm}
      The functor $F_{C^n}^{\on{glob}}: \cC_{\kappa, C^n}^{\on{glob},\lc}\rightarrow D^b(\FS_n)$ is an equivalence.
  \end{thm}
  
\begin{proof}
  It follows by the same induction and gluing arguments as in \cite[\S 8]{g}. To be self-contained, we review the arguments here. 
  
  Namely, by induction, we can assume these two categories are equivalent on the diagonal divisor $\on{Diag}\subset C^n$. Also, note that both of $\cC_{\kappa, C^1}^{\on{glob}, \heartsuit, \lc}$ and $\on{FS}_1$ are constant families of categories on $C^1$ with the fiber categories $\cC_{\kappa}^{\on{glob}, \heartsuit, \lc}$ and $\on{FS}$, so the functor $F_{C^1}^{\on{glob}}$ is an equivalence of categories on $C^1$. Then, according to the factorization category structures of $\cC_{\kappa, C^n}^{\on{glob}, \heartsuit, \lc}$ and $\on{FS}_n$, we obtain that these two categories are equivalent on $\oC^n$. 

  To glue the equivalences to an equivalence on $C^n$, we note that for any object $\cF\in \cC_{\kappa, \oC^n}^{\on{glob}, \heartsuit, \lc}$ and $\cF'\in \cC_{\kappa, \on{Diag}}^{\on{glob}, \heartsuit, \lc}$, we can find a sequence of objects $\cF_i\in \cC_{\kappa, C^n}^{\on{glob}, \heartsuit, \lc}$, such that
  the restriction of $\cF_i$ to the open subspace $\scrM_{\infty\cdot \oC^n}\subset \scrM_{\infty\cdot C^n}$ is $\cF$, and both of $\lim_i\Ext^{i}(\cF_i, \mathfrak{i}_{\on{Diag\to C^n,*}}(\cF'))$ and $\lim_i\Ext^{i}(F_{C^n}^{\on{glob}}(\cF_i), F_{C^n}^{\on{glob}}(\mathfrak{i}_{\on{Diag\to C^n,*}}(\cF')))$ vanish, for $i=0,1$. Here, we can choose $\cF_i$ to be the objects such that the inverse limit of $\cF_i$ goes to the $!$-extension of $\cF$. 
  
  The vanishing of $\Ext^i$ implies that we can calculate $\Ext^i_{\cC_{\kappa, C^n}^{\on{glob}, \heartsuit, \lc}}$ of any two objects in $\cC_{\kappa, C^n}^{\on{glob}, \heartsuit, \lc}$ by a limit and extensions of $\Ext^i$ in $\cC_{\kappa, \on{Diag}}^{\on{glob}, \heartsuit, \lc}$ and $\cC_{\kappa, \oC^n}^{\on{glob}, \heartsuit, \lc}$, similarly for $\Ext^i_{\FS_n}$. Now, the stratawise equivalence implies the desired equivalence over $C^n$.
\end{proof}


\begin{cor}
  The composition of the derived equivalence of~Theorem~\ref{bfsl}:
  $D^b(\Rep_q(\on{F}(4)^{\on{ad}}))\iso D^b(\on{FS})$, with the quasiinverse
  of the derived equivalence of~Theorem~\ref{gaiotto} and the equivalence of~Theorem~\ref{derived}
  gives rise to a braided monoidal equivalence $D^b(\Rep_q(\on{F}(4)^{\on{ad}}))\iso\cC_{\kappa}^{\on{glob},\on{lc}}$. Its further composition with $\pi^!$ gives a braided monoidal equivalence between $D^b(\Rep_q(\on{F}(4)^{\on{ad}}))$ and $(\cC_{\kappa}^{\loc,\lc}, \star, b)$.
\end{cor}


\end{document}

The paper password for this article is: m7nj4